\documentclass[10pt]{article}

\usepackage[full]{textcomp}
\usepackage[osf]{newtxtext}
\usepackage{cabin}
\usepackage[varqu,varl]{inconsolata}
\usepackage[bigdelims,vvarbb]{newtxmath}
\usepackage[cal=boondoxo]{mathalfa}
\usepackage{subfiles}

\usepackage{yfonts}
\usepackage{enumitem}

\usepackage{tikz}
\usepackage{mathrsfs}

\usepackage{mhequ}
\usepackage{mhenvs} 				

\usepackage{diagrams1}
\usepackage{MalliavinDiagrams}
\usepackage{shortcuts}

\usepackage{fancyhdr}
\usepackage{breakurl}
\usepackage{comment}

\usepackage{wasysym}

\usepackage[normalem]{ulem}  		

\newtheorem{assumption}{Assumption}
\newtheorem{example}{Example}			

\usepackage{hyperref}

\begin{document}

\title{Malliavin Calculus and Density for Singular Stochastic Partial Differential Equations
}

\author{Philipp Sch\"onbauer}


\maketitle

\begin{abstract}
We study Malliavin differentiability of solutions to sub-critical singular parabolic stochastic partial differential equations (SPDEs) and we prove the existence of densities for a class of singular SPDEs. Both of these results are implemented in the setting of regularity structures. For this we construct renormalized models in situations where some of the driving noises are replaced by deterministic Cameron-Martin functions, and we show Lipschitz continuity of these models with respect to the Cameron-Martin norm. In particular, in many interesting situations we obtain a convergence and stability result for lifts of $L^2$-functions to models, which is of independent interest. The proof also involves two separate algebraic extensions of the regularity structure which are carried out in rather large generality.
\end{abstract}

\section{Introduction}

We establish Malliavin differentiability and subsequently study the existence of densities of solution to singular stochastic partial differential equations (SPDEs). The equations we have in mind are formally given by systems of the form
\begin{align}\label{eq:singular:SPDE}
\partial_t u_i = \CL_i u_i
+
F_i(u, \nabla u, \ldots)
+
\sum_{j\le n} F_i^j(u, \nabla u, \ldots)\xi_j \,,
\qquad{i \le m}
\end{align}
where each component $u_i$ is in general a distribution on $\R\times \T^d$ for some $d \ge 1$, subject to some initial condition $u_i(0,\cdot)=u_{i,0}$. Here, $\CL_i$  is an elliptic differential operator involving only spatial derivatives, the functions $F_i$ and $F_i^j$ are smooth and allowed to depend on $u=(u_i)_{i\le m}$ and finitely many derivatives of $u$, and 
the random fields $\xi_j$, $j\le m$, are assumed to be jointly Gaussian. 

Equations of type (\ref{eq:singular:SPDE}) have been subject to intensive study in recent years and lead to the development of novel technical approaches  \cite{Hairer2014,GubinelliImkellerPerkowski2014,OttoWeber2016}. While these approaches differ in their scope and technical details, in situations where more then one of them can be applied, they lead to the same notion of solution. For the purpose of this paper we focus on the theory of regularity structures, originally developed in \cite{Hairer2014}, and subsequently extended and generalized in a series of papers \cite{BrunedHairerZambotti2016,ChandraHairer2016,BrunedChandraChevyrecHairer2017}.
Interesting examples that fall under this setting include the generalized KPZ equation \cite{Hairer2016,Ferrari2013,Spohn2013}
\begin{align}\label{eq:example:kpz}
\partial_t u_i = \partial_x^2 u_i + 
\sum_{k,l\le m} f_{k,l}^i(u) (\partial_x u_k) (\partial_x u_l)
+
\sum_{k \le n} g_k(u) \xi_k
\end{align}
in $1+1$ dimensions, the $\Phi^p_d$ \cite{Hairer2014,Weber2016,BrunedChandraChevyrecHairer2017,Tsatsoulis2018}equations
\begin{align}\label{eq:example:phi}
\partial_t u = \Delta u + \sum_{k\le p} a_k u^k + \xi
\end{align}
in $1+d$ dimensions for $d\le 3$ and the generalized PAM equation \cite{Hairer2014,GubinelliImkellerPerkowski2014}
\begin{align}\label{eq:example:pam}
\partial_t u = \Delta u + f(u) + \sum_{i,j \le d} f_{i,j}(u) (\partial_i u)(\partial_j u)
+
g(u)\xi.
\end{align}
in $d=2$ or $d=3$ dimensions. Choosing $\xi$ as white noise, which is the natural choice in these examples, all of these equations have in common that there does not exist a solution in the classical sense. The robust solution theory of  \cite{Hairer2014,BrunedHairerZambotti2016,ChandraHairer2016,BrunedChandraChevyrecHairer2017} instead considers approximate, \emph{renormalized} equations that take the form
\begin{align}\label{eq:singular:SPDE:reg}
\partial_t u^\eps_{i} = \CL_i u^\eps_i
+
F_i(u^\eps, \nabla u^\eps, \ldots)
+
\sum_{j\le n} F_i^j(u^\eps, \nabla u^\eps, \ldots)\xi^\eps_j
+
\sum_{k \le K} c^\eps_k \Upsilon_i^k(u^\eps, \nabla u^\eps, \ldots),
\end{align}
subject to some initial condition $u^\eps_{i}(0,\cdot)=u^\eps_{i,0}$, where $\xi^\eps_j = \xi_j * \rho^{(\eps)}$ for some approximate $\delta$-distribution $\rho^{(\eps)}$. In \cite{BrunedChandraChevyrecHairer2017} it was shown that under some appropriate assumption on the equation there exists a choice of constants $c^\eps_k$ with the property that the sequence of solutions $u^\eps$ converges in probability to some limiting random distribution $u$ as $\eps \to 0$, and we call this limit $u$ the (renormalized) solution to (\ref{eq:singular:SPDE}). 

\begin{remark}
The counter-terms $\Upsilon^k_i$ and the renormalization constants $c^\eps_k$ are given explicitly in \cite{BrunedChandraChevyrecHairer2017} but do not matter much at this stage. We recall their definition in (\ref{eq:counter-term}) below.
\end{remark}

The first purpose of the present article is to establish the existence of continuous path-wise derivatives of the renormalized solution to (\ref{eq:singular:SPDE}) in the direction of Cameron-Martin functions (in the sense of \cite[Def.~3.3.1]{UstunelZakai2000}). This is in particular enough to obtain the existence of a localized version of Malliavin derivative (\cite[Prop.~4.1.3]{Nualart2006}, \cite[Prop.~2.4]{CannizzaroFrizGassiat2015}), which in turn is sufficient for the celebrated Bouleau-Hirsch criterion \cite{BouleauHirsch1986} to apply. The latter gives rather sharp conditions under which densities with respect to Lebesgue measure exist. 

The second purpose of this article is to show that the conditions of the Bouleau-Hirsch criterion are indeed satisfied for an interesting class of equations. The equations for which we can show existence of densities include in particular the stochastic heat equation with multiplicative noise and the $\Phi^p_d$-equations.

The strategy outlined above was already used in \cite{CannizzaroFrizGassiat2015} to show existence of densities for the 2D-PAM equation, and a recent paper \cite{GassiatLabbe2017} treated the case of the $\Phi^4_3$ equation. On the technical level, our approach for showing Malliavin differentiability uses extensions of the regularity structure and is strongly inspired by \cite{CannizzaroFrizGassiat2015}, although the proofs given in the present paper differ in some key aspects, which in particular allows us to obtain statements that are more general. In the second part of the paper we apply the Bouleau-Hirsch criterion by studying the "dual" to the tangent equation, an idea that was already used in \cite{GassiatLabbe2017} to study the existence of densities for~$\Phi^4_3$.

To make things more concrete, we put ourselves in the setting of the "black box" theorem \cite[Thm.~2.13]{BrunedChandraChevyrecHairer2017}. Given non-linearities $F_i$ and $F_i^j$, and a Gaussian noise $\xi$, this theorem establishes explicit formulae for the counterterms and renormalization constants appearing in \ref{eq:singular:SPDE:reg}, and works out concrete assumption on the equations under which the sequence of renormalized solutions converge. 

\begin{assumption}\label{ass:main}
We assume throughout this paper that the assumptions of \cite[Thm.~2.13]{BrunedChandraChevyrecHairer2017} on the equation, the noises and the initial condition are satisfied. To be more precise, we assume \cite[(2.4),Ass.~2.4, Ass.~2.6, Ass.~2.10, (2.18)]{BrunedChandraChevyrecHairer2017},  we  assume that we are given jointly Gaussian random fields $(\xi_j)_{j \le n}$ in the sense of \cite[Def.~2.11]{BrunedChandraChevyrecHairer2017} and we assume that the initial condition can be decomposed as $u_0^\eps = \CS_{\rho,\eps}^-(\xi)(0,\cdot)+\psi^\eps$ as in \cite[(2.17)]{BrunedChandraChevyrecHairer2017} with $\psi^\eps$ converging to some random initial condition $\psi$ in probability in $\CC^\ireg$ as $\eps \to 0$. We refer the reader to Section \ref{sec:singularSPDEs} for a summary of these assumptions and the definition of the space $\CC^\ireg$.
\end{assumption}

For the reader not familiar with these assumptions we recall briefly their purpose. In \cite[(2.4)]{BrunedChandraChevyrecHairer2017} the authors give a rigours meaning to the notion of \emph{sub-criticality}. This is a key assumption which is seen in any of the theories developed in \cite{Hairer2014,GubinelliImkellerPerkowski2014,OttoWeber2016}, and the equations are believed to behave quite differently when this assumption is violated. It also ensures that one can algebraically build a regularity structure adapted to the equation as in \cite{BrunedHairerZambotti2016}. Assumption \cite[Ass.~2.4]{BrunedChandraChevyrecHairer2017} deals with compositions of the solution with smooth functions. It also limits the regularity blow-up at the initial time-slice to ensure that the solution is an actual distribution on the whole space (as opposed to just $\R^+ \times \T^d$). Throughout the solution theory developed in \cite{Hairer2014} the equations are treated in their mild formulation, and Assumption \cite[Ass.~2.6]{BrunedChandraChevyrecHairer2017} guarantees the existence of a Green's function for $\partial_t - \CL_i$, together with suitable analytic estimates. Assumption \cite[(2.10)]{BrunedChandraChevyrecHairer2017} is a technical assumption that ensures that the solution to our equation can always be written as an \emph{explicit} distribution-valued, stationary, random process, plus an \emph{implicit} function-valued random perturbation (by \emph{explicit} we mean that this process is given as a stationary solution to a \emph{linear} equation and polynomial expressions in this solution). The explicit stationary process appearing for the \emph{regularized} noise is denoted by $\CS_{\rho,\eps}^-(\xi)$ and appears in the rather cumbersome way in which the initial conditions are phrased. This is needed, since in general the spaces in which $\CS_{\rho,\eps}^-(\xi)$ converges as $\eps\to 0$ are spaces of space-\emph{time} distributions, and it follows that evaluating the limit process at a fixed time is in general not well defined.
Finally, \cite[(2.18)]{BrunedChandraChevyrecHairer2017} ensures that the analytic BPHZ theorem of \cite{ChandraHairer2016} can be applied, which in particular establishes the existence of a limit model. Most notably, this assumption rules out divergent variances in the "trees" used to build the regularity structure.

\begin{remark}
As will be clear from the proof, the existence of Malliavin derivatives is essentially automatic, as soon as the equation can be lifted to an abstract fixed point problem in some regularity structure, compare Theorem \ref{thm:local:H:diff} below. In particular, this result applies in principle to equations where parts of Assumption \ref{ass:main} are violated. For instance, one could even treat non-Gaussian noises in this framework. Since the setting of \cite{BrunedChandraChevyrecHairer2017} is already quite general, we decide to work under the assumptions introduced there in order not to over-complicate the presentation.
\end{remark}

We recall from \cite[Thm.~2.13]{BrunedChandraChevyrecHairer2017} that under Assumption \ref{ass:main} there exists a unique maximal stopping time $\tau=\tau(\omega)$ and a maximal solution $u = (u_i)_{i\le m}$ on $[0,\tau(\omega)) \times \T^d$ to (\ref{eq:singular:SPDE}). To be more precise, there exists a choice of constants $c^\eps_k$  for $\eps>0$ and a sequence of stopping times $\tau^\eps=\tau^\eps(\omega)$ with $\tau^\eps \to \tau$ in probability as $\eps \to 0$ and such that the classical solution $u^\eps$ to (\ref{eq:singular:SPDE:reg}) with $\Upsilon_i^k$ given as in (\ref{eq:counter-term}) exits almost surely on $[0,\tau^\eps)$, and such that for $T>0$ the sequence $u^\eps$ conditioned on the event $\{ \tau > T\}$ converges as $\eps \to 0$ to $u$ in probability in the space of space-time distributions $\CD'([0,T)\times \T^d)$.
When restricted to positive times, this convergence also takes place in the H\"older-Besov space\footnote{See Section \ref{sec:singularSPDEs} for the definition of the function $\reg:\{1, \ldots m\}\to \R$.}
$\CC^\reg:=\bigoplus_{i \le m} \CC^{\reg(i)}((0,T) \times \T^d)$. Moreover, the stopping time $\tau$ can be chosen maximal, in the sense that the statement above does not hold for any stopping time $\tilde\tau$ such that $\tilde\tau>\tau$ with positive probability.

\subsection{Main Results}

We want to study the finite-dimensional law of the random variables given by either testing the solution $u$ against test-functions or by evaluating at a finite number of space-times points, that is, we study the law of
\begin{align}
(u_i(\phi_l^i))_{l \le L, i\le m} \in \R^{L \times m}
\quad\qquad \text{ resp. } \quad\qquad
(u_i(z_l^i))_{l \le L, i \le m}\in \R^{L \times m}
\end{align}
for some $L \in \N$, test-functions $\phi_l^i \in \CC_c^\infty(\R\times \T^d)$ and space-time points $z_l^i \in \R^+ \times \T^d$. We will establish Malliavin differentiability \cite{Malliavin1997,Nualart2006} of these random variables, and a fortiori study the existence of densities with respect to Lebesgue measure. As has already been observed in \cite{CannizzaroFrizGassiat2015} and later in \cite{GassiatLabbe2017}, the classical notion of Malliavin differentiability is to strong for our purposes, as it imposes moment bounds which are simply not true in general in our setting. Instead, we are lead to use a version of Malliavin differentiability more adapted to this setting, and we borrow the notion of local $H$-Fr\'echet differentiability from \cite[Def.~3.3.1]{UstunelZakai2000}, which we recall in Definition \ref{def:local:H1:diff} below. Denoting by $H^\xi$ the Cameron-Martin space for the jointly Gaussian random fields $\xi=(\xi_i)_{i \le n}$, our main result on Malliavin differentiability reads as follows.

\begin{theorem}\label{thm:main:malliavin}
Under Assumption~\ref{ass:main}, let $u$ be the solution to (\ref{eq:singular:SPDE}) given by \cite[Thm.~2.13]{BrunedChandraChevyrecHairer2017}, let $\tau=\tau(\omega) \in (0,\infty]$ be the time of existence of $u$, let $\psi:=\lim_{\eps\to 0}\psi^\eps$ and assume that $\psi^\eps$ and $\psi$ and are locally $H^\xi$-Fr\'echet differentiable for any $\eps>0$. Then, for any $T>0$ and any $i \le m$ the solution $u_i$ restricted to $(0,T) \times \T^d$ and conditioned on $\{ T < \tau \}$
is locally $H^\xi$-Fr\'echet differentiable with values in $\CC^{\reg(i)}( (0,T) \times \T^d )$ in the sense of Definition~\ref{def:local:H1:diff}. 
The $H^\xi$-derivative $D_h u_i$ of $u_i$ in the direction of $h \in H^\xi$ is given by $v_{i,h}$, where $v_h$ is the renormalized solution to the equation
\begin{multline}\label{eq:tangent:equation:main}
\partial_t v_{i,h} = \CL_i v_{i,h}
+
DF_i(u, \nabla u, \ldots)(v_{h}, \nabla v_{h}, \ldots) 
\\
+ 
\sum_{j \le n} DF_i^j(u, \nabla u, \ldots)(v_{h}, \nabla v_{h}, \ldots) \xi_j
+
\sum_{j \le n} F_i^j(u, \nabla u, \ldots) h_j,
\end{multline}
with initial condition $v_{i,h}(0) = D_{h} u_{i,0}$.
\end{theorem}

We refer the reader to (\ref{eq:tangent:equation:reg}) below for a precise formulation of what we mean by \emph{renormalized solution} to (\ref{eq:tangent:equation:main}).

Local $H$-Fr\'echet differentiability is a powerful tool to establish existence of densities due an argument by Bouleau and Hirsch \cite{BouleauHirsch1986}, see also \cite[Sec.~2.1.3]{Nualart2006} and the references therein. We show existence of densities under some simplifying assumptions which we introduce in Section \ref{sec:densities} below. These assumptions are somewhat technical and we refrain from stating them precisely at this stage. Instead, we refer the reader to the paragraph below Theorem \ref{thm:main:density} for an informal discussion of these assumptions and to Proposition \ref{prop:main:density} for a class of interesting equations for which these assumptions are indeed satisfied. Taking Assumptions \ref{ass:dual:equation:simplicitiy}, \ref{ass:counter-terms} and \ref{ass:density:space-time} from Section \ref{sec:densities} for granted at the moment, our main result concerning densities is the following.

\begin{theorem}\label{thm:main:density}
Assume that Assumptions \ref{ass:dual:equation:simplicitiy}, \ref{ass:counter-terms} and \ref{ass:density:space-time} below hold. Let furthermore $(\xi_i)_{i \le n}$ be a family of jointly Gaussian noises on some probability space $(\Omega,\P)$ with the property that the Cameron-Martin space $H^{\xi}$ of $\xi$ is dense in $L^2((0,\infty) \times \T^d)^n$. Let also $T>0$ and assume that $\P(\{ T<\tau \})>0$.

Then, for any $L\in \N$ and any family $(\phi_l)_{l\le L}$ with $\phi_l \in \CC_c^\infty((0,T) \times \T^d; \R^m)$ of \emph{linearly independent}, smooth, compactly supported $\R^m$-valued test-functions, one has that the $\R^{Lm}$-valued random variable
\[
(\dual{u_1^1,\phi_1}, \ldots, \dual{u_m^1,\phi_1}, \ldots \ldots,
 \dual{u_1^L,\phi_L}, \ldots, \dual{u_m^L,\phi_L})
\]
conditioned on $\{ T < \tau \}$ has a density with respect to Lebesgue measure. 
\end{theorem}

We briefly discuss the assumption of the previous theorem. Assumption \ref{ass:dual:equation:simplicitiy} severly limits the explicit dependence of the right hand side of (\ref{eq:singular:SPDE:reg}) on the derivatives of the solution. This is done mainly for convenience, as it simplifies many computations below. Assumption \ref{ass:counter-terms} ensures that the renormalization constants for the "dualized" tangent equation are identical to the constants appearing in the original tangent equation, which is needed in order to pass to the limit $\eps\to 0$ in the dual equation. We believe that both of these assumptions are not really necessary and it will be the subject of future research to establish a density result that does not require them. Assumption \ref{ass:density:space-time} on the other hand is crucial, as it ensures that in case of multiplicative noise the term multiplying the noise does not make it degenerate. 

Instead of giving the precise assumptions at this stage we will limit ourselves to the following examples in order to demonstrate the scope of Theorem \ref{thm:main:density}.

\begin{proposition}\label{prop:main:density}
Assumptions \ref{ass:dual:equation:simplicitiy}, \ref{ass:counter-terms} and \ref{ass:density:space-time} hold for the $\Phi^p_2$ equations, the $\Phi^4_3$ and the $\Phi^{4}_{4-\eps}$ equation\footnote{The $\Phi^4_{4-\eps}$ equation is of form (\ref{eq:example:phi}) in $1+4$ dimensions with $p=3$, but with noise that has slightly better regularity then space-time white-noise, compare \cite[Sec.~2.8.2]{BrunedChandraChevyrecHairer2017}. This equation becomes critical for space-time white noise.}. It also holds for the multiplicative stochastic heat equation
\begin{align}\label{eq:SHE}
\partial_t u = \Delta u + f(u) + g(u)\xi
\end{align}
in $1+1$ dimension, as soon as the smooth function $g$ does not vanish anywhere on $\R$. In particular, the statement of Theorem \ref{thm:main:malliavin} holds as soon as $\xi$ is a Gaussian noise with Cameron-Martin space dense in $L^2(\R\times \T^d)$.
\end{proposition}

\begin{remark}
We remark that in \cite{GassiatLabbe2017} the authors obtained existence of densities for the $\Phi^4_3$ equation under the same assumptions as above. Additionally, they obtained existence of densities for noises whose Cameron-Martin space is not dense in $L^2$, but are such that they are everywhere "rough  enough" in a certain sense. We expect that it is possible to generalize these arguments to all of examples given above.

The existence of densities for 2D PAM established in \cite{CannizzaroFrizGassiat2015} on the other hand falls out of our setting for two reasons. One is that the authors obtained existence of densities for point evaluations and the other is that PAM is driven by purely spatial white noise. On the other hand, the approach to show non-degeneracy of the Malliavin derivative used in \cite{CannizzaroFrizGassiat2015} draws on a maximum principle and thus uses extensively the specific structure of PAM. Additionally, in \cite{CannizzaroFrizGassiat2015} the authors only study density for the evaluation at a single space-time point.
\end{remark}

\subsection{Application: Multiplicative Stochastic Heat Equation}

We apply our results to the stochastic heat equation (\ref{eq:SHE}) driven by a space-time  dependent Gaussian noise $\xi$ on $\R \times \T$ satisfying the assumptions of Section \ref{sec:gaussian:noises} and vector fields $f,g \in \CC^\infty(\R)$ with $g>0$. We refrain from stating the precise assumptions on the noise at this point as they are somewhat convoluted, but we note that these assumptions allow in particular the case of space-time white-noise. The regularized and renormalized equation is given by
\[
\partial_t u^\eps = \Delta u^\eps + f(u^\eps) + g(u^\eps)\xi^\eps
+
C_\eps^1 g'(u^\eps)g(u^\eps)
+
C_\eps^2 g'(u^\eps)^3 g(u^\eps)^2
+
C_\eps^3 g''(u^\eps) g'(u^\eps) g(u^\eps)^2 
\]
for some constants $C_\eps^i$ for $i=1,2,3$, subjection to (for simplicity deterministic) initial condition $u^\eps(0)=u_0$. For space-time white noise $\xi$, this equation was first derived in \cite{HairerPardoux2015}, where it was also shown that in this case one can choose $C_\eps^2$ and $C_\eps^3$ independent of $\eps$.
For more general noises, it follows from \cite{BrunedChandraChevyrecHairer2017} that given some initial condition $u_0 \in \CC^{\frac{1}{2}}(\T)$ there exists a choice of constants $C_\eps^i$, $i=1,2,3$, such that the regularized solution $u^\eps$ conditioned on $\{ \tau>T \}$ converges to some limit $u$ in probability the space $\CC^{\frac{1}{2}-\kappa}((0,T)\times \T^d)$.

By Theorem \ref{thm:main:malliavin}, the solution is locally $H^\xi$-Fr\'echet differentiable, and its derivative $v_h = D_h u$ satisfies the tangent equation
\[
\partial_t v_h = \Delta v_h + f'(v_h) + g'(v_h)\xi + g(v_h) h.
\]
More precisely, one has $v_h = \lim_{\eps \to 0}v_h^\eps$, where $v_h^\eps$ is the classical solution to
\begin{multline*}
\partial_t v^\eps_h = \Delta v^\eps_h + f'(u^\eps)v^\eps_h + g'(u^\eps)v^\eps_h \xi^\eps + g(u^\eps) h^\eps
\\
\Big(
+
C_\eps^1 \big( g''(u^\eps)g(u^\eps) + g'(u^\eps)^2 \big)
+
C_\eps^2 \big( 3g''(u^\eps)g'(u^\eps)^2 g(u^\eps)^2 + 2g'(u^\eps)^4 g(u^\eps) \big)
\\
+
C_\eps^3 \big( g'''(u^\eps) g'(u^\eps) g(u^\eps)^2 
	+ g''(u^\eps)^2 g(u^\eps)^2 + 2g''(u^\eps) g'(u^\eps)^2 g(u^\eps) \big)
\Big) 
v^\eps_h,
\end{multline*}
subject to the initial condition $v_h^\eps(0) = 0$.

Furthermore, assuming that the Cameron-Martin space $H^\xi$ of $\xi$ is dense in $L^2(\R \times \T)$, then for any family of linearly independent test function $(\varphi_i)_{i \le L}$ with $\vphi_i \in \CC_c^\infty((0,T) \times \T)$ the $\R^L$-valued random variable given by $( \dual{u,\vphi_1}, \ldots , \dual{u, \vphi_L} )$ conditioned on the event $\{ \tau>T \}$ admits a density with respect to Lebesgue measure.

\subsection{Outline of the paper}

In section \ref{sec:reg:structures} and \ref{sec:singularSPDEs} we recall the notations, definitions and results about the theory of regularity structures which were developed in \cite{Hairer2014,BrunedChandraChevyrecHairer2017,ChandraHairer2016,BrunedChandraChevyrecHairer2017}. In Section \ref{sec:gaussian:measure:theory} we review some classical results about Gaussian measure theory in infinite dimensional spaces. 

As in \cite{CannizzaroFrizGassiat2015}, we first construct in an algebraic step an extended regularity structure in Section \ref{sec:ext:reg:structure} by adding for any noise-type $\Xi$ a symbol $\hat \Xi$ that acts as an abstract place-holder  for a fixed Cameron-Martin function. The extended set of trees is then given by allowing any appearance of any noise-type $\Xi$ in any tree to be replaced by $\hat\Xi$. In Section~\ref{sec:ext:models} we perform the main analytic argument which shows that for fixed Cameron-Martin function $h_\Xi$ and Gaussian noise $\xi_\Xi$ we can indeed define a renormalized model that in particular has the property that $\bold\Pi\Xi = \xi_\Xi$ and $\bold\Pi \hat\Xi = h_\Xi$, and this model is locally Lipschitz continuous in $h$. An extended model can then be mapped in a locally Lipschitz continuous way onto a "shifted" model in Section \ref{sec:shifts:models}, which in particular shows that the model behaves in a continuous way under shifting the noise by a Cameron-Martin function. In Sections \ref{sec:abstract:fixed:point} and \ref{sec:abstract:fixed:point:shift} we show how to lift and shift abstract fixed point problems. This will in particular allow us to consider for fixed Cameron-Martin function $h$ the equations driven by $\xi + rh$ for any $r \in \R$ in an $r$-independent model, and is thus suited to study G\^ateaux differentiability of the solution map in Cameron-Martin directions. G\^ateaux and Fr\'echet differentiability are then established in Section~\ref{sec:Gateaux} and~\ref{sec:H:diff}, respectively. Finally, in Section~\ref{sec:H:diff:SPDE} this abstract theory is applied to singular SPDEs of the type (\ref{eq:singular:SPDE}) under Assumption \ref{ass:main}, and we derive in particular the tangent equation~(\ref{eq:tangent:equation:main}).

In order to establish the existence of densities we study the dual equation (\ref{eq:dual:equation}) of the tangent equation (\ref{eq:tangent:equation:main}). We want to lift the dual equation again to an abstract fixed point problem, and since the dual equation is a stochastic PDE going \emph{backward} in time, we are led to construct another extension of the regularity structure, this time extending the set of kernel-types by adding for any type $\ft$ a type $\Kdual \ft$ representing the dualized kernel.  We then derive in Section~\ref{sec:dual:FPP} an abstract fixed point problem for the dual equation and we identify its reconstruction as the actual solution to the dual equation in Section~\ref{sec:dual:id}. This step is not automatic, since it is a-priori not clear that the renormalization constants obtained in these two ways coincide (it is not even clear that they differ by something of order $1$ in a suitable sense, which is the main reason that Section \ref{sec:densities} is less general then the rest of the paper).
This identification relies on Assumption \ref{ass:counter-terms} which basically enforces the identity that we need, and we show in Section \ref{sec:dual:single:equation} that Assumption \ref{ass:counter-terms} is satisfied when considering single equations (as opposed to systems of equations). Finally, we derive the existence of densities in a spirit similar to \cite{GassiatLabbe2017} by showing that the solution to the dual equation does not vanish identically in Section~\ref{sec:existence:densities}.

\paragraph{Acknowledgements.}
The author would like to thank Giuseppe Cannizzaro, Ajay Chandra,  and Martin Hairer for helpful discussions during the preparation of this article. The author acknowledges funding through Martin Hairer's ERC consolidator grant, project 615897.

\section{Setting and Notation}

\subsection{General Conventions on Notation}
\label{sec:notation}

We introduce some notation that is used throughout this article. Given $M\in \N$ we write $[M]:=\{1, \ldots, M \}$. We fix a spatial dimension $d\ge 1$ and write $\domain:= \R\times \T^d$. Given $z \in \domain$ we often write $z=(z_0,z_1,\ldots,z_d)$ with $z_0 \in \R$ and $(z_1,\ldots z_d) \in \T^d$. Given a finite set $A$, a subset $B\ssq A$, and a variable $z \in \domain^A$, we write $z_B := (z_a)_{a \in B}$.
We also fix a space-time scaling $\fs:\{0, \ldots , d\} \to \N$, and we write $|\fs|:=\sum_{i=0}^d \fs(i)$ for the effective space-time dimension. For a multi-index $k \in \N^{\{0, \ldots d \}}$ we write $|k|_\fs := \sum_{i =0}^d \fs(i)k_i$, and for $z \in \domain$ we write $|z|_\fs := \sum_{i=0}^d |z_i|^{\frac{1}{ \fs(i) } }$. We use the convention that sums of the form 
\[
\sum_{|k|_\fs\le r} \cdots
\]
always run over all multi-indices $k \in \N^{\{0,\ldots,d\} }$ with $|k|_\fs \le r$.

We write $\CC_c^\infty(\domain)$ for the space of compactly supported, smooth functions $\phi : \domain \to \R$, we endow this space with the topology given by the system of semi-norms
\[
\|\phi\|_{K,r} := \sup_{z\in K} \sup_{|k|_\fs \le r} |\partial^k \phi(z)|
\]
for $K\ssq \domain$ compact and $r \in \N$, and we write $\CD'(\domain)$ for the dual space of $\CC_c^\infty(\domain)$. We call $\rho \in \CC_c^\infty(\domain)$ a \emph{mollifier} if $\int \rho(x) dx=1$, and in this case we define
\[
\rho^{(\eps)}(z) := \eps^{-|\fs|} \rho( \frac{z_0}{\eps^{\fs(0)} }, \ldots, \frac{z_d}{\eps^{\fs(d)} }).
\]

Finally, the following terminology of multi-sets will be useful. A multiset $\tm$ with values in $A$ is an element of $\N^A$. Given two multisets $\tm , \tn \in \N^A$ we write $\tm\sqcup\tn \in \N^A$ for the multiset given by $(\tm\sqcup\tn)(a) := \tm(a) + \tn(a)$, and we write $\fm \sqsubset \fn$ if $\fm \le \fn$. We also naturally identify a subset $B\ssq A$ with the multiset $\I_B : A \to \N$. 
Given any finite set $I$ and a map $\varphi : I \to A$ we write $[I,\varphi]$ for the multiset with values in $A$ given by
\begin{align}\label{eq:notation:multiset}
[I,\varphi]_a := \#\{ i \in I : \varphi(i) = a \}
\end{align}
for any $a \in A$. 

We sometimes discuss concepts in detail for concrete examples in order to clarify notation. In these cases we use notations of the form $[a,b,c,\ldots]$, to denote multisets. For instance, we write $[a,a,b]:=2\I_a + \I_b$.

Given a multiset $\tm$ as above and a function $f$ on $A$ we also freely use the notation $\sum_{a \in \tm} f(a)$ and $\prod_{a \in \tm} f(a)$. These expression should be interpreted as
\[
\sum_{a \in \tm} f(a) := \sum_{a \in A} \tm(a) f(a)
\qquad\text{ and }\qquad
\prod_{a \in \tm} f(a) := \prod_{a \in A} {f(a)}^{\tm(a)}.
\] 

Sometimes it will be useful to consider functions $f$ whose domain is formally given by $\domain^\tm$ for some multiset $\tm$. Setting $\CM(\tm):=\{(a,k) : a \in A, 1 \le k \le \tm(a) \}$, when we write $f: \domain^\tm \to \R$ we really mean that $f: \domain^{\CM(\tm)} \to \R$ is a function which is symmetric under any permutation $\sigma$ of $\CM(\tm)$ with the property that for any $(a,k) \in \CM(\tm)$ there exits $l \le \tm(a)$ such that $\sigma(a,k) = (a,l)$.

\subsection{Regularity Structures}
\label{sec:reg:structures}

In this section we recall the main notations and results about regularity structures that we will use in the sequel. Throughout this paper we assume we are given a finite set of types $\FL=\FL_-\sqcup \FL_+$. The finite set $\FL_+$ will index the components of the equation, while the finite set $\FL_-$ will index the Gaussian noises appearing on the right hand side of the equation. We assume that $\FL$ is equipped with a homogeneity assignment $|\cdot|_\fs:\FL_\star \to \R^\star$ for $\star \in \{+,-\}$. Recall from \cite[Def.~5.7]{BrunedHairerZambotti2016} that a \emph{rule} $R$ is a collection $(R(\ft))_{\ft \in \FL}$ that assigns to any type $\ft \in \FL$ a set of multisets $R(\ft)$ with values in $\FL\times \N^{d+1}$. We recall the notions of \emph{normal}, \emph{sub-critical} and \emph{complete} from \cite[Def.~5.7, Def.~5.14,  Def.~5.22]{BrunedHairerZambotti2016}. Let us especially recall that a rule is subcritical if there exists a map $\reg : \FL \to \R$ with the property that
\[
\reg(\ft) < |\ft|_\fs + \inf_{N \in R(\ft)} \reg(N),
\]
where we set $\reg(N) := \sum_{ (\fl,k) \in N}  (\reg(\fl) - |k|_\fs)$ for any multiset $N \in \N^{\FL \times \N^{d+1}}$. A rule $R$ is normal if $R(\ft)$ is stable under taking multi-subsets of any $N\in R(\ft)$, and additionally if $R(\ft):=\{\emptyset \}$ for any $\ft \in \FL_-$. Completeness ensures that the set of tree conforming to the rule $R$ (c.f. \cite[Def.~5.8]{BrunedHairerZambotti2016}) is stable under the action of renormalization.

\begin{example}\label{ex:rule}
In case of stochastic heat equation (\ref{eq:SHE}) one has a unique kernel type $\ft$ and a unique noise type $\Xi$ and the rule $R$ is given by $R(\Xi):= \{ \emptyset \}$ and $R(\ft)$ contains all multisets of the form
\[
[(\ft,0),\ldots]\qquad or \qquad[(\Xi,0), (\ft,0), \ldots]
\]
where we write $(\ft,0),\ldots$ denote an arbitrary (possible vanishing) number of occurrences of $(\ft,0)$.
\end{example}

We assume we are given a normal, subcritical, complete rule $R$ and we denote by $\CT^\ex$ the regularity structure constructed as in \cite[Def.~5.26]{BrunedHairerZambotti2016}. We will actually work with a slightly simplified structure as far as the extended decoration is concerned, compare Section~\ref{sec:trees} below. We extend the homogeneity assignment $|\cdot|_\fs$ to $\CT^\ex$ in the usual way, taking into account the extended decoration\footnote{In the notation of \cite[Def.~5.3]{BrunedHairerZambotti2016} this was denoted by $|\cdot|_+$.}, and we write $\CT \ssq \CT^\ex$ for the reduced regularity structure obtained as in \cite[Sec.~6.4]{BrunedHairerZambotti2016}. We will very rarely need the homogeneity assignment that neglects the extended decoration, but in these situations we will denote this by $|\cdot|_-$ as in \cite[Def.~5.3]{BrunedHairerZambotti2016}. We write $\TT^\ex$ and $\TT$ for the set of trees in $\CT^\ex$ and $\CT$, respectively, so that $\CT^\ex$ and $\CT$ are freely generated by $\TT^\ex$ and $\TT$ as linear spaces. We write $\TT^\ex_\alpha$ for the set of trees $\tau \in \TT^\ex$ with the property that $|\tau|_\fs = \alpha$, we write $\CT^\ex_\alpha := \linspace{\TT^\ex_\alpha}$, and for $\gamma \in \R$ we write $\CQ_{<\gamma}$ for the projection of $\CT^\ex$ onto $\bigoplus_{\alpha<\gamma} \CT^\ex_\alpha$.

Finally, we make the following Assumption on the regularity structure, which is needed to apply the results of \cite{ChandraHairer2016}.

\begin{assumption}\label{ass:main:reg}
For any tree $\tau \in \TT$ one has 
\begin{align}
|\tau|_\fs > 
	\big( -\shalf \big) 
	\lor \max_{u \in L(\tau)} |\ft(u)|_\fs
	\lor \big( - |\fs| - \min_{\Xi \in \FL_-} |\Xi|_\fs \big).
\end{align}
\end{assumption}

\begin{remark}
Under Assumption \ref{ass:main} we are indeed in this setting, compare \cite[Sec.~5.5]{BrunedHairerZambotti2016} and \cite[Sec.~3.1]{BrunedChandraChevyrecHairer2017}. In particular, Assumption~\ref{ass:main:reg} follows from \cite[(2.18)]{BrunedChandraChevyrecHairer2017}. However, even though singular SPDEs are the application that we have in mind, most arguments in the proof of Theorem \ref{thm:main:malliavin} are carried out on the level of the regularity structure and do not really require that we are in the setting of an actual singular SPDE.
\end{remark}

Finally, we make the simplifying assumption on the rule that we do not allow products or derivatives of noises to appear on the right hand side of the equation. As was already remarked in \cite{ChandraHairer2016} and \cite{BrunedChandraChevyrecHairer2017}, such an assumption does not seem to be crucial but simplifies certain arguments.

\begin{assumption}\label{ass:noises:derivatives:products}
We assume that for any $\ft \in \FL$ and any $N \in R(\ft)$ there exists at most one pair $(\Xi,k) \in \FL_- \times \N^{d+1}$ such that $N_{(\Xi,k)} \ne 0$, and this case $k=0$ and $N_{(\Xi,0)}=1$. 
\end{assumption}

\subsubsection{Trees}\label{sec:trees}

Trees $\tau \in \TT^\ex$ can be written as typed, decorated trees $\tau=(T^{\fn,\fo}_\fe, \ft)$, where $T$ is a rooted tree with vertex set $V(T)$, edge set $E(T)$ and root $\rho_T$, the map $\ft$ assigns types to edges and is formally a map $\ft:E(T)\to \FL$, and the decorations $\fn,\fe,\fo$ are maps $\fn:N(T) \to \N^d$, $\fe : E(T) \to \N^d$ and $\fo : N(T) \to (-\infty,0]$. We call $\fo$ the \emph{extended decoration}. Here we define the decomposition of the set of edges into $E(T) = L(T) \sqcup K(T)$ with $e \in L(T)$ (resp. $e \in K(T)$) if and only if $\ft(e) \in \FL_-$ (resp. $\ft(e) \in \FL_+$), and we write $N(T) \ssq V(T)$ for the set of $u \in V(T)$ such that there does not exist $e \in L(\tau)$ such that $u=e^\uparrow$.
We will often abuse notation slightly and leave the type map $\ft$ and the root $\rho_\tau$ implicit. We recall that the relation between the homogeneity assignments $|\cdot|_\fs$ and $|\cdot|_-$ is given by
$
|\cdot|_\fs = |\cdot|_- + \sum_{u \in N(T)}\fo(u),
$
so that in particular one has $|\tau|_\fs \le |\tau|_-$ for any tree $\tau \in \TT^\ex$.

On a rooted tree $T$ we define a total order $\le$ on $V(T)$ by setting $u \le v$ if and only if $u$ lies on the unique shortest path from $v$ to the root $\rho_T$ and we write edges $e \in E(T)$ as order pairs $e=(e^\uparrow, e^\downarrow)$ with $e^\uparrow \ge e^\downarrow$. If $u \in V(T) \backslash\{\rho_T \}$, then there exists a unique edge $e \in E(T)$ such that $u = e^\uparrow$, and in this case we write $u^\downarrow := e$.  Recall that it follows from the fact that $R$ is normal (c.f. \cite[Def.~5.7]{BrunedHairerZambotti2016}) that elements $u \in V(T)\backslash N(T)$ are leaves of the tree $T$. 

Given a typed, decorated tree $\tau$ as above, $k \in \N^{d+1}$ and $\ft \in \FL_+$ we write
$
\CJ_\ft^k \tau
$
for the planted, decorated, typed tree obtained from $\tau$ by attaching an edge $e$ to the root with type $\ft$ and $\fe(e)=k$. We write $\TT_{\ft} \ssq \TT$ for the set of trees $\tau \in \TT$ such that $\CJ_\ft^0 \tau \in \TT$.

\begin{example}\label{ex:introduce:tree}
Throughout the paper we will consider examples from stochastic heat equation (\ref{eq:SHE}) whenever we need to clarify notations. In particular, we often consider the tree $\treeExampleSmall$, where we introduce the following graphical conventions:
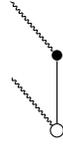
\begin{align*}
\treeExample
\qquad\qquad
\parbox[b] {200pt}{
\drawroot ... root $\rho({\treeExampleSmall})$, element of $N(\treeExampleSmall)$
\\
\drawnode ... node, element of $N(\treeExampleSmall)$
\\
\starOneKernelShort ... edge of kernel type, element of $K(\treeExampleSmall)$
\\
\starOneNoiseShort ... edge of noise type, element of $L(\treeExampleSmall)$
}
\end{align*}
\end{example}

\subsubsection{Algebraic Notation} 

We use the notation $\CT_-^\ex$, $\hat\CT_-^\ex$, $\CT_+^\ex$, $\hat\CT_+^\ex$, $\CG_-^\ex$ and $\CG_+^\ex$ from \cite{BrunedHairerZambotti2016} for the respective spaces defined in \cite[Def.~5.26,~(5.23),~Def.~5.36]{BrunedHairerZambotti2016}, and we write $\CG_-$ for the reduced renormalization group as in \cite[Thm.~6.28]{BrunedHairerZambotti2016}. We recall that $\CT_-^\ex$ and $\CT_+^\ex$ form Hopf algebras and $\CG_-^\ex$ and $\CG_+^\ex$ are defined as their respective character groups.
We use the notation $\Delta_-^\ex$ and $\Delta_+^\ex$ for the co-products for negative and positive renormalization respectively, as in \cite[Cor.~5.32]{BrunedHairerZambotti2016}, and we write $\tilde \CA_-^\ex : \CT_-^\ex \to \hat\CT_-^\ex$ and $\tilde\CA_+^\ex : \CT_+^\ex \to \hat \CT_+^\ex$ for the twisted antipodes defined in \cite[Prop.~6.2,~Prop.6.5]{BrunedHairerZambotti2016}. 

%

\subsubsection{Models}

We assume that for any $\ft \in \FL_+$ we are given a decomposition of the Green's function into $G_\ft=K_\ft + R_\ft$ with $R_\ft \in \CC^\infty(\domain)$ and such that $K_\ft \in \CC_c^\infty(\domain\backslash\{0\})$ satisfies \cite[Ass.~5.1, Ass.~5.4]{Hairer2014}, and given the kernel assignment $(K_\ft)_{\ft \in \FL_+}$ we recall the definition of admissible models \cite[Def.~2.7,~Def.~8.29]{Hairer2014}. We call a model $Z=(\Pi,\Gamma)$ \emph{smooth} if $\Pi_x \tau \in \CC^\infty(\domain)$ for any $\tau \in \TT^\ex$ and any $x \in \domain$, and we call $Z$ \emph{reduced} if $\Pi_x \tau$ does not depend on the extended decoration of $\tau$. Given an admissible tuple $\bold \Pi\tau \in \CC^\infty(\domain)$ for $\tau \in \CT^\ex$ we write $\CZ(\bold \Pi)$ for the model constructed as in \cite[(6.11),(6.12)]{BrunedHairerZambotti2016}, whenever this is well defined, and we write $\CM_\infty$ for the set of smooth, reduced, admissible models for $\CT^\ex$ of the form $\CZ(\bold \Pi)$. We write $\CM_0$ for the closure of $\CM_\infty$ in the space of models, and given a probability space $(\Omega,\P)$ we write $\CM_\infty^\rand$ and $\CM_0^\rand$ for the spaces of $\CM_\infty$ and $\CM_0$ valued random variables on $(\Omega,\P)$, respectively, endowed with the topology induced by convergence in probability. We write $\Omega_\infty:= \Omega_\infty(\FL_-):= \CC^\infty(\domain)^{\FL_-}$ and given $f \in \Omega_\infty$ we write $Z^f= \CZ(\bold \Pi^f)$ for the canonical lift of $f$ to a model $Z^f \in \CM_\infty$, c.f. \cite[Rem.~6.12]{BrunedHairerZambotti2016}.  We finally write $\Omega_0:= \CD'(\domain)^{\FL_-}$.


\begin{remark}
Again, we remark that under Assumption \ref{ass:main} we are in this setting, compare \cite[Ass.~2.6]{BrunedChandraChevyrecHairer2017}.
\end{remark}

\subsubsection{Gaussian Driving Noises}
\label{sec:gaussian:noises}

Given a probability space $(\Omega,\P)$ we write $\SM_\infty:= \SM_\infty(\FL_-)$ for the space of $\Omega_\infty$-valued centred, stationary, jointly Gaussian random fields $\eta$ on $(\Omega,\P)$. We want to introduce a class of $\Omega_0$-valued Gaussian noises that we are going to consider in the sequel, and for this we introduce the following notation.
Given a $\Omega_0$-valued jointly Gaussian, stationary, centred random noise $\eta$ on $(\Omega,\P)$, we denote by $\cov_{\ft,\ft'}^\eta \in \CD'(\domain)$ the distributional covariance of $\eta_\ft$ and $\eta_{\ft'}$ defined via the identity 
\[
\E[\eta_\ft(\varphi) \eta_{\ft'}(\psi)]
=
\cov_{\ft,\ft'}^\eta(\int \varphi(x- \cdot)\psi(x)dx)
\]
for any $\varphi,\psi \in \CC_c^\infty( \domain )$. We note that this is well defined by stationarity. The next definition if motivated by the assumptions made in \cite{ChandraHairer2016}. In order to state it we fix for any $k \in \N^{d+1}$ a function $P_k \in \CC_c^\infty(\domain)$ such that $P_k(x) = x^k$ in a neighbourhood of the origin. 

\begin{definition}
We write $\FC(\FL_-)$ for the space of families of kernels $(\cov_{\ft,\ft'})_{\ft,\ft' \in \FL_-}$ such that $\cov_{\ft,\ft'}^* = \cov_{\ft',\ft}$ in the sense that one has 
\[
\dual{ \cov_{\ft,\ft'} * \phi, \psi }_{L^2} =  \dual{ \phi, \cov_{\ft',\ft} * \psi }_{L^2} 
\]
for any $\ft,\ft' \in \FL_-$ and any choice of test functions $\phi,\psi \in \CC_c^\infty(\domain)$, such that $\cov$ is non-negative definite in the sense that one has
\[
\sum_{\ft} \dual{ \sum_{\ft'} \cov_{\ft,\ft'} * \phi_{\ft'}, \phi_\ft }_{L^2} \ge 0
\] 
for any family $\phi \in \CC_c^\infty(\domain)^ {\FL_-}$, 
such that the singular support of the distribution $\cov_{\ft,\ft'} \in \CD'(\domain)$ is contained in $\{0\}$, and denoting by $\hat\CC_{\ft, \ft'}$ the smooth function representing $\CC_{ \ft, \ft'}$ away from the origin, we require that
\begin{itemize}
\item for any test function $\varphi \in \CC_c^\infty(\domain)$ such that $D^k \varphi(0)=0$ for any $|k|_\fs < -|\ft|_\fs - |\ft'|_\fs-|\fs|$, one has $\CC_{\ft, \ft'}(\varphi) = \int \varphi(x) \hat\CC_{\ft, \ft'}(x) dx$; and
\item there exists $\theta>0$ such that one has $\|\cov\|_{|\cdot|_\fs} <\infty$.
\end{itemize}
Here, we define the quantity
\begin{align}
\|\cov\|_{|\cdot|_\fs} := 
	\td(\cov)+
	\max_{\ft,\ft' \in \FL_-}
	\sup_{x \in \domain\backslash\{0\}}
	\supss{k \in \N^{d+1} \\ |k|_\fs\le 6|\fs|}
	|D^k \hat\cov_{\ft, \ft'}(x) | |x|^{-|\ft|_\fs-|\ft'|_\fs+|k|_\fs-\theta}
\end{align}
with
\[
\td(\cov) := \max_{\ft,\ft' \in \FL_-}
\supss{k \in \N^{d+1} \\ |k|_\fs < -|\ft|_\fs-|\ft'|_\fs-|\fs|} 
|\cov_{\ft,\ft'} (P_k)|.
\]
\end{definition}
We write $\SM_0 = \SM_0(\FL_-)$ for the set of $\Omega_0$-valued, jointly Gaussian, centred, stationary random fields $\eta$ with the property that  $\cov_{\ft,\ft'}^\eta \in \FC(\FL_-)$, and we write $\|\eta\|_{|\cdot|_\fs}:=\|\cov^\eta\|_{|\cdot|_\fs}$.
\begin{remark}
Note that while $\FC(\FL_-)$ is a linear space, $\SM_0$ is not, due to the non-linearity of the map $\eta \mapsto \cov^\eta$.
\end{remark}
Given an element $\eta \in \SM_0$ and a mollifier $\rho \in \CC_c^\infty(\domain)$ with $\int \rho(x) dx = 1$, we call the sequence $\eta^\eps:= \eta * \rho^{(\eps)} \in \SM_\infty$ an \emph{approximation} of $\eta$. We say that a map $X$ from $\SM_\infty$ into a topological space $\CX$ \emph{extends continuously} to $\SM_0$ if there exists a map $\hat X:\SM_0 \to \CX$ that extends $X$ and is such that whenever $\eta^\eps \in \SM_\infty$ is an approximation of $\eta \in \SM_0$ in the sense above, then $X(\eta^\eps) \to X(\eta)$. We finally note that by Kolmogorov's continuity theorem for any $\xi \in \SM_0(\FL_-)$ there exists a version such that $\xi_\ft$ is an element of $\CC_\fs^{|\ft|_\fs}(\domain)$ almost surely for any $\ft \in \FL_-$.

\subsubsection{Modelled Distributions} \label{sec:modelled:dist}

We recall the terminology and notation from \cite[Sec.~3,~Sec.~6]{Hairer2014} of modelled distributions. Given a model $Z \in \CM_0$ and $\gamma>0$, $\eta \in \R$ we write $\CD^{\gamma,\eta}$ for the space of singular modelled distributions defined in \cite[Def.~6.2]{Hairer2014} allowing a singularity at the hyperplane $t=0$. More precisely $\CD^{\gamma,\eta}$ consists of all maps $f:\domain \to \CT_{<\gamma}$ with the property that, setting $P=\{z \in \domain: z_0 = 0\}$, one has that
\begin{align}\label{eq:bound:D:gamma:eta}
\|f\|_{\gamma,\eta,K} := 
\sup_{ z, \bar z }
\sup_{\beta<\gamma}
\bigg(
\frac
	{\|f(z)\|_\beta}
	{ ( |z_0|^{\frac{1}{\fs_0}} \land 1 )^{(\eta-\beta)\land 0}  }
+
\frac{\|f(z) - \Gamma_{z, \bar z} f(\bar z)\|_\beta}
{\|z-\bar z\|_\fs^{\gamma-\beta} (|z_0|^{\frac{1}{\fs_0}}\land |\bar z_0|^{\frac{1}{\fs_0}} \land 1)^{\eta-\gamma} }
\bigg)
\end{align}
is finite for any compact $K \ssq \domain$. Here, the first supremum runs over all $z,\bar z \in K\backslash P$ with the property that $\| z - \bar z \|_{\fs} \le |z_0|^{\frac{1}{\fs_0}}\land |\bar z_0|^{\frac{1}{\fs_0}}$.
Given a sector $V$ of $\CT$ we write $\CD^{\gamma,\eta}_V$ for the space of $f \in \CD^{\gamma,\eta}$ such that $f(x) \in V$ for any $x \in \domain$, and we write $\CD^{\gamma,\eta}_V(Z)$ if we want to emphasise the underlying model. Often we want to consider localized version of these spaces that contain functions $f$ that only live on a bounded time interval $[0,T)$ for some $T>0$. We write $\CD^{\gamma,\eta,T}$ for the space of all functions $f: [0,T)\times \T^d \to \CT_{<\gamma}$ satisfying the bound (\ref{eq:bound:D:gamma:eta}) for any compact $K \ssq [0,T) \times \T^d$. The notation $\CD^{\gamma,\eta,T}_V$ and $\CD^{\gamma,\eta,T}_V(Z)$ then have meanings analogue to above.

On the spaces $\CD^{\gamma,\eta}_V$ for $\gamma>0$ we denote by $\CR:\CD^{\gamma,\eta}_V\to \CC^{\alpha \land \eta}_{|\cdot|_\fs}(\domain)$ the reconstruction operator defined in \cite[Prop.~6.9]{Hairer2014}, provided that $\alpha \land \eta > -|\fs|+\fs(0)$, where $\alpha\le 0$ denotes the regularity of the sector $V$.

Finally, we denote for any $\ft \in \FL_+$ by $\CK_\ft$ the operator constructed in \cite[(5.15)]{Hairer2014} acting between $\CD^{\gamma,\eta}_V$ and $\CD^{\gamma+|\ft|_\fs,(\eta\land\alpha)+|\ft|_\fs}$ for any sector $V$ of regularity $\alpha$ such that one has $\CR \CK_\ft f = K_\ft * \CR f$ for any $f \in \CD^{\gamma,\eta}$. We also define the operator $\CP_\ft f:= \CK_\ft f + R_\ft * \CR f$ for any $f \in \CD^{\gamma,\eta}$. 
Of course, the operators $\CK_\ft$ and $\CP_\ft$ depend slightly on $\gamma$. Since $\gamma$ will always be clear from the context, we leave it implicit in this notation.

In \cite[Sec.~6]{Hairer2014} basic properties of certain maps (multiplication, differentiation, integration, composition) between space of modelled distributions were derived and we summarize them in Proposition \ref{prop:Dgamma:loc:Lipschitz} below. 

\subsubsection{BPHZ Theorem}
\label{sec:BPHZ:theorem}
Given a smooth noise $\eta \in \SM_\infty$ we define as in \cite[(6.23)]{BrunedHairerZambotti2016} a character $g^\eta_-$ on $\hat\CT_-^\ex$ by setting $g^\eta_-(\tau) := \E(\bold \Pi^\eta \tau)(0)$ for any tree $\tau\in\hat\CT_-^\ex$, and extending this linearly and multiplicatively, where $\bold \Pi^\eta$ is such that $Z^\eta = \CZ(\bold \Pi^\eta)$ is the canonical lift of $\eta$. We then define the BPHZ-character $g^\eta_\BPHZ \in \CG_-$ as in \cite[(6.24)]{BrunedHairerZambotti2016} by setting
\[
g^\eta_\BPHZ(\tau) := g^\eta_-(\tilde \CA^\ex_- \tau)
\]
for any $\tau \in \TT_-^\ex:=\{\tau \in \TT^\ex: |\tau|_-<0\}$, and extending this linearly and multiplicatively. For any character $g \in \CG_-^\ex$ we use the notation $M^g:\CT^\ex \to \CT^\ex$ for the linear operator given by
\[
M^g := (g \otimes \Id)\Delta_-^\ex, 
\]
and for a smooth noise $\eta \in \SM_\infty$ we set
\[
\hat{\bold \Pi}_\BPHZ^\eta \tau := \bold\Pi^\eta M^{g_\BPHZ^\eta}\tau
\]
for any $\tau \in \TT^\ex$, and we defined the BPHZ-renormalized model $\hat Z^\eta_\BPHZ:= \CZ(\hat {\bold \Pi}_\BPHZ^\eta)$, compare \cite[Thm.~6.17]{BrunedHairerZambotti2016}.

The following is then a direct consequence of \cite{ChandraHairer2016}.
\begin{theorem}\label{thm:BPHZ}
The map $\eta \mapsto \hat Z^\eta_\BPHZ$ extends continuously to a map from $\SM_0$ into $\CM_0^\rand$.
\end{theorem}
\begin{proof}
See \cite{ChandraHairer2016}.
\end{proof}


\subsection{Singular SPDEs}
\label{sec:singularSPDEs}

In \cite{BrunedChandraChevyrecHairer2017} the authors established a black box theorem for solving a large class of singular SPDEs of the form (\ref{eq:singular:SPDE}). We briefly recall the notations introduced in this paper, as far as we are going to need it later on.
In order to unify the notation, we assume that $\# \FL_+ = n$ and $\# \FL_- = m$ and we write $(u_\ft)_{\ft \in \FL_+}$ and $(\xi_\ft)_{\ft \in \FL_-}$ rather than $(u_i)_{i \le n}$ and $(\xi_j)_{j \le m}$. 
We recall that we assume that for $\ft \in \FL_+$ we are given a differential operator $\CL_\ft$ involving only spatial derivatives and such that $\partial_t - \CL_\ft$ admits a Green's function $G_\ft$ satisfying Assumption \cite[Ass.~2.6]{BrunedChandraChevyrecHairer2017}.
Furthermore, we recall from \cite[Sec.~2.5]{BrunedChandraChevyrecHairer2017} that we assume we are given two functions $\reg: \FL\to \R$ and $\ireg : \FL_+ \to \R$ satisfying \cite[Def.~2.2]{BrunedChandraChevyrecHairer2017} and \cite[Ass.~2.4]{BrunedChandraChevyrecHairer2017}. 
We define
\[
\CC^\ireg := \bigoplus_{\ft \in \FL_+} \CC^{\ireg(\ft)}_\fs (\T^d)
\]
and we write $\CS_{\rho,\eps}^-:\bigoplus_{\ft \in \FL_-} \CC_\fs^{|\ft|_\fs}(\domain) \to \bigoplus_{\ft \in \FL_+} \CC^\infty(\domain)$ for the map given by \cite[(A.11)]{BrunedChandraChevyrecHairer2017}, and we assume that the initial condition $u_0^\eps$ is of the form
\[
u^\eps_0(\xi) = \CS_{\rho,\eps}^-(\xi)(0,\cdot) + \psi^\eps(\xi)
\]
for a sequence of $\CC^\ireg$-valued random fields $\psi^\eps$ such that $\psi^\eps \to \psi$ in probability as $\eps \to 0$.

We now recall the definition of the counter-terms appearing on the right hand side of (\ref{eq:singular:SPDE:reg}). For this we borrow some more notation from \cite{BrunedChandraChevyrecHairer2017}. 

For $\ft \in \FL_+$ we often write $F_\ft^\bullet := F_\ft$ in order to avoid case distinctions. The smooth functions $F^\Xi_\ft$ are allowed to depend on $D^k u_\fl$ where $\fl \in \FL_+$ and where $k \in \N^{d+1}$ ranges over a finite set of multi-indices, say $|k|_\fs\le r$. Consequently, it makes sense to define for any $\fl \in \FL_+$ and $k\in N^{d+1}$ the derivative $D_{(\fl,k)}F_\ft^\Xi$ of $F_\ft^\Xi$ in the direction of $(D^k u_\fl)$. We will reserve the symbol $\partial$ for derivatives in direction of space-time variables.

For any tree $\tau=(T^{\fn,\fo}_\fe,\ft) \in \CT^\ex$ and node $\mu \in N(\tau)$ we write $\Xi[\mu]:=\ft(e)$ if there exists a (necessarily unique) edge $e \in L(\tau)$ with $e^\downarrow = \mu$. We write $\Xi[\mu] := \bullet$ otherwise. Moreover, we write $n[\mu]:=\# \{e \in K(\tau): e^\downarrow = \mu \}$, and we write $e_j[\mu]$, $j\le n[\mu]$ for the $n[\mu]$ distinct edges $e \in K(\tau)$ such that $e^\downarrow = \mu$. Note that this is uniquely defined up to order of $e_j[\mu]$. Finally, we write $\ft_j[\mu] := \ft(e_j[\mu])$ and $k_j[\mu]:=\fe(e_j[\mu])$ for the type and derivative decoration of $e_j[\mu]$, respectively. Note that any tree $\tau \in \CT$ can now be written in the form
\begin{align}\label{eq:tree:normal:form}
\tau = 
X ^{ \fn( \rho_\tau ) }
\Xi[\rho_\tau]
\prod_{j=1}^{n[\rho_\tau]}
\CJ_{\ft_j[\rho_\tau]} ^{ k_j[ \rho_\tau ] }
[\tau_j]
\end{align}
for some decorated, typed trees $\tau_j \in \CT$.

\begin{definition}\label{def:ft:non:vanishing}
For $\fl \in \FL_+$ we say that a tree $\tau \in \TT$ is \emph{$\fl$-non-vanishing} if for any $\mu \in N(\tau)$ one has that
\begin{align}\label{eq:ft:non:vanishing}
\Big(\partial ^{\fn(\mu)} \prod_{j=1} ^{n[\mu]}
	D_{(\ft_j[\mu], k_j [\mu])}
\Big) 
		F_{\ft(\mu)}^{\Xi[\mu]}
(u, \nabla u, \ldots)
\end{align}
does not vanish identically for any smooth function $u:\domain \to \R^{\FL_+}$, where we set $\ft(\rho_\tau):=\fl$, and $\ft(\mu) := \ft(\mu^\downarrow)$ if $\mu \ne \rho_\tau$.
We write $\TT_\fl^F$ for the set of tree $\tau \in \TT$ that are $\fl$-non-vanishing and are such that $\CJ_{\fl}^0[\tau] \in \TT$, and we write $\TT_{\fl,-}^F$ for the set of $\tau \in \TT_\fl^F$ such that $|\tau|_\fs <0$, and we write $\tilde \TT_\fl^F := \TT_{\fl,-}^F \sqcup\{ \bullet \}$.
\end{definition}

It follows from a straight forward inductive argument that the definition of $\fl$-non-vanshing given above coincides with \cite[Def.~2.9]{BrunedChandraChevyrecHairer2017}. We now define the counter-terms appearing in the renormalized equation.

\begin{definition}
For $\fl \in \FL_+$ and $\tau \in \tilde \TT_\fl^F$ we define the function
\begin{align}\label{eq:counter-term}
\Upsilon_\fl^F[\tau]
:=
\prod_{\mu \in N(\tau)}
\Bigg(
\partial ^{\fn(\mu)} \prod_{j=1} ^{n[\mu]}
	D_{(\ft_j[\mu], k_j [\mu])} 
		F_{\ft(\mu)}^{\Xi[\mu]}
(u, \nabla u, \ldots)
\Bigg)
\end{align}
where we use the notation from Definition \ref{def:ft:non:vanishing} for the type $\ft(\mu) \in \FL_+$.
\end{definition}

\begin{example}
In case of stochastic heat equation (\ref{eq:SHE}) one has $F_\ft^\bullet = f$ and $F_\ft^\Xi = g$. Moreover, for the tree $\tau = \treeExampleSmall$ from Example \ref{ex:introduce:tree} one has
\[
\Upsilon_\ft^F[\treeExampleSmall](u)
=
g'(u)g(u).
\]
\end{example}

Given a character $g \in \CG_-$ and a smooth noise $\eta \in \SM_\infty(\FL_-)$ we define the $g$-renormalization of (\ref{eq:singular:SPDE}) by 
\begin{align}\label{eq:singular:SPDE:g:renorm}
\partial_t u_{\ft} = \CL_\ft u_\ft
+
F_\ft(u, \nabla u, \ldots)
+
\sum_{\Xi \in \FL_-}
	F_\ft^\Xi( u , \nabla u , \ldots ) \eta_\Xi
+
\sum_{\tau \in \TT_{\ft,-}^F}
	\frac{ g(\tau) }{ S(\tau) }
		\Upsilon_\ft^F[\tau]( u, \nabla u, \ldots).
\end{align}

Here $S(\tau) \in \N$ is a symmetry factor explicitly given in \cite[(2.16)]{BrunedChandraChevyrecHairer2017}.
We are thus in the setting of (\ref{eq:singular:SPDE:reg}) with $K = \# \TT_{\ft,-}^F$.
Given initial conditions as above and letting $u^\eps$ be the solution to (\ref{eq:singular:SPDE:g:renorm}) with $\eta$ replaced by $\xi^\eps$ and $g$ replaced by $g_\BPHZ^\eps$, then the statement of \cite[Thm.~2.13]{BrunedChandraChevyrecHairer2017} precisely says that $u^\eps$ converges to some limit $u$ in probability as $\eps \to 0$.

\subsection{Gaussian Measure Theory}\label{sec:gaussian:measure:theory}

In this section we review basic facts about Gaussian measures on infinite dimensional spaces as far as it is needed for the purpose of this paper. We follow in this section mostly the lecture notes \cite{Hairer2009notes} for basic properties of the Cameron-Martin space. For more details we refer the reader to standard literature \cite{Nualart2006}. 
Given a separable Fr\'echet space $\CX$ we call a centered probability measure $\mu$ on $\CX$ equipped with the Borel sigma field \emph{Gaussian} if all finite dimensional projections of $\mu$ are Gaussian. 
A Gaussian measure $\mu$ is uniquely determined by its covariance operator $\cov_\mu : \CX^* \to \CX$ defined via the identity $l^*(\cov_\mu(k^*))= \E^\mu[k^*l^*]$ for any $k^*,l^* \in \CX^*$. 
We denote the image of the covariance operator by $\topcirc{H_\mu} \subseteq \CX$ and we equip this space with a scalar product given by $\langle h,k \rangle_\mu := \E^\mu[ k^* l^*]$ where $h^*,k^* \in \CX^*$ are such that $ \CC _\mu (h^*)=h$ and $ \CC _\mu (k^*)=k$. The closure $H_\mu$ of $\topcirc{H_\mu}$ under the norm induced from this scalar product is a Hilbert space know as \emph{Cameron-Martin space.} It is well known that $H_\mu \hookrightarrow \CX$ continuously, the space $H_\mu$ determines the Gaussian measure $\mu$ uniquely, and one has the following classical result due to Cameron and Martin \cite{CameronMartin1944}.
\begin{theorem} \emph{(Cameron-Martin)}
	Let $\mu$ be a Gaussian measure on some Fr\'echet space $\CX$ and define for $h \in \CX$ the operator $T_h : \CX \to \CX$ by $T_h (x):= x + h$. Then one has $(T_h)_* \mu \ll \mu$ if and only if $h \in H_\mu$ and in this case one has  $(T_h)_* \mu \simeq \mu$.
\end{theorem}

We want to use non-degeneracy of the pathwise derivative in Cameron-Martin directions to establish existence of densities. The classical Malliavin derivative imposes moment bounds on this derivative which are not available in our setting. An approach more adapted to our situation is the notion of local Cameron-Martin Fr\'echet differentiability, in the form introduced in \cite[Def.~3.3.1]{UstunelZakai2000}. 
\begin{definition}\label{def:local:H1:diff}
Let $\CY$ be a Banach space. We call a $\CY$-valued random variable $X$ on $(\CX,\mu)$ is \emph{locally $H_\mu$-Fr\'echet differentiable} if there exists a $\mu$-null set $N$ such that for any $\omega \in N^c$ the map $H_\mu \to \CY$, $h \mapsto X( \omega + h )$, is Fr\'echet differentiable in a $H_\mu$-neighbourhood of the origin.
\end{definition}

We note that this version of local $H$-Fr\'echet differentiability was also used in \cite[Def.~2.2]{CannizzaroFrizGassiat2015} and \cite[Def.~4.1]{GassiatLabbe2017}.
If $X$ is locally $H_\mu$-Fr\'echet differentiable we denote its $H_\mu$-Fr\'echet derivative by $DX$. The main motivation for this definition is the criterion by Bouleau-Hirsch \cite{BouleauHirsch1986} for the existence of densities.
In order to deal with situations in which we the solution does not exist globally, we use a slightly generalized version, and for this we make the following construction. Let $U \ssq \CX$ be a measurable subset of $\CX$. We say that $U$ is \emph{$H_\mu$-open} if for any $ x \in U$ there exists $\eps>0$ such that for any $h \in H_\mu$ with $\|h\|_{H_\mu}<\eps$ one has $x+h \in U$. We fix an $H_\mu$-open set $U$ and we define for $\eps>0$ the set $U_\eps$ as the $\eps$-involution of $U$ in $H_\mu$, i.e.
\[
U_\eps:= \{ x \in U : \forall h \in H_\mu \text{ with } \|h\|_{H_\mu} \le \eps \text{ one has } x+h \in U \}.
\]
We then assume that there exists a sequence of locally $H_\mu$-Fr\'echet differentiable random variables $\varphi_\eps: \CX \to \R$ for $\eps>0$ that approximates the indicator function $\I_U$ from the inside in the sense that one has $0 \le \vphi_\eps \le 1$, and $\vphi_\eps(x)=1$ for any $x \in U_\eps$ and $\vphi_\eps(x)=0$ for $x \in U^c$. If such a sequence exists, we say that $U$ can be \emph{approximated from the inside}.

\begin{theorem}\label{thm:Bouleau:Hirsch}
\emph{(Bouleau-Hirsch)}
Let $X$ be an $\R^n$-valued random variable on a Gaussian probability space $(\CX,\mu)$ with separable Cameron-Martin space $H_\mu$, and let $U \ssq \CX$ be an $H_\mu$-open measurable subset of $\CX$ such that $U$ can be approximated from the inside. Let moreover $\Omega_X \ssq \CX$ be the event that $DX$ has full rank, and assume that $\mu(U \cap \Omega_X)>0$. Then $X$ conditioned on the event $U \cap \Omega_X$ admits a density with respect to Lebesgue measure.
\end{theorem}
\begin{proof}
See \cite[Prop.~2.4]{CannizzaroFrizGassiat2015}.
\end{proof}

We will see that the time of existence $\tau$ is lower semi-continuous with respect to $H_\mu$-translations, see Remark \ref{rem:lower:semi:continuity} below, which implies in particular that the event $\{ \tau > T \}$ is $H_\mu$-open. The fact that $\{ \tau>T \}$ can be approximated from the inside can be shown exactly as in \cite[Lem.~5.3]{CannizzaroFrizGassiat2015}, and this leads to the following.

\begin{corollary}\label{cor:Bouleau:Hirsch}
Under Assumption \ref{ass:main}, let $T>0$, let $u \in \CC^\reg((0,T)\times \T^d)$ be the solution to a singular SPDE of the form (\ref{eq:singular:SPDE}) on a Gaussian probability space $(\Omega,\P)$, let $X: \CC^\reg((0,T)\times\T^d) \to \R^n$ be a $C^1$ map, and assume that $\P( \Omega_{X(u)} \cap \{ \tau>T \})>0$. Then $X(u)$ restricted on $\Omega_{X(u)} \cap \{ \tau>T\}$ admits a density with respect to Lebesgue measure.
\end{corollary}

We now implement the general constructions from this section in the situation that the underlying Fr\'echet space is given by $\Omega=\bigoplus_{\ft \in \FL_-}\CD'(\domain)$, and the Gaussian probability measure $\P$ on $\Omega$ is a stationary, centred Gaussian measure such that its covariance $\cov_\P$ is an element of $\FC(\FL_-)$. Note in particular that in this case the random field 
$\xi$ which $\P$-almost surely agrees with the identity on $\CX$ is an element of $\SM_0$. Given such a Gaussian measure $\P$, we will henceforth use the convention that $\xi$ denotes this particular random Gaussian field, while $\eta$ will still be used to denote more general random fields on $\Omega$ whose laws under $\P$ are Gaussian. We will usually leave the measure $\P$ implicit in the notation and one should always think of $\P$ as arbitrary but fixed.
It is straightforward to see that the space $\topcirc H^\xi$ is then given by 
\[
\topcirc H^\xi =
\set{
	(\sum_{\ft' \in \noit} \cov_{\ft,\ft'}^\xi*\varphi_{\ft'})
		_{ \ft \in \noit} 
	:
	\varphi_{\ft'} \in \CC_c^\infty{(\domain)} \text{ for any } \ft' \in \noit
}.
\]
It follows in particular that one has $\topcirc H^\xi \subseteq \detnoise{\noit}$.
 
We finish this section by introducing the following terminology.

\begin{definition}
Given a kernel $\cov\in \FC(\FL_-)$ 
 we write $ \topcirc H[\FL_-,\cov]$ for the space given by 
\[
\topcirc H [\FL_-,\cov] =
\set{
	(\sum_{\ft' \in \noit} \cov_{\ft,\ft'}*\varphi_{\ft'})
		_{ \ft \in \noit} 
	:
	\varphi_{\ft'} \in \CC_c^\infty(\domain) \text{ for any } \ft' \in \noit
},
\]
endowed with the scalar product  $\dual{\cdot,\cdot}_{H[\FL_-,\cov]}$ given by
\begin{align}\label{eq:HcovNorm}
 \dual{h,k}_{H[\FL_-,\cov]} := \sum_{\ft, \ft' \in \FL_-} 
\int_{\domain\times \domain} dx dy
\,
 \cov_{\ft,\ft'}(x-y) \varphi_{\ft'}(x)\psi_\ft(y)
\end{align}
for any $h,k \in \topcirc H[\FL_-,\cov]$, where $\varphi,\psi \in \CC_c^\infty(\domain)$ are such that $h=(\sum_{\ft' \in \noit} \cov_{\ft,\ft'}*\varphi_{\ft'})_{ \ft \in \noit}$ and $k=(\sum_{\ft' \in \noit} \cov_{\ft,\ft'}*\psi_{\ft'})_{ \ft \in \noit}$.
Finally, we write $H[\FL_-,\cov]$ for the Hilbert space given as the closure of $\topcirc H[\FL_-,\cov]$ under the induced norm.
\end{definition}

It is then not hard to see that for any $\xi$ as above one has that $\topcirc H[\FL_-,\cov^\xi]$ and $H[\FL_-,\cov^\xi]$ agree with $\topcirc H^\xi$ and $H^\xi$, respectively.

\section{Extension and Translation of Models}

In this section we introduce the main technical tools and show key estimates needed to prove Theorem \ref{thm:main:malliavin}.
On the technical level, pathwise differentiability (in Cameron Martin directions) of solutions to singular SPDEs can be effectively studied by introducing an extended regularity structure. The basic idea, which was already used in \cite{CannizzaroFrizGassiat2015}, is to extend the regularity structure $\CT^\ex$ by adding for any noise type $\Xi$ a new noise type $\hat \Xi$, which plays the role of an abstract place holder for a fixed Cameron-Martin function. We perform this extension in two separate steps. We first construct in a purely  algebraic step, using the formalism developed in \cite{BrunedHairerZambotti2016}, an extended regularity structure. Afterwards we show in an analytic step, building up on the result from \cite{ChandraHairer2016}, that any fixed Cameron-Martin function $h$ can be indeed be "lifted" to a renormalized extended model, and, crucially, this lift is locally Lipschitz continuous in $h$.

For any $r \in \R$ the original regularity structure then maps into the extended structure via a map $\SS_r$, which is essentially the multiplicative extension of the map $\Xi \mapsto \Xi + r\hat\Xi$. Conversely, any extended model maps onto a model for the original structure via the "dual" map $\SS_r^*$, which can be viewed as implementing this shift on an analytic level.

At the end of this section we are going to lift abstract fixed point problems to the extended regularity structure, and, using the shift operator, we also make sense of shifted fixed point problems.

\subsection{Extension of the Regularity Structure }\label{sec:ext:reg:structure}

In the sequel it will be useful to consider general extensions of the set of noise types, and we are led to make the following general construction. Given a finite set $I$ we define a new set of noise types by $\noisetra{I}:=\FL_- \times I$, and we write $\extxt{ \FL}_-^I:=\noisetypes \sqcup \noisetra{I}$. We call $\extxt{ \FL}_-^I$  the \emph{extended} set of noise types.
We also define the \emph{extended} set of types by $\extxt\FL^I:=\extxt\FL_-^I \sqcup \FL_+$.
There exists a natural map $\cq: \extxt{\FL}^I \to \FL$, which acts as the identity on $\FL$ and removes the second variable on $\hat\FL_-^I$, i.e. one has $\cq(\Xi,i):= \Xi$ for any $\Xi \in \FL_-$ and any $i \in I$. We extend the homogeneity assignment $|\cdot|_\fs$ to $\noiseextxt{I}$ by setting $|(\Xi,i)|_\fs:=|\Xi|_\fs$ for any $\Xi \in \FL_-$ and any $i \in I$. In order to avoid case distinctions, we will sometimes add a distinct element $\star$ to the index set by setting $I^\star:= I \sqcup\{\star \}$, and we identify $\extxt\FL_-^I$ with $\FL_- \times I^\star$.

Starting from the set of noise-types $\noiseextxt{I}$ and the (unchanged) set of kernel-types $\FL_+$, we can consider an extension of the rule $R$ to a rule $\extxt {R}^I$, which is defined by allowing any appearance of any noise types $\Xi \in \noisetypes$ being replaced by any extended noise type of the form $(\Xi,i)$ for $i \in I$. 

To be more precise, 
with the notation $\extxt \CN^I:= \N^{\extxt\FL^I \times \N^{d+1}}$, we define a rule $\extxt{R}^I: \extxt\FL^I \to 2^{\extxt\CN^I}\backslash\{\emptyset \}$ by setting
\begin{align}\label{eq:rule:extended}
\extxt R^I (\ft) := \{ N \in \extxt\CN^I : \cq N \in R(\ft)  \}
\end{align}
for any $\ft \in \FL_+$, and $\extxt R^I(\ft) := \{\emptyset \}$ for $\ft \in \extxt\FL_-^I$. Here, we define $\cq N \in \N^{\FL\times \N^{d+1}}$ by setting
\begin{align}\label{eq:qN}
(\cq N)_{(\ft, k)} := 
\sum_{\cq \tilde \ft = \ft} N_{(\tilde \ft,k)}
\end{align}
for any $\ft \in \FL$ and any $k \in \N^{d+1}$, where the sum runs over all $\tilde \ft \in \extxt\FL^I$ with $\cq \tilde \ft = \ft$. 
The following Lemma shows that one can construct a regularity structure $\extxt\CT^\ex[I]$ starting from the extended rule $\extxt R^I$ as in \cite{BrunedHairerZambotti2016}.
\begin{lemma}\label{lem:rule:extension}
For any finite set $I$ the rule $\extxt R ^I$ is a complete and subcritical rule. In particular, we can define the extended regularity structure $\extxt \CT^\ex[I]$ as in \cite[Sec.~5.5]{BrunedHairerZambotti2016}. Then $\extxt\CT^\ex[I]$ coincides with the span of all decorated trees $\tau = (T^{\fn,\fo}_\fe,\ft)$ with $\ft : E(T) \to \extxt\FL$ and with the property that $(T^{\fn,\fo}_\fe, \cq\ft) \in \CT^\ex$.
\end{lemma}
\begin{proof}
In order to see that $\extxt R^I$ is subcritical, recall from \cite[Def.~5.14]{BrunedHairerZambotti2016} and the fact that $R$ is subcritical that there exists a function $\reg: \FL \to \R$ with the property that
\begin{align}\label{eq:reg:extxt}
\reg(\ft) < |\ft|_\fs + \inf_{N \in R(\ft)} \reg(N)
\end{align}
for any $\ft \in \FL$. We extend $\reg$ to a function $\reg : \extxt\FL^I \to \R$ by setting $\reg(\ft) := \reg(\cq \ft)$ for any $\ft \in \extxt\FL^I$.
Then one has $\reg (N) = \reg (\cq N)$, where $\cq N$ is as in (\ref{eq:qN}), and thus the fact that (\ref{eq:reg:extxt}) holds for any $\ft \in \extxt\FL^I$ is a trivial consequence from the respective bound for $\cq \ft$ and the fact that $|\cq\ft|_\fs = |\ft|_\fs$.
Completeness (c.f. \cite[Def.~5.20]{BrunedHairerZambotti2016}) is a little tedious to verify, but completely straight forward.  
\end{proof}

Note that we could always consider the completion of $\extxt R^I$ as in \cite[Prop.~5.21]{BrunedHairerZambotti2016}, so that showing completeness is not really crucial.
The construction in \cite[Sec.~5.5,Sec.6]{BrunedHairerZambotti2016} results in a number of spaces which are all completely determined by the rule $\extxt R^I$. We adopt the convention that we use the notation $\extxt X[I]$ to denote a space $X$ constructed from $\extxt R^I$, and we sometimes drop  $[I]$ from the notation, whenever $I$ is clear from the context. In particular, we write $\extxt{\hat\CT}_-^\ex[I]$ and $\extxt{\hat\CT}_+^\ex[I]$ for algebras constructed in \cite[Def.~5.26]{BrunedHairerZambotti2016}, we write $\extxt{\CT}_-^\ex[I]$ and $\extxt{\CT}_+^\ex[I]$ for the Hopf algebras constructed in \cite[(5.23)]{BrunedHairerZambotti2016}, and we write $\extxt{\CG}_{+}^\ex[I]$ and $\extxt{\CG}_{-}^\ex[I]$ for the character group of $\extxt\CT_+^\ex[I]$ and $\extxt\CT_-^\ex[I]$, respectively, compare \cite[Def.~5.36]{BrunedHairerZambotti2016}.

For $\tau \in \extxt\CT^\ex[I]$ we write $\hat L(\tau) := \{u \in L(\tau) : \ft(u) \in \hat \FL_-^I\}$.
One has the obvious embedding $\CT^\ex \hookrightarrow \extxt\CT^\ex[I]$, and multiplicatively extended we obtain a Hopf algebra monomorphism $\CT_-^\ex \hookrightarrow \extxt{\CT}_-^\ex[I]$. This embedding between the Hopf algebras induces a natural group monomorphism between their character groups $\CG_-^\ex \hookrightarrow \extxt{\CG}_-^\ex[I]$, which is defined by extending any character $g \in \CG_-^\ex$ in such a way that $g(\tau)$ vanished for any tree $\tau$ outside of $\CT_-^\ex$. We will use all of these embeddings implicitly, so that in particular we view $\CT_-^\ex$ as a sub Hopf algebra of $\extxt \CT_-^\ex[I]$, and we view $\CG_-^\ex$ as a subgroup of $\extxt{\CG}_-^\ex$.

Given an admissible family $\extxt{\bold \Pi}\tau \in \CC^\infty(\domain)$ for any $\tau \in \extxt\CT^\ex$ we write $\CZ(\extxt {\bold \Pi})$ for the admissible model constructed as in \cite[(6.11),(6.12)]{BrunedHairerZambotti2016} and we write similar to before $\extxt\CM_\infty^I$ and $\extxt\CM_0^I$ for the set of smooth reduced admissible models of the form $\CZ(\extxt{\bold \Pi})$ and for their closure in the space of models, respectively.
For any finite set $I$ we consider the space $\SM_\infty(\noisetypes) \times \detnoise{\noisetra{I}}$, and we write elements of this set as tuples $(\eta,h)$ or $(\eta,(h_i)_{i\in I})$, with $h \in \Omega_\infty(\hat\FL_-^I)$ and $h_i \in \CC^\infty(\domain)$, depending on the situation. Note that one has $(\eta,h) \in \detnoise{\noiseextxt{}}$ almost surely, and thus the canonical lift of $(\eta,h)$ to a random admissible model $\extxt Z^{\eta,h} \in \extxt\CM_\infty^I$ is well defined. Finally, we denote by 
\[
{\extxt {\hat Z}}_\BPHZ^{\eta,h} := M^{g^\eta_\BPHZ} \extxt Z^{\eta,h}
\]
the BPHZ renormalization. Here, we denote by $g^\eta_\BPHZ \in \CG_-^\ex \subseteq \extxt{\CG}_-^\ex$ the BPHZ character for the smooth stationary noise $\eta$, and we use the convention introduced above to view any character $g \in \CG_-^\ex$ also as a character of $\extxt\CT_-^\ex$.

The particular case that $I=\{1\}$ will play the most important role in the sequel. In this case we use the shorter notation $\hat \Xi:=(\Xi,1)$ for any $\Xi \in \FL_-$, 
and we write $\ext\CT^{\ex}:=\extxt{\CT}^{\ex}[\{1\}]$ and $\ext{\CG}_-^\ex:=\extxt{\CG}_-^{\ex}[\{1\}]$, and similar for the other spaces defined above.  We call $\ext\CT^\ex$ the \emph{onefold extension} of $\CT^\ex$.
More generally, if $I=\{1, \ldots , m\}$, then we call $\extxt\CT^{\ex}[I]$ the \emph{$m$-fold} extension of $\CT^\ex$.

Extensions of the set of noise-types can be used to conveniently encode shifts and differences between canonical lifts of smooth functions to models. This construction will allow us later in particular to almost automatically obtain Lipschitz bounds from uniform bounds applied to extended regularity structures. We make two constructions that we will use throughout the paper. For this, let $\CT^\ex$ be a regularity structure obtained from some noise-type set $\FL_-$, and let as above $\ext\CT^\ex$ be its onefold extension. 

The first constructions concerns "shifts" of models. For this, we introduce the operator $\SS: \CT^\ex \to \ext \CT^\ex$ by setting
\begin{align}\label{eq:shift:operator:A}
\SS (\tau,\ft) := \sum_{\tilde \ft} (\tau,\tilde \ft)
\end{align}
for any typed, decorated tree $(\tau,\ft) \in \CT^\ex$, where the sum over $\tilde\ft$ runs over all maps $\tilde \ft: E(\tau) \to \ext{\FL}_-$ such that $\cq \tilde\ft = \ft$. The shift operator algebraically encodes a binomial expansion of a tree $\tau$ when it is interpreted for the "shifted" noise $f+h$. The following Lemma shows in particular that this binomial expansion interacts  nicely with the action of renormalization. 
\begin{lemma}\label{lem:shift:Pi}
For any $h,k \in \Omega_\infty(\FL_-)$ write $Z^{h+k}=(\Pi^{h+k},\Gamma^{h+k})$ and $\ext Z^{h,k}= (\ext \Pi^{h,k}, \ext \Gamma^{h,k})$. Then, for any $g \in \CG_-^\ex$ one has the identity
\begin{align}\label{eq:shift:Pi}
\Pi^{h+k}_z M^g  = \ext \Pi^{h,k}_z M^g \SS .
\end{align}
on $\CT^\ex$ for any $z \in \domain$.
\end{lemma}
\begin{example}
As an example, consider the tree $\treeExampleSmall$ form Example \ref{ex:introduce:tree}. We graphically represent the shifted noise type $\hat\Xi$ by $\starOneNoiseShiftShort$, so that one has $\SS \treeExampleSmall = \treeExampleSmall + \treeExampleSmallShiftA + \treeExampleSmallShiftB + \treeExampleSmallShiftC$. The left and the right hand side of equation (\ref{eq:shift:Pi}) for $g = \bold 1^\star$ then read respectively $\Pi^{h+k}_z \treeExampleSmall = ((h+k)*K)(z) (h(z)+k(z))$ and $\ext\Pi^{h,k}_z \SS\treeExampleSmall = (h*K)(z) h(z) + (h*K)(z) k(z) + (k*K)(z) h(z) +(k*K)(z) k(z)$.
\end{example}
\begin{proof}
For the identity element $g = \bold 1^*\in \CG_-^\ex$, this follows directly from the definition of the canonical lift and the definition of the shift operator (\ref{eq:shift:operator:A}). 
Indeed, since the canonical lift is multiplicative, it suffices to show (\ref{eq:shift:Pi}) on planted trees $\tau$. If $\tau = \Xi$ for some $\Xi \in \FL_-$, then $\SS \tau = \Xi + \hat\Xi$, and one has $\Pi^{h+k}_z \tau = h+k = \extxt\Pi^{h,k}_z \Xi + \extxt\Pi^{h,k}_z \hat \Xi = \extxt\Pi^{h,k}_z \SS \tau$.
Finally, if $\tau = \CJ^{k}_\ft \sigma$ with $\ft \in \FL_+$, then one has $\SS \tau = \CJ^k_\ft \SS \sigma$, and the result follows inductively from the respective identity for $\sigma$ and the admissibility condition \cite[(8.19)]{Hairer2014}, see also \cite[(6.13)]{BrunedHairerZambotti2016}.

If $g\ne \bold 1^*$ we use the fact that $\SS$ commutes with the co-product (see Lemma \ref{lem:properties:shift:operator}) below) and the fact that by our convention $g$ vanishes outside of $\CT_-^\ex$, which by definition of $\SS$ implies in particular that one has the identity $g\SS = g$ on $\CT_-^\ex$. It follows that one has
\[
M^g \SS=
(g \otimes \Id) \Delta_-^\ex \SS =
(g \SS \otimes \SS) \Delta_-^\ex =
(g \otimes \SS) \Delta_-^\ex,
\]
on $\CT^\ex$. Comparing this with $M^g = (g \otimes \Id) \Delta_-^\ex$, the result follows by applying the first part of the proof to the right components of these tensor products.
\end{proof}

The second construction we are carrying out in this section concerns differences between shifts of canonical lifts. 
To this end we consider the sets $I:=\{\ch\}$ and $J:=\{ \ch, \ck, \ch-\ck \}$ and the corresponding extended noise-type sets $\extxt\FL_-^I=\FL_- \sqcup \FL_-\times I$ and $\noiseextxt{J}= \FL_- \sqcup \FL_- \times J$. We write elements in $\Omega_\infty({\noiseextxt{I}})$ and $\Omega_\infty({\noiseextxt{J}})$ as pairs $(h,h^\ch)$ and quadrupels  $(h,h^\ch, h^\ck, h^{\ch-\ck})$, respectively.

In order to construct the next operator, it will be helpful to fix for any tree $\tau \in \extxt \CT^\ex[I]$ a total order $\preceq$ on $\hat L(\tau)$.
We then define an linear operator $\AA:\extxt \CT^\ex[I] \to \extxt\CT^\ex[J]$ by setting for any typed, decorated tree $(\tau,\ft) \in \extxt\CT^\ex[I]$
\[
\AA (\tau,\ft) := \sum_{u \in \hat L(\tau)} (\tau, \ft[u]),
\]
where we define for any $u \in  \hat L(\tau)$ the type map $\ft{[u]}:E(\tau) \to \extxt\FL^J$ by setting
\begin{align}\label{eq:ft[u]}
\ft{[u]}_e :=
\begin{cases}
(\ft_e^{(1)}, \ch	)			&\text{ if }e \in \hat L(\tau), \, e\prec u	\\
(\ft_e^{(1)},\ch-\ck)		&\text{ if }e \in \hat L(\tau), \, e= u		\\
(\ft_e^{(1)},\ck)			&\text{ if }e \in \hat L(\tau), \, u\prec e	\\
\ft_e						&\text{ if }e \in E(\tau) \backslash \hat L(\tau)
\end{cases}
\end{align}
for any edge $e \in E(\tau)$. Here, $\ft_e^{(1)}$ denotes the first component of $\ft_e$, so that $\ft_e = (\ft_e^{(1)},\ch)$ for any $e \in \hat L (\tau)$. We write $\AA^\preceq$ and $\ft^\preceq[u]$ if we want to highlight the underlying family of total orderings on the sets $L(\tau)$ for $\tau \in \TT_-^\ex$ used to construct $\AA$. The point of this total ordering is that we want to expand the difference $\extxt\Pi^{f,h}_z - \extxt \Pi^{f,k}_z$ into a telescoping sum, and the statement below will be valid for any total order of $\hat L(\tau)$. We will use this telescopic sum later on in order to obtain almost automatically local Lipschitz bounds from uniform bounds.
The next Lemma shows in particular that this telescopic sum  interacts nicely with the action of renormalization.
\begin{lemma}\label{lem:diff:Pi}
For any $f,h,k \in \Omega_\infty(\FL_-)$ write $\extxt Z^{f,h}:= (\extxt \Pi^{f,h}, \extxt \Gamma^{f,h})$ and $\extxt Z^{f,h,k,h-k} := (\extxt \Pi^{f,h,k,h-k}, \extxt \Gamma^{f,h,k,h-k})$ for the canonical lifts of $(f,h)$ and $(f,h,k,h-k)$ to models in $\extxt\CM_\infty^I$ and $\extxt \CM_\infty^J$, respectively. Then, for any $g \in \CG_-^\ex$ one has the identity
\begin{align}\label{eq:diff:Pi}
\extxt\Pi^{f, h}_z M^g- \extxt\Pi^{f, k}_z M^g
=
\extxt \Pi^{f, h, k, h-k }_z  M^g \AA
\end{align}
on $\extxt \CT^\ex[I]$ for any $z \in \domain$.
\end{lemma}
\begin{proof}
Assume first that $g= \bold 1^*$ and fix $\tau \in \extxt \CT^\ex[I]$.
The proof is somewhat complicated by the fact that $\AA$ is not multiplicative with respect to the tree product due the arbitrary choice of $\preceq$. 
However, using an argument similar to the one given in the proof of Lemma \ref{lem:shift:Pi}, going inductively in the size of the tree, it is straight forward to see that one has the identity
\[
\extxt\Pi^{ f, h, k, h-k }_z (\tau, \ft^{\preceq}[u]) = 
\extxt\Pi^{ f, h, k, h-k }_z (\tau, \ft^{\preceq}_+[u] ) 
-  
\extxt\Pi^{ f, h, k, h-k }_z (\tau, \ft^{\preceq}_-[u] )
\]
for any $u \in \hat L(\tau)$, where we define $\ft^{\preceq}_+[u]$ and $\ft^{\preceq}_-[u]$ as in (\ref{eq:ft[u]}), but in the case $e=u$ with right hand side replaced by
\[
\ft^{\preceq}_+[u]_u := (\ft_u^{(1)},\ch)
\qquad\text{ and }\qquad
\ft^{\preceq}_-[u]_u := (\ft_u^{(1)},\ck).
\]
Now, directly from the definition we get that whenever $u,v \in \hat L(\tau)$ are adjacent with respect to $\preceq$ and such that $u \preceq v$, then one has that $\ft^{\preceq}_+[u] = \ft^{\preceq}_-[v]$, so that $\extxt\Pi^{ f, h, k, h-k }_z \AA\tau$ turns into a telescopic sum, which is equal to 
\[
\extxt\Pi^{ f, h, k, h-k }_z (\tau, \ft^{\preceq}_+[\max_{\preceq} L(\tau)] ) 
-  
\extxt\Pi^{ f, h, k, h-k }_z (\tau, \ft^{\preceq}_-[\min_{\preceq} L(\tau)] )
=
\extxt\Pi^{f, h}_z \tau - \extxt\Pi^{f, k}_z \tau. 
\]

For a general character $g \in \CG_-^\ex$, note that the first part of the proof implies in particular that $\extxt \Pi^{f, h, k, h-k }_z  \AA^{\preceq}$ is independent of $\preceq$. Since moreover by our convention the character $g$ vanished outside of $\CT_-^\ex$ it follows with an argument almost identically to the one given in Lemma \ref{lem:shift:Pi} for the shift operator $\SS$ that $\AA$ satisfies the following identity
\[
(g \otimes \extxt \Pi^{f, h, k, h-k }_z) \Delta_-^\ex \AA
=
(g \otimes \extxt \Pi^{f, h, k, h-k }_z \AA) \Delta_-^\ex
\]
on $\extxt{\CT}^\ex[I]$ 
\footnote{Note however that such an identity is \emph{not} true on the purely algebraic level, that is, without hitting the right component on both sides with the operator $\extxt\Pi^{ f, h, k, h-k }_z$. Indeed, since a tree $\tau$ can contain multiple subtrees which are \emph{different} as subtrees in $\tau$ but \emph{identical} as algebraic objects, there is in general no choice of total order $\preceq$ on $\hat L(\tau)$ for any tree $\tau \in \extxt\CT^\ex$ with the property that such a statement becomes true.}
. We conclude as in the proof of Lemma \ref{lem:shift:Pi}.
\end{proof}

\subsection{Extension of Models} \label{sec:ext:models}

We now assume that we are given a partition of the set of noise-types into $\FL_- = \FL_-^\rand \sqcup \FL_-^\det$. We want to consider noises which are random, centred, stationary and Gaussian for $\Xi \in \FL_-^\rand$ and deterministic for $\Xi \in \FL_-^\det$. To this end, we introduce the notation that given a pair $(\eta,f) \in \SM_\infty(\FL_-^\rand) \times \Omega_\infty(\FL_-^\det)$ we write $Z^{\eta,f}$ fo the canonical lift of the tupel $\eta\sqcup f$ to a random model. In such a situation, we furthermore want to consider a modification of negative renormalization that only takes into account diverging subtrees $\tau$  which have the property that all leaves $u \in L(\tau)$ have types $\ft(u) \in \FL_-^\rand$. Denoting the set of trees $\tau \in \TT_-^\ex$ with this property by $\TT_-^\ex[\FL_-^\rand]$, we define the character $g_\BPHZ^\eta \in \CG_-^\ex$ on trees $\tau \in \TT_-^\ex$ by setting $g_\BPHZ^\eta(\tau):= \E \bold \Pi^{\eta} \tilde \CA^\ex_- \tau (0)$ if $\tau\in \TT^\ex_-[\FL_-^\rand]$, and $g_\BPHZ^\eta(\tau)=0$ otherwise, and extending this linearly and multiplicatively. Finally, we define the renormalized model by
\[
\hat Z_\BPHZ^{\eta,f} := M^{g^\eta_\BPHZ} Z^{\eta, f}.
\]

We write $L^\rand(\tau)$ and $L^\det(\tau)$ for the set of $u \in L(\tau)$ with the property that $\ft(u) \in  \FL_-^\rand$ and $\ft(u) \in \FL_-^\det$, respectively.

\begin{remark}
We are mainly going to be interested in the setting where $\CT^\ex$ is itself a one-fold extension with noise types $\FL_- \sqcup \hat \FL_-$, and one has $\FL_-^\rand = \FL_-$ and $\FL_-^\det = \hat \FL_-$. Note that in this case the notation of canonical lifts $Z^{\eta,f}$ and the BPHZ character $g^\eta_\BPHZ$ introduced above coincides with the notation introduced in Section \ref{sec:ext:reg:structure}.
\end{remark}

We recall that, following arguments similar to \cite[Thm.~7.8]{Hairer2014}, see also \cite{ChandraHairer2016}, convergence of $\hat Z^{\eta^\eps}_\BPHZ=:(\hat \Pi^{\eta^\eps}, \hat \Gamma^{\eta^\eps})$ in $\CM_0^\rand$ can essentially be reduced to bounds of the form
\begin{align}\label{eq:bound:trees}
\E | (\hat \Pi_z^{\eta^\eps} \tau)(\phi^\lambda_z)|^2 \lesssim \lambda^{2|\tau|_\fs+\kappa}
\end{align}
uniformly in $\lambda \in (0,1)$, for some $\kappa>0$, any $\tau \in \TT$ of negative homogeneity and any $\phi \in \CC_c^\infty(\domain)$, compare \cite[Thm.~7.8]{Hairer2014}.
These moments can be conveniently represented as a finite sum over BPHZ-renormalized evaluations of graphs, each obtained via a "pairing" of the leaves of two disjoint copies of $\tau$, and the bound (\ref{eq:bound:trees}) follows from bounding each of these contractions separately. This was carried out in \cite{ChandraHairer2016} by applying a purely analytical BPHZ theorem for (hyper-)graphs. In the sequel we will need to work with a slightly different formulation of this analytical bound and in order to state it we introduce some notation. We begin with a simple lemma.

\begin{lemma}
Given $w,z \in \domain$ and a tree $\tau \in \TT$  there exists a unique locally integrable function $\Lambda_{z;w}\tau : \domain^{L(\tau)} \to \R$, smooth away from the big diagonal\footnote{The \emph{big diagonal} contains all $x \in \domain^{L(\tau)}$ such that there exists distinct $u,v \in L(\tau)$ such that $x_u = x_v$.} and away from $w,z$, symmetric under any permutation $\sigma$ of $L(\tau)$ with the property that $\ft \circ \sigma = \ft$, and such that one has
\[
\Pi^\eta_w \tau (z) = \int_{\domain^ {L(\tau)} } dx
\,
\Lambda_{z;w}\tau (x_{L(\tau)})
\prod_ {u \in L(\tau)} \eta_{\ft(u)}(x_{u}).
\]
for any $\eta \in \Omega_\infty$.
\end{lemma}
\begin{proof}
The proof is straightforward using induction over the number of edges of $\tau$, the fact that the canonical lift is an admissible model,  and the fact that $\Pi_w^\eta$ is multiplicative for the tree product.
\end{proof}

Furthermore, for any $f \in \Omega_\infty(\FL_-^\rand)$ and $\cov \in \FC(\FL_-^\det)$ we define the function $ \Lambda_{z;w}^{f, \cov} : \domain ^{ L ^\det (\tau ) } \to \R$ by setting
\[
 \Lambda_{z;w} ^{f,\cov} \tau ( x _{ L ^\det ( \tau ) } )
:=
\int_{\domain ^{ L ( \tau ) } } dy_{ L (\tau ) }
\;
	\Lambda_{z;w} \tau ( y _{ L ( \tau ) } )
	\prod_{ u \in L ^ \rand (\tau ) } f _{\ft(u)} (y_u)
	\prod_{ u \in L ^ \det ( \tau ) } \cov_{\ft(u)}(x_u - y_u).
\]
It follows that $\Lambda_{z;w} ^{f, \cov}$ is symmetric under any permutation $\sigma$ of $L ^\det (\tau)$ with the property that $\ft \circ \sigma = \ft$ on $L^\det ( \tau )$, so that we can naturally view the domain of definition of $\Lambda_{z;w} ^{f,\cov}$ as $\D ^ {\noisems^\det(\tau)}$, where $\noisems^\det(\tau)$ is the multiset given by $\noisems^\det(\tau):=[L^\det(\tau),\ft]$, compare the notation introduced in Section \ref{sec:notation}. We will switch frequently between these two pictures in the sequel.

Given a $\FL_-^\det$-valued multiset $\tm$ and a kernel $\cov \in \FC(\FL_-^\det)$, we let $\topcirc H_{\tm, \cov}$ be the space given by all functions $F \in \CC_c^\infty(\domain^\tm)$ which can be written in the form
\begin{align}\label{eq:F:barF}
F(x_\tm) := 
\int_{\domain^\tm} d\bar x_{\tm}
	G( \bar x_\tm)
	\prod_{u \in \tm} \cov_{\ft(u)}(x_u  - \bar x _u)
\end{align}
for some $G \in \CC_c^\infty(\domain^{\tm})$, endowed with the scalar product
\begin{align}\label{eq:F:barF:norm}
\dual{ F, \bar F }_{\tm, \cov} := 
\int _{ \domain 
	^{ \tm } \times \domain ^{ \tm   }  
} 
dx _{ \tm  } d \bar x _{ \tm  }
G( x_ \tm ) \bar G ( \bar x _ \tm)
\prod_{u \in \tm } \cov_{ \ft( u ) } ( x_u - \bar x_{u} ),
\end{align}
where $G$ and $\bar G$ are as in (\ref{eq:F:barF}) for $F$ and $\bar F$ respectively, and we write $H_{\tm,\cov}$ for the closure of $\topcirc H_{\tm,\cov}$ under the induced norm.
We also write $\CT^\ex[\tm] \ssq \CT^\ex$ for the linear subspace spanned by trees $\tau \in \TT^\ex$ with the property that one has $\tm^\det(\tau) = \tm$, and we note that for any $f \in \Omega_\infty(\FL_-^\rand)$ and any $w,z \in \domain$ one has
\[
 \Lambda^{f,\cov}_{z;w} : \CT^\ex[\tm] \to H_{\tm, \cov}.
\] 

Finally, it follows directly from the definition of the coproduct $\Delta_-^\ex$ and the character $g_\BPHZ^\eta$ that for any $\eta \in \SM_\infty( \FL_-^\rand )$ one has 
\[
M^{g_\BPHZ^\eta} = (g_\BPHZ^\eta \otimes \Id) \Delta_-^\ex : \CT^\ex[\tm] \to \CT^\ex[\tm]
\]
for any multiset $\tm$, so that it makes sense to define the random variable
\[
\Lambda^{\eta, \cov}_{z, \bar z; w} \tau
:=
\dual{ 
	 \Lambda_{z; w} ^{\eta , \cov} M^{g_\BPHZ^\eta} \tau 
\, , \,
	\Lambda_{\bar z;w} ^{\eta , \cov} M^{g_\BPHZ^\eta}  \tau } _{\tm, \cov}
\]
for any $\tau \in \CT^\ex[\tm]$.
The following theorem contains the key analytic bound on which the analysis below is bases on. 

\begin{theorem}\label{thm:BPHZ:graphs:from:trees}
For any $C>0$ and any compact $K \ssq \domain$ one has the bound
\begin{align}\label{eq:bound:E:Lambda}
\sup_{w \in K}
\E
\Big |
\int_{ \domain \times \domain } dz d\bar z 
\,
\Lambda^{\eta ,\cov}_{z, \bar z; w} \tau
\,
\phi^\lambda_w(z) \phi^\lambda_w(\bar z)
\Big | 
\lesssim 
\lambda^{2 |\tau|_\fs + \theta}
\end{align}
uniformly over $\lambda \in (0,1)$, $\cov \in \FC( \FL_- ^\det )$ with $\|\cov\|_{|\cdot|_\fs}\le C$, and $\eta \in \SM_\infty( \FL_- ^\rand )$ with $\|\eta\|_{|\cdot|_\fs}\le C$, for $\theta>0$ small enough and for any $\tau \in \TT$ of negative homogeneity and any $\phi \in \CC_c^\infty(\domain)$.
\end{theorem}

\begin{remark}
The reason for writing the above expression in this unusual way is that we are going to apply Cauchy-Schwarz estimates in the Hilbert space $H_{\fm,\cov}$.
\end{remark}

\begin{proof}
First note that we can replace $\cov_{\ft,\ft'}$ with a regularization $\cov_{\ft,\ft'}^\eps:=\cov_{\ft,\ft'}*\rho^{(\eps)}$ for some symmetric,  non-negative definite mollifier $\rho$ and then take the limit $\eps \to 0$, so that it suffices to show (\ref{eq:bound:E:Lambda}) for smooth kernels $\cov \in \FC(\FL_-^\det)$. Furthermore, by definition $\cov \in \FC(\FL^\det_-)$ is non-negative definite when viewed as an integral operator on $L^2(\domain)^{\FL_-^\det}$, and as a consequence there exists a (unique) Gaussian, centred, stationary noise $\bar \eta \in \SM_\infty(\FL_-^\det)$ with the property that 
\[
\E[\bar \eta_\ft(\phi) \bar \eta_{\ft'}(\psi)] = 
\delta_{\ft,\ft'}
\int_{\domain \times \domain} dx dy
\cov_{\ft}(x-y) \phi(x) \psi(y) 
\]
for any choice of test functions $\phi,\psi \in \CC_c^\infty(\domain)$ and any $\ft,\ft' \in \FL_-^\det$.
Enlarging the probability space $(\Omega,\P)$ if necessary, we can additionally assume that $\Omega=(\Omega_\rand\times \Omega_\det)$ and $\P=\P_\rand\otimes \P_\det$ for some probability measures $\P_\rand$ on $\Omega_\rand $ and $\P_\det$ on $\Omega_\det$, respectively, and such that $\eta$ respectively $\bar \eta$ is a collection of random fields on $(\Omega_\rand,\P_\rand)$ respectively $(\Omega_\det,\P_\det)$. In particular, one has that $\eta$ and $\bar \eta$ are independent, and we write $\xi := \eta \sqcup \bar \eta \in \SM_\infty(\FL_-)$. We denote as usual the BPHZ character for $\xi$ by $g^\xi_\BPHZ \in \CG_-^\ex$ and we denote the BPHZ renormalized canonical lift of $\xi$ to a model by $\hat Z_\BPHZ^\xi=(\hat \Pi^\xi, \hat \Gamma^\xi)$.

We first assume that the tree $\tau \in \TT$ has the property that for any noise type $\Xi \in \FL_-^\det$ there exists at most one $u \in L(\tau)$ such that $\ft(u) = \Xi$. We claim that in this case one has
\begin{align}\label{eq:identity:Lambda:Pi}
\E
\Big |
\int_{ \domain \times \domain } dz d\bar z 
\,
\Lambda^{\eta ,\cov}_{z, \bar z; w} \tau
\,
\phi^\lambda_w(z) \phi^\lambda_w(\bar z)
\Big | 
=
\E
\Big|
\int_{ \domain } dz
\hat \Pi^\xi_{w}(z)
\phi^\lambda_w(z) 
\Big|^2,
\end{align}
from which (\ref{eq:bound:E:Lambda}) follows from Theorem \ref{thm:BPHZ}. In order to see (\ref{eq:identity:Lambda:Pi}) note that from the assumption on $\tau$ and the fact that $\xi_\ft$ and $\xi_{\ft'}$ are independent for any $\ft \in \FL_-^\det$ and any $\ft' \in \FL_-^\rand$ it follows that $g^\xi_\BPHZ$ vanishes on subtrees $\sigma \ssq \tau$ with the property that $\sigma \notin \CT^\ex_-[\FL_-^\rand]$, which in turn implies that one has the identity
\[
M^{g_\BPHZ^\eta} \tau = M^{g_\BPHZ^\xi} \tau.
\]
A fortiori it follows that one has
\[
\Lambda^{\eta,\cov}_{z, \bar z; w}
=
\E^{\P_\det} |\hat \Pi^\xi_w(z) \hat \Pi^\xi_w(\bar z)|
\qquad\quad
\P_\rand-a.s.,
\]
and (\ref{eq:identity:Lambda:Pi}) follows.

In the general case, define for $\Xi \in \FL_-^\det$ the number $m(\Xi) \in \N\cup\{0\}$ as the number of noise-type edges $u \in L(\tau)$ with the property that $\ft(u) = \Xi$, and let $m := \max _{\Xi \in \FL_-^\det} m(\Xi)$. We then consider the $m$-fold extension $\extxt\CT^\ex[m]$ of $\CT^\ex$. 

We define $\extxt\FL_-^{[m],\rand}:= \FL_-^\rand \times [m]^\star$ and $\extxt \FL_-^{[m],\det}:= \FL_-^\det \times [m]^\star$, and we define an element $\tilde\cov \in \FC(\extxt\FL_-^{[m],\det})$ by setting
\[
\tilde\cov_{\Xi,\Xi'} := \tilde\cov_{\Xi,(\Xi',k')}:=0
\qquad\text{ and }\qquad
\tilde\cov_{(\Xi,k),(\Xi',k')}
:=
\delta_{k.k'} \cov_{\Xi,\Xi'}
\]
for any $\Xi,\Xi' \in \FL_-$ and any $k,k'\in [m]$.
Moreover, we define $\Phi$ as the set of all type maps $\tilde\ft : E(\tau) \to \extxt\FL_-^{[m]}$ such that $\tilde\ft(e) = \ft(e)$ for any kernel-type edge $e \in K(\tau)$ and such that the following holds.
\begin{itemize}
\item For any noise-type edge $u \in L^\rand(\tau)$ one has $\tilde \ft(u) = \ft(u)$.
\item For any noise-type edge $u \in L^\det(\tau)$ one has $\tilde\ft (u) = (\ft(u), k(u))$ for some $k(u) \in [m]$.
\item For any noise-type $\Xi \in \FL_-^\det$ the map $k$ restricted to $L_\Xi(\tau):=\{u \in L(\tau) : \ft(u) = \Xi \}$ is a bijection from $L_\Xi(\tau)$ onto $[m(\Xi)]$.
\end{itemize}
We note the following consequences of this definition: For any $\tilde\ft \in \Phi$ one has $(\tau,\tilde\ft) \in \extxt\CT^\ex[m]$ and $(\tau,\tilde\ft)$ satisfies the assumption of the first part of the proof. It remains to apply the results of the first part to any of the trees $(\tau,\tilde\ft)$ individually, note that one has
\[
 \Lambda_{z; w}^{\eta,\cov} \tau (x_{L^\det(\tau)})
=
\frac{1}{S}
\sum_{\tilde \ft \in \Phi}
	\Lambda_{z;w} ^{\eta,\cov} (\tau,\tilde\ft) (x_{L^\det(\tau)}),
\]
and use a Cauchy-Schwarz estimate.
Here $S$ is a symmetry factor given by
$
S = \prod_{\Xi \in \FL_-^\det }{m(\Xi)!}.
$
\end{proof}

In order to continue, we first note a relation between the norms of the spaces ${ H[\FL_-^\det, \cov]}$ defined in (\ref{eq:HcovNorm}) and $\|\cdot\|_{\fm,\cov}$ defined in (\ref{eq:F:barF:norm}). 
\begin{lemma}\label{lem:fh}
Let $\tm$ be any multi-set with values in $\FL_-^\det$, let $\cov \in \FC(\FL_-^\det)$, let $h \in H[\FL_-^\det, \cov]$, and  
define $\fh \in H_{\tm,\cov}$ by setting
\[
\fh := \bigotimes_{u \in \tm} h_{\ft(u)}.
\]
Then one has
\[
\|\fh\|_{\tm,\cov} \lesssim (\|h\|_{H[\FL_-^\det,\cov]})^{\#\tm}.
\]
\end{lemma}
\begin{proof}
By definition (\ref{eq:F:barF:norm}), it suffices to show the statement for $\# \tm = 1$. In this case write $\ft \in \FL_-^\det$ for the type such that one has $\tm=\{\ft\}$ and note that, writing $\pi_\ft:H[\FL_-^\det,\cov] \to H[\FL_-^\det,\cov]$ for the projection given by $(\pi_\ft(k))_{\ft'} := k_{\ft} \delta_{\ft,\ft'}$ for any $k \in H[\FL_-^\det,\cov]$, we have the equality $\|\fh\|_{\tm,\cov} = \|\pi_\ft h\|_{H[\FL_-^\det,\cov]}$. The result now follows from the fact that the $\pi_\ft$ is a continuous projection, since both kernel and range of $\pi_\ft$ are closed in $H[\FL_-^\det,\cov]$.
\end{proof}

Recall that for $\tau \in \CT^\ex$ we defined the multiset $\tm^\det(\tau)=[L^\det(\tau),\ft]$, and with this notation we define for any $h \in H[\FL_-^\det,\cov]$ the quantities
\begin{align}\label{eq:fh}
\fh^h_\tau := \bigotimes_{u \in \tm^\det(\tau)} h_{\ft(u)} 
\qquad\text{ and }\qquad
[h]_\tau := \|\fh^h_\tau\|_{\tm^\det(\tau),\cov}.
\end{align}
We stress that $[\cdot]_\tau$ fails to be a semi-norm unless $\#\tm^\det(\tau)=1$.
Our ultimate goal is to show that $(\eta,f) \mapsto \hat Z^{\eta,f}_\BPHZ$ extends continuously to $\eta \in \SM_0(\FL_-^\rand)$ and $f \in H[\FL_-^\det,\cov]$ for any kernel $\cov \in \FC(\FL_-^\det)$.
As a preparation for this statement, we show the following result.

\begin{proposition}\label{prop:canonical:lift:cameron:martin}
Let $\CT^\ex$ be a regularity structure constructed as in Section \ref{sec:reg:structures} satisfying Assumption \ref{ass:main:reg}, and let $\FL_-=\FL_-^\det$ and $\FL_-^\rand = \emptyset$. Let moreover $\cov \in \FC(\FL_-)$ be a kernel. Then the map $\Omega_\infty(\FL_-) \to \CM_\infty(\FL_-)$, $h \mapsto Z^h$, extends to a continuous map from $H[\FL_-,\cov]$ into $\CM_0(\FL_-)$ which is locally Lipschitz continuous in the sense that for any $\gamma \in \R$, any compact $K \ssq \domain$ and any $R>0$ one has
\begin{align}\label{eq:canonical:lift}
\fancynorm{Z^h; Z^k}_{\gamma,K} \lesssim \|h-k\|_{H[\FL_-,\cov]}
\end{align}
uniformly over all $h,k \in H[\FL_-,\cov]$ with $\|h\|_{H[\FL_-,\cov]} \lor \| k \|_{H[\FL_-,\cov]}<R$.
\end{proposition}
\begin{proof}
For any fixed $\gamma>0$ and compact $K \ssq \domain$ the pseudo metric $\fancynorm{\cdot;\cdot}_{\gamma,K}$ induces a complete metric space $\CM_0(K)$ via metric identification, so that  it is sufficient to show that one has the bound $(\ref{eq:canonical:lift})$ for any $h,k \in \topcirc H[\FL_-,\CC]$; note for this that $\topcirc H[\FL_-,\CC]$ is dense in $H[\FL_-,\CC]$, and thus any such local Lipschitz map has a unique extension to $H[\FL_-,\CC]$ which is again locally Lipschitz, and this concludes the proof. Following arguments identical to \cite[Thm. 10.7]{Hairer2014}, it suffices to show that for any tree $\tau \in \CT$ of negative homogeneity and any $\phi \in \CC_c^\infty(\domain)$ with $\supp \phi \ssq B_1(0)$ there exists $\theta>0$ such that one has the bound
\begin{align}\label{eq:canonical:lift:2}
|(\Pi^{h}_z \tau - \Pi^k_z \tau)(\phi_z^\lambda)| \lesssim \lambda ^{|\tau|_\fs+\theta} \|h-k\|_{H[\FL_-,\cov]} 
\end{align}
uniformly over $h,k \in \topcirc H[\FL_-,\cov]$ with $\|h\|_{H[\FL_-,\cov]} \lor \|k\|_{H[\FL_-,\cov]}\le R$, $z \in K$ and $\lambda \in (0,1)$. 

We first show that one has the bound
\begin{align}\label{eq:canonical:lift:3}
|\Pi^{h}_z \tau(\phi_z^\lambda)| \lesssim \lambda ^{|\tau|_\fs+\theta} 
[h]_\tau.
\end{align}
uniformly over $h \in \topcirc H[\FL_-,\cov]$ with $\| h \|_{H[\FL_-,\cov]} \le R$, $z \in K$ and $\lambda \in (0,1)$. 
For this we use the identity 
\begin{align} \label{eq:canonical:lift:4}
\Pi_z^h \tau(y) = 
\dual{
	\treekpos_{y;z}^{0, \cov}\tau , \fh^h_\tau
}_{\tm^\det(\tau), \cov}
\end{align}
from which a Cauchy-Schwarz estimate on the Hilbert space $H_{\noisems^\det(\tau), \cov}$ shows that the left hand side of (\ref{eq:canonical:lift:3}) can be estimated by
\begin{align*}
\|\fh^h_\tau\|_{\noisems(\tau), \cov}
\,
\|
F_z
\|_{ \noisems(\tau), \cov },
\end{align*}
where $F_z(x_{L(\tau)} ): = 
\int_{\domain}dy \,
	\treekpos_{y;z} \tau ( x_{L(\tau)})	\phi_z^\lambda(y) $.
Comparing the second term in this expression with Theorem~\ref{thm:BPHZ:graphs:from:trees}, the estimate (\ref{eq:canonical:lift:3}) follows. 

The bound (\ref{eq:canonical:lift:2}) is now an almost immediate consequence of (\ref{eq:canonical:lift:3}) applied to extended regularity structure and Lemma~\ref{lem:diff:Pi} applied for $g = \bold 1^*$. Indeed, first note that one has for any typed tree $(\tau,\ft) \in \CT^\ex$ the identity \footnote{Recall the notation $\ext\FL_- = \FL_- \sqcup \hat\FL_-$.}
\[
\Pi_z^h (\tau, \ft) = \ext \Pi _z ^{0,h} (\tau, \hat \ft),
\]
where on the right hand side we denote the canonical lift of $(0,h)$ to a model for the onefold extension $\ext \CT^\ex$ of $\CT^\ex$, and we write $\hat \ft := \ft$ for any $\ft \in \FL_+$. By Lemma \ref{lem:diff:Pi} it follows that one has
\[
\Pi_z^h (\tau, \ft) - \Pi_z^k (\tau, \ft)
=
\extxt\Pi_z^{0,h,k,h-k}  \AA (\tau,\hat \ft).
\]
Applying (\ref{eq:canonical:lift:3}) to each tree on the right hand side of this identity, we obtain the desired bound. Note for this that Lemma \ref{lem:fh} implies in particular that $[h]_\tau \lesssim \|h-k\|_{H[\FL_-,\cov]}$ uniformly over all $h,k \in H[\FL_-,\cov]$ with $\|h\|_{H[\FL_-,\cov]} \lor \| k \|_{H[\FL_-,\cov]}<R$.
\end{proof}

The main stochastic ingredient for the proof of Theorem \ref{thm:extension:operator} below is the following bound, for which we introduce the notation  $\hat Z^{\eta,h}_\BPHZ = (\hat \Pi^{\eta,h}, \hat \Gamma^{\eta,h})$
for any $\eta \in \SM_\infty(\FL_-^\rand)$.
We then have the following.
\begin{proposition}\label{prop:CM:BPHZ:estimate:}
Let $\CT^\ex$ be a regularity structure constructed as in Section \ref{sec:reg:structures} satisfying Assumption \ref{ass:main:reg}, and assume that we are given the decomposition $\FL_- = \FL_-^\det \sqcup \FL_-^\rand$, and let $\cov \in \FC(\FL_-^\det)$. 
 Then there exists $\kappa>0$ such that one has for any compact $K \ssq \domain$, any $p\ge 1$, any $C>0$, and any $\tau \in \CT^\ex$, and any test function $\phi \in \CC_c^\infty(\domain)$ the bound
\begin{align}\label{eq:ext:operator:2}
\sup_{z \in K}
\E 
\Big[
\sup_{ h \in H[\FL_-^\det, \cov] }
\frac{1}{[h]_\tau}
|(\hat\Pi^{\eta,h}_z \tau) (\phi^\lambda_z)|
\Big]
^{2p} 
\lesssim 
\lambda^{2p|\tau|_\fs + 2p\kappa}
\end{align}
uniformly over $\lambda \in (0,1)$ and $\eta \in \SM_\infty(\FL_-^\rand)$ with $\|\eta\|_{|\cdot|_\fs} \le C$.
\end{proposition}
\begin{proof}
First note that $\hat Z^{\eta,h}_\BPHZ$ is almost surely well defined for any $\eta \in \SM_\infty(\FL_-^\rand)$  and any $h \in H[\FL_-^\det,\cov]$ by Proposition \ref{prop:canonical:lift:cameron:martin}. Furthermore, using Gaussian hypercontractivity it suffices to show this proposition for $p=1$.  

We now write $\fh^h_\tau \in H_{\tm^\det(\tau),\cov}$ for the function in (\ref{eq:fh}),  
and we use the identity
\[
\hat \Pi_w^{\eta, h}\tau(z) = 
	\Pi^{\eta, h}_w M^{g_\BPHZ^\eta} \tau (z) =
	\dual{
		\hat \Lambda^{\eta,\cov}_{z; w} M^{g_\BPHZ^\eta} \tau 
		\, , \,
		\fh^h_\tau
	}_{\tm^\det(\tau), \cov}
\]
so that by Cauchy-Schwarz we obtain the estimate
\[
|(\hat\Pi^{\eta,h}_z \tau) (\phi^\lambda_z)|
\le
\Big |
\int_{ \domain \times \domain } dz d\bar z 
\,
\Lambda^{\eta ,\cov}_{z, \bar z; w} \tau
\,
\phi^\lambda_w(z) \phi^\lambda_w(\bar z)
\Big | ^{\frac{1}{2} }
\|\fh^h_\tau\|_{\tm,\cov}
\]

The result now follows from Theorem \ref{thm:BPHZ:graphs:from:trees}.
\end{proof}

We demonstrate the idea of the proof of Proposition \ref{prop:CM:BPHZ:estimate:} with the following example.
\begin{example}
We consider a situation similar to Example \ref{ex:introduce:tree}, but this time we assume we have two noise types $\FL_-=\{\starOneNoiseShort,\,\starOneNoiseShiftShort\}$ with $\starOneNoiseShort \in \FL_-^\rand$ and \,$\starOneNoiseShiftShort \in \FL_-^\det$, and we write $H_{\,\starOneNoiseShiftShort}:=H[\{\,\starOneNoiseShiftShort\,\},\cov_{\starOneNoiseShiftShort }]$. We want to derive in detail the bound (\ref{eq:ext:operator:2}) for $\treeExampleSmallShiftA$ and $p=1$. First, we can write
\[
\E[ \sup_{h} \frac{1}{\| h_{\,\starOneNoiseShiftShort}\|_{H_{\,\starOneNoiseShiftShort}}} | \hat\Pi^{\eta,h}_z \treeExampleSmallShiftA (\phi^\lambda_z)| ]^2
=
\E \| \hat\Pi^{\eta, \delta_y}_z \treeExampleSmallShiftA (\phi^\lambda_z) \|_{H_{\,\starOneNoiseShiftShort}} ^2
\]
where on the right hand side we take the $\|\cdot\|_{H_{\,\starOneNoiseShiftShort}}$-norm of the map $y \mapsto \hat\Pi^{\eta, \delta_y}_z \treeExampleSmallShiftA (\phi^\lambda_z)$, where $\delta_y$ denotes the $\delta$-distribution centred around $y$. This is a slight abuse of notation, explicitly this expression is equal to
\[
\hat\Pi^{\eta, \delta_y}_z \treeExampleSmallShiftA (\phi^\lambda_z)
=
\int dx \,\eta(x) K(x-y) \phi^\lambda_z(y).
\]
Now, let $\tilde\eta \in \SM_0$ be a Gaussian noise independent of $\eta$ and such that $\E[\tilde\eta(x), \tilde\eta(y)]=\cov_{\,\starOneNoiseShiftShort}(x-y)$. It then follows from the definition of the Hilbert space $H_{\,\starOneNoiseShiftShort}$ that we have the identity
\[
\E \| \hat\Pi^{\eta, \delta_y}_z \treeExampleSmallShiftA (\phi^\lambda_z) \|_{H_{\,\starOneNoiseShiftShort}} ^2
=
\E
|\hat\Pi^{\eta, \tilde\eta}_z \treeExampleSmallShiftA (\phi^\lambda_z)|^2.
\]
From this expression we obtain the bound (\ref{eq:ext:operator:2}) from \cite{ChandraHairer2016}.
\end{example}

The main result of this section is the following theorem.

\begin{theorem}\label{thm:extension:operator}
Assume we are given the decomposition $\FL_- = \FL_-^\det \sqcup \FL_-^\rand$, and let $\cov \in \FC(\FL_-^\det)$. Then the map $\SM_\infty(\FL_-^\rand) \times \Omega_\infty(\FL_-^\det) \to \CM_\infty^\rand$, $(\eta,h) \mapsto \hat Z^{\eta, h}_\BPHZ$ extends to a continuous map from $\SM_0(\FL_-^\rand) \times H[\FL_-^\det, \cov]$ into $\CM_0^\rand$.

Moreover, for any $\eta \in \SM_0(\FL_-^\rand)$ there exists a null set $N$ such that for any $\omega \in N^c$ one has that the map $h \mapsto \hat Z^{\eta, h}_\BPHZ(\omega)$ is locally Lipschitz continuous in the sense for any $\gamma \in \R$, any compact $K \ssq \domain$ and any $R>0$ one has
\begin{align}\label{eq:extension:lipschitz}
\fancynorm{ \hat Z^{\eta,h}_\BPHZ(\omega) ; \hat Z^{\eta,k}_\BPHZ(\omega) }_{\gamma,K}
\lesssim
\|h-k\|_{H[\FL_-^\det, \cov]}
\end{align}
uniformly over all $h,k \in H[\FL_-^\det,\cov]$ with $\|h\|_{H[\FL_-^\det,\cov]} \lor \|k \|_{H[\FL_-^\det,\cov]}<R$. 

Finally, given an approximation $\eta^\eps \in \SM_\infty(\FL_-^\rand)$ of $\eta$ there exists a subsequence $\eps \to 0$ and a null set $N$ with the property that $\eta^\eps(\omega) \to \eta(\omega)$ and $\hat Z^{\eta^\eps,h}_\BPHZ(\omega) \to \hat Z^{\eta,h}_\BPHZ(\omega)$ for any $h \in H[\FL_-^\det,\cov]$ and $\omega \in N^c$.
\end{theorem}

\begin{remark}
At this stage we point out a key difference between the present approach and \cite{CannizzaroFrizGassiat2015}. In the latter the authors obtain a stronger statement as they construct a \emph{deterministic} continuous extension operator which is such that it maps $(\hat Z^\eta_\BPHZ(\omega),h)$ onto $\hat Z^{\eta,h}_\BPHZ(\omega)$ for almost every fixed $\omega$. This however comes at the expense of analytical difficulties that could only be carried out in a very special case. In the current paper in contrast, the analytic difficulties are bypassed by constructing the extension \emph{stochastically}, which in particular allows us to use the results of \cite{ChandraHairer2016} to show the necessary analytic estimates. This comes at the expense of a somewhat weaker statement but has the advantage of being immediately completely general.
\end{remark}

\begin{proof}
We first note that as a consequence of Proposition \ref{prop:CM:BPHZ:estimate:} for any $\eta \in \SM_0(\FL_-^\rand)$ and $h \in H[\FL_-^\det, \cov]$ there exists a random variable $\hat Z^{\eta,h}_{\BPHZ}$ taking values in ${\CM_0}$ such that for some $\theta>0$ one has
\begin{align}\label{eq:ext:operator:1}
\E \fancynorm{{\hat Z}^{\eta^\eps, h}_{\BPHZ}}_{\gamma, K}^p \lesssim 1 
\quad\text{ and }\quad
\E \fancynorm{{\hat Z}^{\eta^\eps, h}_{\BPHZ}; {\hat Z}^{\eta,h}_{\BPHZ}}_{\gamma, K}^p 
\lesssim \eps^{\theta p}
\end{align}
uniformly in $\eps>0$, where $\eta^\eps \in \SM_\infty(\FL_-^\rand)$ is an approximation of $\eta$. 
To see this, note that the first bound above is a consequence of Proposition \ref{prop:CM:BPHZ:estimate:} using an argument identical to the one given in \cite[Thm.~10.7]{Hairer2014}.
The existence of the extension $\hat Z^{\eta, h}_{\BPHZ}$ and the second bound in (\ref{eq:ext:operator:1}) follow now along the usual lines, see for instance \cite[Sec.~10.4,10.5]{Hairer2014} or the proof of \cite[Thm.~2.34]{ChandraHairer2016}.

To see (\ref{eq:extension:lipschitz}), let $\eta \in \SM_0(\FL_-^\rand)$ and let $\eta^\eps \in \SM_\infty(\FL_-^\rand)$ be an approximation of $\eta$. Let moreover $\tilde H \ssq H[\FL_-^\det,\cov]$ be a dense and countable subset. It suffices to fix a subsequence $\eps\to 0$ such that ${\hat Z}^{\eta^\eps, h}_{\BPHZ} \to {\hat Z}^{\eta, h}_{\BPHZ}$ almost surely for any $h \in \tilde H$ as $\eps \to 0$, and to show that (\ref{eq:extension:lipschitz}) holds almost surely with $\eta$ replaced by $\eta^\eps$ uniformly over $\eps > 0$ and $h,k \in \tilde H$ with $\|h\|_{} \lor \|k\|\le R$, where we write $\|h\|:=\|h\|_{H[\FL_-^\det,\cov]}$ from now on in order to simplify notation. Upon choosing a sub-subsequence, this follows from the estimate
\begin{align}\label{eq:extension:lipschitz:2}
\E\Big[
 \supss{h,k \in \tilde H \\ \|h\| \lor \|k\| \le R}
\|h-k\|_{H[\FL_-^\det, \cov]}^{-1}
\fancynorm{ \hat Z^{\eta^\eps,h}_\BPHZ ; \hat Z^{\eta^\eps,k}_\BPHZ }_{\gamma,K}
\Big]
\lesssim
1
\end{align}
uniformly over $\eps>0$. Using again an argument identical to the one given in \cite[Thm.~10.7]{Hairer2014} we see that (\ref{eq:extension:lipschitz:2}) is implied once we show for any fixed $\phi \in \CC_c^r(\domain)$ and any $\tau \in \CT^\ex$ the bound
\begin{align}\label{eq:ext:operator:estimate}
\sup_{x \in K}
\E
\Big[
\supss{h,k \in \tilde H \\ \|h\| \lor \|k\| \le R}
\| h - k \|_{H[\FL_-^\det,\cov]}^{-2p}
|	(\hat \Pi_x^{\eta^\eps,h} \tau - 
	\hat \Pi_x^{\eta^\eps,k} \tau)
	(\phi_x^\lambda)|^{2p}
\Big]
 \lesssim
 \lambda^{2p|\tau|_\fs+ 2p\theta}
\end{align}
uniformly over $\lambda \in (0,1)$ for any $p\ge 1$.

We will show this by applying (\ref{eq:ext:operator:2}) to the extended regularity structure constructed in Lemma \ref{lem:diff:Pi}. For this consider the set $J := \{ \ch, \ck, \ch - \ck \}$ and we denote by $\extxt\FL_-^J = \FL_- \sqcup \hat \FL_-^J$ the corresponding extended set of noise-types. The set $\hat \FL_-^J$ comes with a natural decomposition into $\hat \FL_-^{J,\det}:= \FL_-^\det \times J$ and $\hat \FL_-^{J,\rand} := \FL_-^\rand \times J$. 

We define $\extxt\eta^\eps \in \SM_\infty (\FL_-)$ and, given $h \in \tilde H$, we define $\extxt{h} \in \Omega_0( \FL_-\times\{\ch\} )$ by setting 
\begin{align*}
\extxt\eta^\eps _\ft:= 
\begin{cases}
	\eta^\eps_\ft \quad	& \text{ if }\ft \in \FL_-^\rand \\
	0					& \text{ if }\ft \in \FL_-^\det
\end{cases}
\,,\qquad\quad
\extxt h _{(\ft,\ch)}:= 
\begin{cases}
	0					& \text{ if }\ft \in \FL_-^\rand \\
	h_\ft 		\quad	& \text{ if }\ft \in \FL_-^\det
\end{cases}
\end{align*}
Now we note that one has for any $(\tau,\ft) \in \CT^\ex$ the identity $\hat \Pi_x^{(\eta^\eps,h)} (\tau,\ft) = \extxt{ \hat \Pi}_x^{\extxt \eta^\eps, \extxt h} (\tau, \tilde \ft)$, where we set $\tilde \ft := \ft$ if $\ft \in \FL_-^\rand \sqcup \FL_+$ and $\tilde \ft := (\ft, \ch)$ if $\ft \in \FL_-^\det$. By Lemma \ref{lem:diff:Pi}  it follows that one has the identity
\[
(\hat \Pi_x^{\eta^\eps,h} \tau - 
	\hat \Pi_x^{\eta^\eps,k} \tau)
=
(\hat \Pi^{\eta^\eps, h, k, h-k}_x \AA \tau),
\]
and by definition of $\AA$ we conclude that (\ref{eq:ext:operator:2}) implies (\ref{eq:ext:operator:estimate}).

The fact that the null set $N$ such that $\hat Z^{\eta^\eps,h}_\BPHZ(\omega) \to \hat Z^{\eta,h}_\BPHZ(\omega)$ for $\omega \in N^c$ can be chosen independently of $h \in H[\FL_-^\det,\cov]$ is now a direct consequence of the uniform bound of the local Lipschitz constants.
\end{proof}

The main application of the previous theorem is the situation in which the regularity structure is itself a onefold extension with noise types given by $\FL_- \sqcup \hat \FL_-$, and one has that $\FL^\rand_- = \FL_-$ and $\FL^\det_- = \hat\FL_-$.
This leads to the following immediate corollary.

\begin{corollary}\label{cor:extension:operator}
Let $\eta \in \SM_0(\FL_-)$ and let $\eta^\eps \in \SM_\infty(\FL_-)$ be an approximation of $\xi$.
Then, there exists a subsequence $\eps \to 0$ and a null set $N \ssq \Omega$ such that the following holds. For any $\omega \in N^c$ and $h\in H^\eta$  one has that ${\ext {\hat Z}}_\BPHZ^{\eta^\eps,h}(\omega)$ converges to ${\ext {\hat Z}}_\BPHZ^{\eta,h}(\omega)$ as $\eps \to 0$, and one has that for any $\omega \in N^c$ the map $h\mapsto \ext {\hat Z}_\BPHZ^{\eta,h}(\omega)$ is locally Lipschitz continuous from $H^{\eta}$ into $\ext{\CM}_0$.
\end{corollary}

\subsection{Shift of Models}\label{sec:shifts:models}

We recall the map $\cq: \ext\FL \to \FL$ from Section \ref{sec:ext:reg:structure} defined as the identity on $\FL$ and mapping $\hat \Xi$ onto $\Xi$ for any $\Xi \in \FL_-$. We extend this map to a projection $\cq: \ext{\CT}^\ex \to \CT^\ex$ by defining for any typed tree $\tau = (\tau,\ft )\in\shift\CT^\ex$ with type map $\ft:E(T) \to \ext\FL$  the tree 
\[
\cq(\tau,\ft):=(\tau,\cq\ft),
\]
and extending $\cq$ linearly to all of $\ext\CT^\ex$.
An important role is then played by the following operator, which generalizes the operator defined in (\ref{eq:shift:operator:A}). 
\begin{definition}
For $\tau\in\CT^{\ex}$ we denote by $\SS[\tau]\subseteq\shift\CT^\ex$ the set of trees $\sigma\in\shift\CT^\ex$ such that $\cq\sigma=\tau$. For any $r \in \R$ we define the linear operator $\SS_r: \CT^\ex \to \shift\CT^\ex$ by setting for any tree $\tau\in\TT^\ex$
\[
\SS_r \tau:= \sum_{\sigma\in\SS[\tau]} r^{m(\sigma)}\sigma.
\]
where $m(\sigma):= \# \hat L(\sigma)$. We call $\SS_r$ the \emph{shift operator}, and we write $\SS_r[\CT^\ex] \ssq \ext{\CT}^\ex$ for the image of the shift operator.
\end{definition}

\begin{example}
On the tree $\treeExampleSmallShiftA$ from Example \ref{ex:introduce:tree}, writing $\Xi = \starOneNoiseShort$ and $\hat\Xi= \starOneNoiseShiftShort$, we obtain the formula
$\SS_r \treeExampleSmall = \treeExampleSmall + r\,\treeExampleSmallShiftA + r\,\treeExampleSmallShiftB + r^2\,\treeExampleSmallShiftC$.
\end{example}
We extend the shift operator linearly and multiplicatively to the algebras $\hat\CT_-^\ex$ and $\hat\CT_+^\ex$, as well as to the Hopf algebras $\CT_-^\ex$ and $\CT_+^\ex$.
The following is a simple consequence of the definition.

\begin{lemma}\label{lem:properties:shift:operator}
For any $r \in \R$ the shift operator $\SS_r$ commutes with the action of both $\Delta_-^\ex$ and $\Delta_+^\ex$ on $\CT^\ex$. In particular, its image $\SS_r[\CT^\ex]$ forms a sector in $\ext{\CT}^\ex$ and $\SS_r$ commutes with the action of renormalization $M^g$ for any $g \in \CG_-^\ex$. Moreover, $\SS_r$ maps $\CT$ into $\ext\CT$, is multiplicative under the tree product, and commutes with the operation of compositions with smooth functions, integration and differentiation.
\end{lemma}
\begin{proof}
The fact that $\SS_r$ commutes with the action of the co-products is tedious to verify, but straightforward.

To see that $\SS_r[\CT^\ex]$ forms a sector, let $\Gamma=(\Id \otimes \gamma)\Delta^+_\ex$ for some character $\gamma \in \ext{\CG}_+^\ex$ as in \cite[(6.12)]{BrunedHairerZambotti2016}. Then one has
\[
\Gamma \SS_r \tau = (\SS_r \otimes \gamma \SS_r) \Delta_\ex^+ \tau \in \SS_r[\CT^\ex]
\]
for any $\tau \in \CT^\ex$. Similarly, for any $g \in \CG_-^\ex$ one has
\[
\SS_r M^g = (g \otimes \SS_r) \Delta_\ex^- = (g \otimes \Id) (\SS_r \otimes \SS_r) \Delta_\ex^- = M^g \SS_r.
\]
Here we used the fact that by definition of the embedding $\CG_-^\ex \hookrightarrow \ext{\CG}_-^\ex$ one has that  $g = g \SS_r $ for any $g \in \CG_-^\ex$.

The remaining statements of the lemma are a simple consequence of the definitions. 
\end{proof}

It follows directly from the definition that the operator $\SS_r$ is one to one, so that we can define its inverse $\SS_r^{-1}$ on the sector $\SS_r[\CT^\ex]$.
We now denote by\footnote{We write $\CL(X,Y)$ for the space of continuous linear maps from $X$ to $Y$.} $\SS^*_r: \CL(\shift\CT^\ex,\CC^\infty(\domain)) \to \CL(\CT^\ex,\CC^\infty(\domain))$
and $\SS^*_r:\shift\CM_\infty\to\CM_\infty$ the maps given respectively by
\[
(\SS^*_r \bold\Pi)\tau:=\bold\Pi\SS_r\tau
\qquad\text{ and }\qquad
\SS^*_r(\CZ(\bold\Pi)):=\CZ(\SS^*_r\bold\Pi).
\]

\begin{lemma}
The map $\SS^*_r : \ext{\CM}_\infty \to \CM_\infty$ is well defined and extends to a locally Lipschitz map from $\ext{\CM}_0$ onto $\CM_0$. Moreover, writing $\ext Z=(\ext \Pi,\ext \Gamma ) \in \ext{\CM}_0$ and $\SS^*_r(\ext Z)=(\Pi,\Gamma)\in\CM_0$, one has the identities
\begin{align}\label{eq:shift:operator:adjoint}
\Pi_x \tau = \ext\Pi_x \SS_r \tau
\quad\text{ and }\quad
\Gamma_{xy} \tau = \SS_r^{-1} \ext\Gamma_{xy} \SS_r
\end{align}
for any $x,y \in \domain$.
\end{lemma}
\begin{proof}
First note the right hand side of second identity in (\ref{eq:shift:operator:adjoint}) makes sense, since by Lemma \ref{lem:properties:shift:operator} the image $\SS_r[\CT]$ forms a sector, so that it is invariant under $\ext \Gamma_{x,y}$.

In order to see that the map $\SS^*_r : \ext{\CM}_\infty \to \CM_\infty$ is well defined, note first that it is clear from the definition that $\SS_r^*$ maps admissible  $\bold \Pi$ onto admissible $\SS_r^* \bold\Pi$, so that it remains to show the required analytic bounds. Noting that the map $\SS_r$ leaves homogeneities of trees invariant, these analytic bounds follow once we show that the identity (\ref{eq:shift:operator:adjoint}) holds for any $\ext \Pi$ and $\ext \Gamma$ given by the expressions \cite[(6.11)]{BrunedHairerZambotti2016} and \cite[(6.12)]{BrunedHairerZambotti2016}, respectively, but with $\bold \Pi$ replaced by $\SS_r^* \bold \Pi$. But as a direct consequence of Lemma \ref{lem:properties:shift:operator} it follows that $\SS_r$ commutes with the positive twisted antipode $\tilde \CA_+^\ex$, so that it follows that one has $f_z(\SS_r^* \bold \Pi)= f_z(\bold \Pi)\SS_r$, where $f_z$ is as in \cite[(6.11)]{BrunedHairerZambotti2016}. Plugging this into the respective formulae for $\ext \Pi$ and $\ext \Gamma$, the identities in (\ref{eq:shift:operator:adjoint}) follow.

Finally, the fact that $\SS_r^*$ extends to a locally Lipschitz map on $\ext{\CM}_0$ follows straight forwardly from the identity (\ref{eq:shift:operator:adjoint}) and the definition of the metric in the space of models \cite[(2.17)]{Hairer2014}.
\end{proof}

With this notation we have the following Theorem.

\begin{theorem}\label{thm:shift:operator}
Let $N\subseteq \Omega$ be the null set constructed in Corollary \ref{cor:extension:operator}. Then one has for any $\omega \in N^c$, any $r\in \R$, and any $h \in H^\xi$ the identity\footnote{Recall our convention that $\xi$ agrees with the identity on some Gaussian measure space almost surely, i.e. one has $\xi(\omega)=\omega$.}
\begin{align}\label{eq:shift:operator}
 {\hat Z}_{\BPHZ}^\xi(\omega + r h)
=
\SS^*_r(\ext {\hat Z}_\BPHZ^{\xi,h}(\omega)).
\end{align}
As a consequence, there exists a fixed subsequence $\eps \to 0$ (independent of $\omega$ and $h$) such that one has $\hat Z_\BPHZ^{\xi^\eps} (\omega + h)\to \hat Z_\BPHZ^{\xi}(\omega + h)$ for any $\omega \in N^c$ and $h \in H^\xi$.
Finally, for any $\omega \in N^c$ the maps $H^\xi \to \CM_0$ and $H^\xi \times H^\xi \to \ext \CM_0$ given by $h \mapsto \hat Z_\BPHZ^\xi(\omega + h)$ and $(h,k) \mapsto \hat Z_\BPHZ^{\xi,h}(\omega+k)$ are locally Lipschitz continuous.
\end{theorem}

\begin{proof}
By Theorem \ref{thm:extension:operator} it suffices to show (\ref{eq:shift:operator}) with $\xi$ replaced by $\xi^\eps$ for any $\eps>0$. For this in turn it is sufficient to show the identity $\SS^*_r (R^g \shift Z^{f,h})=R^g Z^{f + rh}$ for any $r\in \R$, any $g \in \CG_-$ and any $f \in \ext{\Omega}_\infty$, which is a simple consequence of Lemma \ref{lem:properties:shift:operator}. The rest of Theorem \ref{thm:shift:operator} is then a direct consequence of Theorem \ref{thm:extension:operator} and the fact that $\SS_r^*$ is locally Lipschitz.
\end{proof}

\subsection{Lifts of Abstract Fixed Point Problems}
\label{sec:abstract:fixed:point}

We are going to describe a class of abstract fixed point problems on the spaces $\CD^{\gamma,\eta}_V$ that we are going to look at in the sequel.

Let $\sector \ssq \CT$ be a sector spanned by a set of trees $\hat\sector$. Then the space $\ext \sector \ssq \ext \CT$ spanned by all trees $\sigma \in \SS[\tau]$ for some $\tau \in \hat\sector$ forms a sector in $\ext \CT$.
More generally assume we are given for any $\ft \in \FL_+$ a sector $\sector_\ft$ in $\CT$ spanned by sets of trees $\hat\sector_\ft$, we write $\ext \sector_\ft$ for the sector in $\ext\CT$ constructed as above, and we write
\[
\sector:=  \bigoplus_{\ft \in \FL_+} \sector_\ft
\qquad\text{ and }\qquad
\ext \sector:=  \bigoplus_{\ft \in \FL_+} \ext \sector_\ft.
\]
Moreover, assume we are given additionally exponents $\gamma=(\gamma_\ft)_{\ft \in \FL_+}$ and $\eta=(\eta_\ft)_{\ft \in \FL_+}$. For any model $Z \in \CM_0$ we then define the space $\CD^{\gamma,\eta}_\sector(Z)$ and a system of semi norms $\fancynorm{\cdot}_{\gamma,\eta,K}$ for any compact $K \ssq \domain$ by
\[
\CD^{\gamma,\eta}_\sector(Z) := \bigoplus_{\ft \in \FL_+}\CD^{\gamma_\ft,\eta_\ft}_{\sector_\ft} (Z)
\qquad \text{ and } \qquad
 \fancynorm{f}_{\gamma,\eta,K} := \sum_{\ft \in \FL_+}  \fancynorm{f_\ft}_{\gamma_\ft,\eta_\ft,K}.
\]

We fix from now on families of sectors $\CV_\ft$ and $\bar \CV_\ft$ in $\CT$ for $\ft \in \FL_+$, both spanned by sets of trees, and families of exponents $\gamma=(\gamma_\ft)_{\ft \in \FL_+}$, $\eta=(\eta_\ft)_{\ft \in \FL_+}$, $\bar \gamma=(\bar \gamma_\ft)_{\ft \in \FL_+}$ and $\bar \eta=(\bar \eta_\ft)_{\ft \in \FL_+}$ with $\gamma_\ft>0$ and $\bar\gamma_\ft \ge \gamma_\ft - |\ft|_\fs$ and $\eta_\ft, \bar \eta_\ft \in \R$. 
We then recall the following terminology from \cite[Sec.~7.3]{Hairer2014}
\begin{definition}
Given a model $Z \in \CM_0$, we call a map $F: \CD^{\gamma, \eta}_\sector(Z) \to \CD^{\bar \gamma,\ \bar \eta}_{\bar \sector}(Z)$ \emph{locally Lipschitz} if for any compact set $K\ssq \domain$ and any $R>0$ one has the bound
\begin{align}\label{eq:D:gamma:locally:Lipschitz}
\fancynorm{F(f) - F(g)}_{\bar\gamma, \bar\eta, K} \lesssim \fancynorm{f - g}_{\gamma,\eta,K}
\end{align}
uniformly over all $f,g \in \CD^{\gamma, \eta}_\sector(Z) $ with $\fancynorm{f}_{\gamma,\eta,K} \lor \fancynorm{g}_{\gamma,\eta,K} \le R$.
Given a locally Lipschitz map $F^Z: \CD^{\gamma,\eta}_\sector(Z) \to \CD^{\bar \gamma,\bar \eta}_{\bar \sector}(Z)$ for any model $Z \in \CM_0$ we call $F$ \emph{strongly locally Lipschitz} if for any $Z \in \CM_0$ there exists a neighbourhood $\CU$ of $Z$ in $\CM_0$ such that for any compact set $K\ssq \domain$, any $R>0$  one has the bound
\begin{align}\label{eq:D:gamma:strongly:locally:Lipschitz}
\fancynorm{F^Z(f) ;F^{\bar Z}(g)}_{\bar\gamma, \bar\eta, K} \lesssim \fancynorm{f; g}_{\gamma,\eta,K} + \fancynorm{Z, \bar Z}_{\gamma, \bar K}
\end{align}
uniformly over all models $\bar Z \in \CU$ and $f \in \CD^{\gamma, \eta}_\sector(Z) $ and $g \in \CD^{\gamma, \eta}_\sector(\bar Z)$ such that $\fancynorm{f}_{\gamma,\eta,K} \lor \fancynorm{g}_{\gamma,\eta,K} \le R$. Here, we set $\bar K:= \{x \in \domain: \mathrm{dist}(x,K)\le 1\}$.
\end{definition}

Following our usual convention, we will drop the dependence on the model $Z$ from the notation whenever there is no room for confusion. We say that $F$ is a strongly locally Lipschitz family for $(\sector,\bar \sector)$ if we want to emphasise the underling sectors. We want to consider a class of strongly locally Lipschitz families that admit lifts to the extended regularity structure as described in the next definition. 

\begin{definition}
Let $F$ be a strongly locally Lipschitz family for $(\sector, \bar \sector)$. Then we call a family $\ext F^{\ext Z,(r)}$ for $\ext Z \in \ext \CM_0$ and $r \in \R$ a \emph{lift} of $F$ if for any fixed $r \in \R$ the family $\ext F^{(r)}=(\ext F^{\ext Z,(r)})_{\ext Z \in \ext \CM_0}$ is a strongly locally Lipschitz family for $(\ext \sector, \ext {\bar\sector})$,
one has that $\ext F^{\ext Z,(r)}(f)$ is jointly Lipschitz continuous in $(f,\ext Z,r)$, i.e. one can strengthen (\ref{eq:D:gamma:strongly:locally:Lipschitz}) to
\begin{align}\label{eq:lift:lipschitz}
\fancynorm{F^{\ext Z,(r)}(f) ;F^{\ext {\bar Z},(s)}(g) }_{\bar\gamma, \bar\eta, K} \lesssim \fancynorm{f; g}_{\gamma,\eta,K} + \fancynorm{\ext Z, \ext {\bar Z}}_{\gamma, \bar K} + |r-s|
\end{align}
uniformly additionally for $|r|\lor |s|\le R$, and one has the identity
\begin{align}\label{eq:lift}
\SS_r F = \ext F^{(r)}\SS_r
\end{align}
on $\sector$ for any $r \in \R$ and $\ext Z \in \ext \CM_0$. Here, on the left hand side of (\ref{eq:lift}) we apply $F$ for the model $\SS_r^* \ext Z$ and on the right hand side we apply $\ext F^{(r)}$ for $\ext Z$. We call $\ext F$ a \emph{$C^1$-lift} if additionally one has that for any fixed model $\ext Z \in \ext\CM_0$ the map $(r,f) \mapsto F^{(r)}(f)$ is a Fr\'echet differentiable map from $\R \times \CD^{\gamma,\eta}_{\ext \sector}(\ext Z)$ into $\CD^{\bar \gamma,\bar \eta}_{\ext {\bar \sector}} (\ext {Z})$ with strongly locally Lipschitz continuous derivatives.
In the case that such a ($C^1$-) lift exists, we say that $F$ \emph{admits} a ($C^1$-) lift.
\end{definition}

\subsection{Shift of Abstract Fixed Point Problems}
\label{sec:abstract:fixed:point:shift}

In this section, if not explicitly otherwise stated, we make the notational convention that given $\ext Z \in \ext{\CM}_0$ and $r\in \R$ we write
\begin{align}\label{eq:notation}
\ext Z = (\ext \Pi, \ext \Gamma) 
\qquad\text{ and }\qquad
\SS_r^* \ext Z = Z = (\Pi,\Gamma).
\end{align}
We will show how to use lifts of strongly locally Lipschitz continuous non-linearities to relate abstract fixed point problems for the model $Z \in \CM_0$ to abstract fixed point problems for $\ext Z \in \ext \CM_0$ and consequently how to "shift" these fixed point problems in directions of Cameron-Martin functions. We start with the following Lemma. 

\begin{lemma}\label{lem:shift:D:gamma}
Fix $r\in \R$ and let $\sector  \ssq \CT$ be a sector of regularity $\alpha$. For any $\gamma>0$, $\eta \in \R$ and any $\ext{Z} \in \ext{\CM}_0$ the map $\SS_r$ is a Lipschitz continuous map from $\CD^{\gamma,\eta}_\sector  (\SS^*_r \ext{Z})$ into $\CD^{\gamma,\eta}_ {\ext \sector }(\ext Z)$, and provided that $\eta \le \gamma$ and $\alpha \land \eta >-\fs_0$ one has the identity
\begin{align}\label{eq:shift:D:gamma}
\CR^{\SS^*_r \ext{Z}} f = \CR ^ {\ext Z} \SS_r f
\end{align}
for any $f \in \CD^{\gamma,\eta}_\sector  (\SS^*_r \ext{Z})$.

Finally, $\SS_r$ maps strongly locally Lipschitz families $(F^{Z})_{Z \in \CM_0}$ onto strongly locally Lipschitz families $(\SS_r F^{\SS_r^* \ext Z})_{\ext Z \in \ext \CM_0}$.
\end{lemma}
\begin{proof}
In order to see that $\SS_r$ maps $\CD^{\gamma,\eta}_\sector (\SS^*_r \ext{Z})$ into $\CD^{\gamma,\eta} _{\ext \sector }(\ext Z)$, it suffices to use the identity $\ext \Gamma_{xy} \SS_r = \SS_r \Gamma_{xy}$ given by Lemma \ref{lem:properties:shift:operator}
and to note that $\SS_r$ does not change homogeneities of trees.
The identity (\ref{eq:shift:D:gamma}) is a direct consequence from the properties of the reconstruction operator, in particular \cite[(3.3)]{Hairer2014}, and the first identity in (\ref{eq:shift:operator:adjoint}).
\end{proof}

Let now $F$ be a strongly locally Lipschitz family for $(\sector,\bar \sector)$ and let $\ext F^{(r)}$ be a lift of $F$. We assume from now on that the pairs of sectors $(\sector_\ft, \bar\sector_\ft)$ are chosen such that $\CP_\ft:\bar\sector_\ft \to \sector_\ft$ for any $\ft\in\FL_+$. 
We also fix a strongly locally Lipschitz family $W^Z \in \CD^{\gamma,\eta}_\sector(Z)$ for $Z \in \CM_0$, and we define for any $\ext Z \in \ext \CM_0$ and $r \in \R$ the function
\begin{align}\label{eq:W:ext}
\ext W^{(r), \ext Z} := \SS_r W^{\SS_r^* \ext Z}
\end{align}
so that $\ext W^{(r),\ext Z} \in \CD^{\gamma,\eta}_{\ext \sector}(\ext Z)$ by Lemma \ref{lem:shift:D:gamma}.

We consider the fixed point problems for $U$ and $\ext U_\ft^{(r)}$ given by
\begin{align}\label{eq:FPP}
U_\ft = \CQ_{< \gamma_\ft}\CP_\ft[ \I_+ F_\ft(U)] + W_\ft
\quad\text{ and }\quad
\ext U_\ft^{(r)} = \CQ_{< \gamma_\ft} 
	\CP_\ft[ \I_+ {\ext {F}^{(r)} _\ft}(\ext U^{(r)})] + \ext W^{(r)}_\ft.
\end{align}
in $\CD^{\gamma,\eta}_\sector( Z )$ and $\CD^{\gamma,\eta}_{\ext \sector}(\ext Z)$ respectively for any model $Z \in \CM_0$ and $\ext Z \in \ext\CM_0$ and any $r\in \R$.
Since $\bar\gamma_\ft \ge \gamma_\ft - |\ft|_\fs$ and the right hand sides are locally Lipschitz continuous by definition, it follows from \cite[Thm.~7.8]{Hairer2014} (see also \cite[Thm.~6.21]{BrunedChandraChevyrecHairer2017})
 that there exist unique maximal solutions $U$ and $\ext U^{(r)}$ to these equations. We denote the maximal time of existence by $T(Z)$ and $T(\ext Z,r)$, respectively. We also define the stopping time $\tau(\omega) := T(\hat Z^\xi_\BPHZ(\omega))$.
In this setting, we have the following result.

\begin{proposition}\label{prop:FPP:shift}
Fix $\ext Z \in \ext \CM_0$ and $r \in \R$, and let $U$ and $\ext U^{(r)}$ be the unique solutions to the fixed point equations (\ref{eq:FPP}) for the models $\SS^*_r \ext Z$ and $\ext Z$,
respectively. Then one has 
\begin{align}\label{eq:FPP:shift}
\ext U^{(r)} = \SS_r U.
\end{align}
\end{proposition}
\begin{proof}
Since the solution to these fixed point equations are unique,  we only need to show that $\SS_r U$ satisfies the second equation in (\ref{eq:FPP}).  For this note that 
\[
{\ext {F}^{(r)}_\ft}(\SS_r U) = \SS_r F_\ft(U),
\]
and since it follows directly from the definition that one has $\CP_\ft \I_+\SS_r = \SS_r \CP_\ft \I_+$ the result follows.
\end{proof}

We stress at this point that the function $U$ on the right hand side of (\ref{eq:FPP:shift}) depends on~$r$ through the model $\SS_r^* \ext Z$.

\begin{remark}\label{rem:lower:semi:continuity}
In \cite[Thm.~7.3]{Hairer2014} it was shown that if $W^Z$ is locally Lipschitz continuous in the model, then the time of existence $T = T(Z)$ of the solution $U$ is a lower semi-continuous map in the model $Z$. Proposition \ref{prop:FPP:shift} shows that the time of existence of $\ext U^{(r) , \hat Z^{\xi,h}_\BPHZ(\omega)}$ is additionally locally Lipschitz continuous in $r \in \R$.
\end{remark}

\begin{remark}
The modelled distribution $\ext W^{(r), \ext Z}$ will be used in order to deal with the initial condition.
The assumption that $W^Z$ is locally Lipschitz continuous in the model $Z \in \CM_0$ is unreasonably strong, as it ends up imposing that the initial condition is a locally Lipschitz continuous map of the model. 
For the existence of the local $H^\xi$-Fr\'echet derivative it is sufficient to assume that for any $\omega \in N^c$ the map $H^\xi \times \R \to \CD^{\gamma,\eta}_{\ext \sector}(\ext Z)$ given by $$(h,r) \mapsto \ext W^{(r),\hat Z^{\xi,h}_\BPHZ(\omega)}$$ is locally Lipschitz continuous, which is a trivial consequence of the assumption that the initial condition is locally $H^\xi$-Fr\'echet differentiable.
Under this less restrictive assumption, the statement of Remark \ref{rem:lower:semi:continuity} is no longer true. However, this assumption is still sufficient to ensure in the same way as above that the time of existence is lower semi-continuous with respect to $r$, which ultimately ends up ensuring that the time of existence of the solution $u$ is lower semi-continuous with respect to Cameron-Martin shifts. A similar Remark applies to the functions $\ext V^{\ext Z,(r)}$ and $r_0( \ext Z, T)$ introduced in Proposition \ref{prop:diff:abstract} below.
\end{remark}

\section{The Malliavin Derivative}

We show in this section that the reconstructed solutions to the abstract fixed point problems considered in Sections \ref{sec:abstract:fixed:point} and \ref{sec:abstract:fixed:point:shift} admit a local $H^\xi$-Fr\'echet derivative. In Section \ref{sec:H:diff:SPDE} we apply this abstract result to singular SPDEs of the form (\ref{eq:singular:SPDE}).

\subsection{Differentiability of the Solution to the Abstract Fixed Point Problem}
\label{sec:Gateaux}

We show that $\ext U^{(r)}$ is differentiable in $r \in \R$ with values in $\CD^{\gamma,\eta}_\sector$. For this let $F$ be a $C^1$-liftable, strongly locally Lipschitz family for a pair of sectors $(\sector,\bar \sector)=(\CV_\ft,\bar \CV_\ft)_{\ft \in \FL_+}$ with the property that $\CP_\ft : \bar\CV_\ft \to \CV_\ft$ for any $\ft \in \FL_+$, and let $\ext F^{(r)}$ be a $C^1$-lift of $F$. 
Let moreover $W^Z \in \CD^{\gamma,\eta}_\sector(Z)$ for $Z \in \CM_0$ be a family as in the previous section, and assume that additionally the map $r \mapsto \ext W^{(r),\ext Z}$ defined as in (\ref{eq:W:ext}) is Fr\'echet differentiable as a map from $\R$ into $\CD^{\gamma,\eta}_{\ext \sector}(\ext Z)$ for any $\ext Z \in \ext \CM_0$. Finally, for any $\ext Z \in \ext \CM_0$ let $U$ and $\ext U^{(r)}$ be the solutions to (\ref{eq:FPP}) for $\SS_r^*\ext Z$ and $\ext Z$, respectively. 
We then have the following.

\begin{proposition}\label{prop:diff:abstract}
Under the assumption at the beginning of this section, for any $\ext Z \in \ext \CM_0$ and any $T< T(\ext Z, 0)$ there exists $r_0 = r_0(\ext Z,T)> 0 $ such that the map $r \mapsto \ext U^{(r)}$ is $C^1$ as a map from $(-r_0,r_0)$ into $\CD^{\gamma,\eta,T}_{\ext \sector}(\ext Z)$ and its derivative $\ext V^{(r)}$ satisfies the fixed point equation
\begin{align}\label{eq:FPP:derivative}
\ext V^{(r)} _\ft= D \ext F^{(r)}_\ft(\ext U^{(r)}) \ext V^{(r)} + (\partial_r \ext F ^{(r)}_\ft) (\ext U^{(r)})
+
\partial_r \ext W^{(r)}_\ft.
\end{align}
for any $\ft \in \FL_+$.
Moreover, the function $\ext V^{(r)}$ is strongly locally Lipschitz continuous in the sense that for any $\ext Z \in \ext \CM_0$ and any $T< T(\ext Z, 0)$ there exists a neighbourhood of $\CU$ of $\ext Z$ such that one has
\begin{align}\label{eq:FPP:derivative:lipschitz}
\fancynorm{ \ext V^{\ext Z, (r)}; \ext { V}^{\ext{ \bar Z},(s)}}_{\gamma,\eta,T}
\lesssim
\fancynorm{ \ext Z ; \ext {\bar Z} }_{\gamma,[-1,T+1]\times \T^d} + |r-s|
\end{align}
uniformly over $\ext {\bar Z} \in \CU$ and $r,s \in (-r_0+\theta, r_0-\theta)$ for any $\theta>0$.

Finally, $r_0(\ext Z,T) \in (0,\infty]$ can be chosen as the supremum over all $r_0>0$ such that $T<T(\ext Z, s)$ for any $s \in (-r_0, r_0)$, and with this choice for any fixed $T$ the function $\ext Z \mapsto r_0(\ext Z, T)$ is lower semicontinuous in $\ext Z \in \ext \CM_0$. 
\end{proposition}

\begin{proof}
We fix for the first part of the proof a model $\ext Z \in \ext \CM_0$. By definition one has that $(r, \ext U) \mapsto \ext{F}_\ft^{(r)}(\ext U)$ is a Fr\'echet differentiable map from $\R \times {\CD} _{\ext \sector}^{\gamma,\eta, T} ( \ext Z )$  into ${\CD}_{\ext {\bar \sector}}^{\bar \gamma , \bar \eta, T} ( \ext Z )$ for any model $\ext Z \in \ext \CM_0$.
In order to see that $r \mapsto \ext U^{(r)}$ is Fr\'echet differentiable, we make use of the implicit function theorem
\footnote{Compare \cite[Prop.~4.7]{CannizzaroFrizGassiat2015} and references therein where the idea to use the implicit function theorem to show differentiability of the solution map was also present.}
on the map $\Phi:\R \times {\CD}_{\ext \sector}^{\gamma,\eta, T}(\ext Z) \to {\CD}_{\ext \sector}^{\gamma,\eta, T} (\ext Z)$ given by 
$$
\Phi:(r,\ext U) \mapsto \ext U - (\CQ_{<\gamma_\ft}\CP_\ft[\I_+ \ext{F}^{(r)}_\ft(\ext U)] + \ext W^{(r)})_{\ft \in \FL_+}
$$
in a neighbourhood of $(0,\ext {U}^{(0)})$. 
Note that  by (\ref{eq:FPP}) one has $\Phi(0,\ext{U}^{(0)})=0$. Since $\CQ_{<\gamma_\ft}\CP_\ft\I_+$ is a bounded linear operator from $\CD^{\bar \gamma_\ft, \bar \eta_\ft, T}_{\ext{\bar \sector}_\ft}$ into $\CD^{\gamma_\ft, \eta_\ft, T}_{\ext \sector_\ft}$ by Proposition \ref{prop:Dgamma:loc:Lipschitz} it also Fr\'echet differentiable, and since by assumption one has that $r \mapsto \ext W^{(r)}$ is Fr\'echet differentiable as well, it follows that $\Phi$ is a $C^1$ map from $\R \times {\CD}_{\ext \sector}^{\gamma,\eta, T}$ into $\CD^{\gamma, \eta, T}_{\ext \sector}$.
We now show that the derivative of $\Phi$ with respect to $\ext U$ is an isomorphism $D_2 \Phi(r,\ext U) : {\CD}_{\ext \sector}^{\gamma,\eta, T}(\ext Z) \to {\CD}_{\ext \sector}^{\gamma,\eta, T}(\ext Z)$ for any $(r,\ext U) \in \R \times {\CD}_{\ext \sector}^{\gamma,\eta, T}$. 
By definition, the derivative $D_2 \Phi(r,\ext U)$ is a bounded linear operator between these spaces, which is given by
\begin{align}\label{eq:diff:abstract:2}
D_2 \Phi(r,\ext U)(\ext V)
=
\ext V
-
(\CQ_{<\gamma_\ft}\CP_\ft[\I_+ D\ext{F}^{(r)}_\ft(\ext U)\ext V_\ft])_{\ft \in \FL_+},
\end{align}
so that we are left to show that this expression is invertible. This is equivalent to solving, for any fixed $X \in {\CD}_{\ext \sector}^{\gamma,\eta,T}(\ext Z)$ the equation $D_2 \Phi(r,\ext U)(\ext V) = X$ for $\ext V \in {\CD}_{\ext \sector}^{\gamma,\eta, T}(\ext Z)$, which is in the form of a fixed point problem and admits a unique solution by \cite[Thm.~7.8]{Hairer2014}. This follows from the fact that the map $\ext V_\ft \mapsto D\ext{F}^{(r)}_\ft(\ext U)\ext V_\ft$ from ${\CD}_{\ext \sector}^{\gamma,\eta,T}(\ext Z)$ into ${\CD}_{\ext {\bar \sector}}^{\bar \gamma,\bar \eta,T}(\ext Z)$ is linear and continuous, and thus it is also Lipschitz continuous.

It now follows that there exists a neighbourhood $\CU$ of $\ext{U}^{(0)}$ and $r_0>0$ and a (unique) $C^1$-function $r\mapsto \ext u(r) \in \CU$ with $\ext u(0) = \ext U^{(0)}$ such that
\[
\Phi(r, \ext u(r))=0
\qquad \text{ for } r \in (-r_0,r_0).
\]
By uniqueness of solutions to the fixed point problem (\ref{eq:FPP}) we infer that one has necessarily $\ext u(r)= \ext U^{(r)}$, so that it follows in particular that $\ext U^{(r)}$ is $C^1$ in $(-r_0,r_0)$.
In order to see the identity (\ref{eq:FPP:derivative}) for the derivative, note that at this point all functions appearing in (\ref{eq:FPP}) are Fr\'echet differentiable in $r\in (-r_0, r_0)$, so that (\ref{eq:FPP:derivative}) follows by differentiating the right hand side of this identity.

In order to see (\ref{eq:FPP:derivative:lipschitz}) it suffices to show local Lipschitz continuity in $r$ and $\ext Z$ separately. The former follows from arguments identically to above, which shows the stronger statement of Fr\'echet differentiability of $V^{(r)}$ in $r$. For the latter, by \cite[Thm.~7.8]{Hairer2014} it suffices to show that $\ext V_\ft \mapsto  D\ext{F}^{(r)}_\ft(\ext U^{(r)})\ext V_\ft + (\partial_r \ext F ^{(r)}_\ft) (\ext U^{(r)}) + \partial_r \ext W ^{ (r) }_\ft$ is strongly locally Lipschitz continuous between the spaces $ \CD^{\gamma,\eta,T}_{\ext \sector}$ and $\CD^{\bar\gamma, \bar \eta, T}_{\ext{ \bar \sector}}$. This however follows by combining the fact that by definition the Frecht derivatives $D{\ext F}_\ft^{(r)}$ and $\partial_r \ext F_\ft^{(r)}$ are strongly locally Lipschitz continuous and the solution $\ext U^{(r)}$ is strongly locally Lipschitz continuous as a consequence of \cite[Thm.~7.8]{Hairer2014}.

For the last part of the theorem assume that $r_0$ has been chosen maximally, and assume that for some $r_1>r_0$ one has that $T(\ext Z, s)<T$ for all $r\in (-r_1, r_1)$. We can then redo the arguments in the first part of the proof with $(0, \ext U^{(0)} )$ replaced by $(r_0, \ext U^{(r_0)})$ to obtain $s_0>0$ such that $s \mapsto \ext U^{(r_0+s)}$ is $C^1$ as a map from $(-s_0, s_0)$ into $\CD^{\gamma,\eta,T}_{\ext V}(\ext Z)$, which shows that $U^{(r)}$ is a $C^1$ map on $(-r_0, r_0 + s_0)$. A similar argument shows that the lower bound can be improved and yields a contradiction. The lower semicontinuity of $r_0$ is now a consequence of the lower semicontinuity of $T(\ext Z,s)$ in $(\ext Z,s)$.
\end{proof}

\subsection{Local $H$-differentiability of the Solution}
\label{sec:H:diff}
As a consequence of Proposition \ref{prop:diff:abstract} we can show that the reconstructed solution map $u=\CR U$ is Gateaux differentiable in $H^\xi$ directions. 

\begin{lemma}
Let $\omega \in N^c$ and let $h \in H^\xi$. Then for any $T< \tau(\omega)$ there exists $r_0>0$ such that the map
\[
r \mapsto u(\omega + r h)
\]
with values in $\bigoplus_{\ft \in \FL_+} \CC^{\alpha_\ft}((0,T) \times \T^d)$, where $\alpha_\ft$ denotes the regularity of the sector $\sector_\ft$, is Fr\'echet differentiable differentiable on $(-r_0,r_0)$.
\end{lemma}
\begin{proof}
Let $r_0 = r_0 ({ \ext {\hat Z}^{\xi, h}_\BPHZ (\omega) }, T )$ be as in Proposition \ref{prop:diff:abstract}.
By Theorem \ref{thm:shift:operator} and Proposition \ref{prop:FPP:shift} one has that 
\[
u(\omega + rh) 
= \CR^{\SS_r^* \ext{\hat Z}_\BPHZ^{\xi,h}(\omega) } U
= \CR^{ \ext{\hat Z}_\BPHZ^{\xi,h}(\omega) } \ext{U}^{(r)}.
\]
Since $r \mapsto \ext U^{(r)}$ is $C^1$ with values in $\CD^{\gamma,\eta, T}_{\ext \sector}$ by Proposition \ref{prop:diff:abstract}, the result follows from the fact that for any fixed model $\ext Z$ the reconstruction operator $\CR^{\ext Z}$ is a bounded linear map on $\CD^{\gamma,\eta,T}_{\ext \sector}$ and thus Fr\'echet differentiable.
\end{proof}

Finally, we can show the following.

\begin{theorem}\label{thm:local:H:diff}
For any $T>0$ the solution $u$ restricted on $(0,T)\times\T^d$ and conditioned on the event $\{ \tau > T \}$ with values in $\bigoplus_{\ft \in \FL_+} \CC^{\alpha_\ft}((0,T) \times \T^d)$ is locally $H^\xi$-Fr\'echet differentiable in the sense of Definition \ref{def:local:H1:diff}.
\end{theorem}
\begin{proof}
For fixed $\omega \in N^c$ and $h \in H^\xi$ let $r_0(\omega,h) := r_0 ({ \ext {\hat Z}^{\xi, h}_\BPHZ (\omega) } )$ be as in Proposition~\ref{prop:diff:abstract}.
Then one has the identity
\[
\partial_r u(\omega + rh)
=
\CR^{ \ext{\hat Z}_\BPHZ^{\xi,h}(\omega) } \ext V^{(r)}_h,
\]
for $r \in (-r_0(\omega,h), r_0(\omega,h) )$, where $\ext V^{(r)}_h$ denotes the derivative of $U^{(r)}$ in the direction of $r$ as in Proposition \ref{prop:diff:abstract} for the model $ \ext{\hat Z}_\BPHZ^{\xi,h}(\omega)$. 
Since $h \mapsto \ext{\hat Z}_\BPHZ^{\xi,h}(\omega)$ is locally Lipschitz continuous by Corollary \ref{cor:extension:operator}, it follows that for any fixed $\omega \in N^c$ the map $h \mapsto r_0(\omega, h)$ is lower semi-continuous. Since furthermore one has $r_0(\omega, 0)=+\infty$, there exists $\mu >0$ and a ball  $B_\mu(0) \ssq H^\xi$ around the origin such hat one has $r_0(\omega,h)>1$ for any $h \in B_\mu(0)$. Now it follows from (\ref{eq:shift:operator}) that for any $h,k \in B_\mu(0)$ one has 
\[
u(\omega+h + rk) 
=
\CR^{ {\hat Z}^{\xi}_\BPHZ (\omega + h + rk) } U
=
\CR^{ {\ext {\hat Z}}^{\xi,k}_\BPHZ (\omega + h) } \ext U^{(r)}
\]
so that it follows in particular from Proposition \ref{prop:diff:abstract} that $h \mapsto u(\omega+h)$ is Gateaux  differentiable in $B_\mu(0)$ with G\^ateaux derivative given by 
\[
\partial_r |_{r=0} u(\omega + h + rk) = 
\CR^{ {\ext {\hat Z}}^{\xi,k}_\BPHZ (\omega + h) } \ext V^{(0)},
\]
so that it remains to show that this expression is continuous in $(h,k) \in H^\xi \times H^\xi$.
This follows from the fact that $\CR$ is strongly Lipschitz continuous, the map $(h,k) \mapsto {\ext {\hat Z}}^{\xi,k}_\BPHZ (\omega + h)$ is locally Lipschitz continuous and by (\ref{eq:FPP:derivative:lipschitz}) one has that $\ext V^{(0)}$ is strongly locally Lipschitz continuous.
\end{proof}

\subsection{Application to subcritical SPDEs}
\label{sec:H:diff:SPDE}

We now apply the result of the previous section to abstract fixed point problems arising from singular SPDEs and conclude the proof of Theorem \ref{thm:main:malliavin}. That is, we show that under the assumption introduced in Theorem \ref{thm:main:malliavin}, the solution $u$ to the singular SPDE (\ref{eq:singular:SPDE}) admits a local $H$-Fr\'echet derivative, and in this situation we can furthermore derive a "tangent equation" (\ref{eq:tangent:equation:main}) for this Fr\'echt derivative, which is informally given by differentiating the original equation with respect to the noise. The precise meaning of (\ref{eq:tangent:equation:main}) is that $v_h$ can be written as a limit $v_h = \lim_{\eps \to 0} v_h^\eps$, where the random smooth function $v_h^\eps = (v_{\ft,h}^\eps)_{\ft \in \FL_+}$ is the unique classical solutions to the system of equations
\begin{multline}\label{eq:tangent:equation:reg}
\partial_t v_{\ft,h}^\eps = \CL_\ft v_{\ft,h}^\eps
+
DF_\ft(u^\eps, \nabla u^\eps, \ldots)(v_{h}^\eps, \nabla v_{h}^\eps, \ldots) 
\\
+ 
\sum_{\Xi \in \FL_-} DF_\ft^\Xi(u^\eps, \nabla u^\eps, \ldots)(v_{h}^\eps, \nabla v_{h}^\eps, \ldots) \xi^\eps_\Xi
\\
+
\sum_{\Xi \in \FL_-} F_\ft^\Xi(u^\eps, \nabla u^\eps, \ldots) h_\Xi^\eps 
\\
+
\sum_{\tau \in \TT_{\ft,-}(F)} c^\eps_\tau D\Upsilon_\ft^\tau(u^\eps, \nabla u^\eps, \ldots)(v_{h}^\eps, \nabla v_{h}^\eps, \ldots)
\end{multline}
with initial condition $v_{\ft,h}^\eps (0) = D_{h} u_{\ft,0}^\eps$. It is not hard to see that $v_h^\eps$ is the $H^\xi$-Fr\'echet derivative $D_h u^\eps$ of the solution $u^\eps$ to the regularized and renormalized equation (\ref{eq:singular:SPDE:reg}) in the direction of $h \in H^\xi$. Note that both $\psi^\eps$ and $\CS_{\rho,\eps}^-(\xi)(0,\cdot)$ are locally $H^\xi$-Fr\'echet differentiable (for the former this follows by assumption, while for the latter this follows from the explicit definition of $\SS_{\rho,\eps}^-(\xi)$ in \cite[(A.11),(6.10)]{BrunedChandraChevyrecHairer2017}, which imply in particular that $\CS_{\rho,\eps}^-(\xi)(0,\cdot)$ takes values in some inhomogeneous Wiener chaos), so that the same is true for $u_0^\eps = \psi^\eps + \CS_{\rho,\eps}^-(\xi)(0,\cdot)$.

\begin{remark}
The tangent equation (\ref{eq:tangent:equation:main}) is in the form of a singular SPDE, however it does not fall under the setting of \cite{BrunedChandraChevyrecHairer2017} since it involves a source $h_\Xi$ which is deterministic and in general not smooth. The fact that $h_\Xi$ is not necessarily smooth was the main reason that the analysis of Section \ref{sec:ext:models} was necessary in the first place. However, if $h_\Xi$ happens to be smooth for any $\Xi \in \FL_-$, then one can treat the tangent equation directly in the framework of \cite{BrunedChandraChevyrecHairer2017}, and in this case the regularized and renormalized equation derived in \cite{BrunedChandraChevyrecHairer2017} coincides with (\ref{eq:tangent:equation:reg}). 
\end{remark}

The solution $u$ constructed in \cite{BrunedChandraChevyrecHairer2017} is given as the reconstruction of $\bar U_\ft + \CP_\ft \tilde U_\ft$ (c.f. \cite[Prop.~6.22]{BrunedChandraChevyrecHairer2017}), where $\tilde U_\ft$ is the constant modelled distributions in $\CD^{\infty,\infty}$ explicitly given in \cite[(6.10)]{BrunedChandraChevyrecHairer2017}. 
We now introduce an abstract differentiation operator $\CD: \CT \to \ext{\CT}$ as the derivative of $\SS_r$ at $r=0$, so that one has
\[
\CD \tau := \partial_r|_{r=0}  \SS_r \tau = \sum_{\sigma \in \SS[\tau], m(\sigma)=1 } \sigma.
\]
With this notation we can see that the $H^\xi$-Fr\'echet derivative of the function $\CR \CP_\ft \tilde U_\ft$ in the direction of $h \in H^\xi$ is  given by $D_h \CR^{\hat Z^\xi_\BPHZ} \CP_\ft \tilde U_\ft = \CR^{\ext{\hat Z}^{\xi,h}_\BPHZ} \CD \CP_\ft\tilde U_\ft$.
The fact that $u_\ft$ is locally $H^\xi$-Fr\'echet differentiable is now equivalent to showing that $\CR\bar U_\ft$ is locally $H^\xi$-Fr\'echet differentiable, which at this stage is an application of Theorem \ref{thm:local:H:diff} to the abstract fixed point problem \cite[(6.16)]{BrunedChandraChevyrecHairer2017} for $\bar U$. The main step here is to show that the right hand side of this fixed point problem admits $C^1$ lifts, and since this is largely a technical issue, we postpone the proof to Appendix \ref{sec:application:to:SPDEs} below. 

It remains to derive the tangent equation. For this consider a regularization $\xi^\eps \in \SM_\infty(\FL_-)$ of $\xi$ given by $\xi^\eps_\Xi := \xi_\Xi * \rho^{(\eps)}$ and let $h_\Xi^\eps := h_\Xi*\rho^{(\eps)}$, let $v_{h,\ft}^\eps := D_h u_{\ft}^\eps$ be the local $H^\xi$-Fr\'echet derivative of the solution $u_\ft^\eps$ to (\ref{eq:singular:SPDE:reg}) in the direction of $h \in H^\xi$.
The fact that one has $v_{h,\ft} = \lim_{\eps \to 0 }v_{h,\ft}^\eps$ follows simply from the fact that both sides of this equation are given as the reconstruction of the modelled distribution $\ext V^{(0)}+\CD\tilde U$, where $\ext V^{(0)}$ denotes the solution to the abstract tangent equation (\ref{eq:FPP:derivative}), for the model $\ext {\hat Z}_\BPHZ^{\xi,h}$ and $\ext {\hat Z}_\BPHZ^{\xi^\eps,h^\eps}$, respectively, and the fact that the latter converges to the former as $\eps \to 0$ by Theorem \ref{thm:extension:operator}.
It is now sufficient to show that $v_{h,\ft}^\eps$ solves (\ref{eq:tangent:equation:reg}). But in the regularized case the map $r \mapsto \xi^\eps(\omega + rh)$ is smooth with derivative given by $\partial_r |_{r=0} \xi^\eps(\omega + rh) = h^\eps$. Furthermore, $(r,x) \mapsto \xi^\eps(\omega+rh)(x)$ is a smooth function $\R\times\domain \to \R^{\FL_-}$ and since both $F$ and $\Upsilon$ are smooth, the former by assumption and the latter by construction, c.f. (\ref{eq:counter-term}), it follows readily from standard Schauder estimates that the map $(r,x) \mapsto u^\eps(\omega +rh)(x)$ is smooth as well. This is sufficient to argue that we are allowed to commute the differentiation operators $\partial_t$ and $\CL_\ft$ with $\partial_r$ in (\ref{eq:singular:SPDE:reg}), and since per definitionem one has that $v_{h,\ft}^\eps = \partial_r|_{r=0} u_\ft^\eps(\omega + rh)$, we obtain~(\ref{eq:tangent:equation:reg}) by a direct computation.

\section{Density of Solutions to singular SPDEs}
\label{sec:densities}

In this section we always fix a time $T>0$. We want to derive conditions such that random variables of the form $X= \dual{u,\phi}:=\sum_{\ft \in \FL_+} \dual{u_\ft, \phi_\ft}$, for 
some tuple of test functions $\phi \in \bigoplus_{\ft \in \FL_+}\CC_c^\infty((0,T)\times \T^d)$, conditioned on the event $\{\tau>T \}$ admit a density with respect to Lebesgue measure. By Theorem \ref{thm:main:malliavin} the random variable $X$ is locally $H^\xi$-Frechet differentiable with derivative in direction of $h \in H^\xi$ given by $\dual{v_{h}, \phi}$, where $v_{h}$ solves (\ref{eq:tangent:equation:main}). By the Bouleau Hirsch criterion, Corollary \ref{cor:Bouleau:Hirsch}, we are lead to study non-degeneracy of this local $H^\xi$-Fr\'echet derivative.

For simplicity we make in this section the following additional assumption.
\begin{assumption} \label{ass:dual:equation:simplicitiy}
We assume that the following is satisfied.
\begin{itemize}
\item The renormalization constants are given by the BPHZ character $c_\tau^\eps = g_\BPHZ^\eps(\tau)$.
\item For any $\Xi \in \FL_-$ and any $\ft \in \FL_+$ one has that $F_\ft^\Xi=F_\ft^\Xi(u)$ and $F_\ft=F_\ft(u)$ depend only on the solution (and not on its derivatives).
\item For any noise type $\ft \in \FL_-$, any $\tau \in \TT_{\ft,-}^F$ and any $\eps>0$ with the property that $c_\tau^\eps \ne 0$ one has either $\Upsilon^\tau_\ft= \Upsilon^\tau_\ft(u)$ depends only on the solution (and not on its derivatives), or $\Upsilon^\tau_\ft(u, \nabla u, \ldots) = \partial_i u_\ft$ for some $i\in \{0,  \ldots , d \}$. We write $\tau \in \TT_{\ft,-}^{F,\nond}$ in the first case and $\tau \in \TT_{\ft,-}^{F,i}$ for $i\in \{0, \ldots d \}$ in the second case.
\end{itemize}
\end{assumption}

The first assumption is merely a convenience and could be easily dropped with a little more algebraic effort later on, compare Remark \ref{rem:c=BPHZ} below. The second and the third assumption greatly simplify the computation; we believe that the statements below are still true without these assumptions, but the proofs given in this paper do not seem to easily generalize to this case.
Under Assumption \ref{ass:dual:equation:simplicitiy} equation (\ref{eq:tangent:equation:reg}) simplifies to
\begin{multline}\label{eq:tangent:equation:reg:simplified}
\partial_t v_{\ft,h}^\eps = \CL_\ft v_{\ft,h}^\eps
+
DF_\ft(u^\eps)v_{h}^\eps 
+ 
\sum_{\Xi \in \FL_-} DF_\ft^\Xi(u^\eps)v_{h}^\eps \xi^\eps_\Xi
+
\sum_{\Xi \in \FL_-} F_\ft^\Xi(u^\eps) h_\Xi^\eps 
\\
+
\sum_{\tau \in \TT_{\ft,-}^{F,\nond}} c^\eps_\tau 
D\Upsilon_\ft^\tau(u^\eps) v_{h}^\eps 
+
\sum_{i=1}^d \sum_{\tau \in \TT_{\ft,-}^{F,i}} c^\eps_\tau \partial_i v_{\ft,h}^\eps.
\end{multline}
We denote by $\CL_\ft^*$ the dual operator to $\CL_\ft$, which is again a differential operator involving only spatial derivatives, and we
consider the equation dual to (\ref{eq:tangent:equation:reg:simplified}), which is a backward random PDE given by
\begin{multline}\label{eq:dual:equation}
-\partial_\ft w_{\ft,\phi}^\eps = \CL_\ft^* w_{ \ft , \phi }^\eps
+
DF_\ft( u^\eps ) w_{\phi}^\eps 
+ 
\sum_{\Xi \in \FL_-} DF_\ft^\Xi(u^\eps) w_{\phi}^\eps \xi^\eps_\Xi
+
\phi_\ft
\\
+
\sum_{\tau \in \TT_{\ft,-}^{F,\nond}} c^\eps_\tau 
D\Upsilon_\ft^\tau(u^\eps) w_{\phi}^\eps 
-
\sum_{i=1}^d \sum_{\tau \in \TT_{\ft,-}^{F,i}} c^\eps_\tau \partial_i w_{\ft,\phi}^\eps,
\end{multline}
on $(0,T)\times \T^d$ with finial condition $w_{\ft,\phi}^\eps(T,\cdot)= 0$. The following lemma is a straight-forward computation.
\begin{lemma}
Let $T>0$ and let $\phi_\ft \in \CC_c^\infty((0,T)\times \T^d)$ for any $\ft \in \FL_+$. Then for any $h \in L^2(\domain) ^{\FL_-}$ one has the identity
\begin{align} \label{eq:dual:identity:reg}
\dual{ v_h^\eps , \phi }_{L^2(\domain)} = \dual{ \sum_{\Xi \in \FL_-} F^\Xi(u^\eps) h_\Xi^\eps  , w_\phi^\eps }_{L^2(\domain)} + 
\dual{v^\eps_h(0,\cdot), w^\eps_\phi(0, \cdot )}_{L^2(\T^d)}
\end{align}
between random variables conditioned on the event $\{ T < \tau \}$.
\end{lemma}

\subsection{A Regularity Structure Adapted to the Dual Equation}
\label{sec:Kext:reg:structure}

Our goal is to understand the behaviour of $w^\eps_{\ft,\phi}$ in the limit $\eps \to 0$. To this end we want to interpret $w^\eps_{\ft,\phi}$ as the reconstructed solution to an abstract fixed point problem, which can be viewed as the dualization of the abstract tangent equation (\ref{eq:FPP:derivative}). The equation for $w$ can be written in its mild formulation, and it is not hard to see that the Greens function for $-\partial_t-\CL_\ft^*$ is given by $x\mapsto G_\ft(-x)$ for any $\ft \in \FL_+$, where $G_\ft$ is the Greens function for $\partial_t-\CL_\ft$. It follows that the kernel types present in our regularity structure $\CT^\ex$ are not rich enough to encode the dual equation, so that as a first step we are lead to build an extension $\Kext\CT^\ex$ of the regularity structure $\CT^\ex$ in which one can consider abstract fixed point problems associated to $G_\ft(-\cdot)$. 

\begin{remark}
Note that in Section \ref{sec:ext:reg:structure} we considered an extension of the noise-types $\FL_-$, whereas in this section we will consider an extension of the kernel-types $\FL_+$. The extension constructed in Section \ref{sec:ext:reg:structure} plays no role in this section, so that from now on we use the symbol $\Kext \CT^\ex$ for the regularity structure constructed below, and we refer to this structure as the \emph{extended regularity structure} form now on.
\end{remark}

To this end we extend the set of kernel types to a set $\Kext\FL_+:=\FL_+ \sqcup \Kdual {\FL_+}$, where $\Kdual {\FL_+} := \{\Kdual \ft : \ft \in \FL_+\}$ is a disjoint copy of $\FL_+$, and we let $|\Kdual \ft|_\fs:=|\ft|_\fs$ for any $\ft \in \FL_+$. One should think of $\Kdual\ft$ as representing the "dual" integral operator to $\ft$. In particular, it will represent the Greens function for a parabolic differential operator going backward in time. 

Given the extended set of types $\Kext \FL := \Kext \FL_+ \sqcup \FL_-$ we extend $\reg$ to a function $\reg : \Kext\FL \to \R$ by setting $\reg( \Kdual \ft ):= \theta$ for any $\ft \in \FL_+$ and some $\theta>0$ small enough, and we define an extension $\Kext R$ of the rule $R$ by allowing any kernel type $\ft \in \FL_+$ to be replaced by $\Kdual \ft$, and additionally allowing an arbitrary number of types $\Kdual{\FL_+}$. To be more precise, we define $\KdualCNzero := \N^{ \Kdual {\FL_+} \times \{0\}}$ and $\Kext\CN:=\N^{\Kext\Psi}$, with $\Kext\Psi:=\Kext{\FL} \times \N^{d+1}$, and we set
\begin{align}
\label{eq:Kext:rule}
\Kext R(\ft)
&:= \{ N \sqcup M \in \Kext \CN : 
\cq N \in R(\cq\ft), M \in \KdualCNzero\} 
\\
\label{eq:Kext:rule:dual}
\Kext R(\Kdual \ft)
&:=
\{ N \in \Kext R(\ft) : \exists (\fl,k) \in \FL_+ \times \N^{d+1} 
	: N \sqcup \{(\fl,k)\} \in \Kext R(\ft) \}
\end{align}
for any $\ft \in  \FL_+$, and as usual
$\Kext R(\ft) := \{\emptyset\}$ for $\ft \in \FL_-$. Here we define $\cq\ft:=\ft$ for any $\ft \in \FL$, $\cq\Kdual\ft:=\ft$ for any $\ft \in \FL_+$ and 
\[
(\cq N)_{(\ft,k)}:=\sum_{\cq\tilde\ft=\ft}N_{(\tilde\ft,k)}
\]
where the sum runs over all $\tilde\ft \in \Kext\FL$ with $\cq\tilde\ft=\ft$. 
\begin{remark}
The fact that we allow for arbitrary $M \in \KdualCNzero$ instead of just $M=0$ in (\ref{eq:Kext:rule}) has the advantage that $\Kext R$ satisfies \cite[Ass.~3.7]{BrunedChandraChevyrecHairer2017} as soon as $R$ satisfies this assumption. This simply ends up ensuring that one can build arbitrary products of $U_\ft$ for any $\ft\in \Kext\FL_+$ with $\reg(\ft)>0$. As was already remarked below \cite[Ass.~3.7]{BrunedChandraChevyrecHairer2017}, this is not a restriction at all, since any subcritical rule can be trivially extended in such a way that this assumption holds, and we will assume from now that \cite[Ass.~3.7]{BrunedChandraChevyrecHairer2017} holds for $R$ and thus also for $\Kext R$.
\end{remark}

\begin{remark}
It might appear more natural to set $\Kext R( \Kdual \ft ) = \Kext R (\ft)$ for any $\ft \in \FL_+$, in which case $\ft$ and $\Kdual \ft$ would be end up to be completely interchangeable in the extended regularity structure constructed from $\Kext R$. The present formulation is more restrictive and leads to a smaller regularity structure, but as we shall see, this structure is still rich enough to lift the equation (\ref{eq:dual:equation}) for $w$ to an abstract fixed point problem. The reason we choose the present formulation is that the natural sector $\CW$ in which the solution to this abstract fixed point problem takes values in is function-like, compare Lemma~\ref{lem:sector:dual:function:like} below.
\end{remark}

\begin{example}
Continuing Example \ref{ex:rule} of the stochastic heat equation, the set $\Kext R(\ft)$ is equal to $\Kext R(\Kdual \ft)$ for any $\ft \in \FL_+$ and contains those multisets $\fm \in \N^{\Kext \Psi}$ with the property that $\fm(\Xi,0)\le 1$ and $\fm(\fl,k) = 0$ for any $\fl \in \FL$ and $k \in \N^{d+1} \backslash\{0\}$.

An example where $\Kext R(\ft)$ and $\Kext R(\Kdual\ft)$ do not coincide is given by the $\Phi^4_3$ equation (\ref{eq:example:phi}). In this case, writing again $\FL_+=\{\ft\}$ and $\FL_-=\{\Xi\}$, the set $\Kext R(\ft)$ is given by the set of all multisets which can be written as
\begin{align*}
\emptyset,\, [(\Xi,0)],\, 
[(\Kdual\ft,0),\ldots],\, 
[(\ft,0),(\Kdual\ft,0),\ldots],\,
\ldots
,\,\text{ or }\,
[(\ft,0),(\ft,0),(\ft,0), (\Kdual\ft,0),\ldots],
\end{align*}
where $(\Kdual\ft,0),\ldots$ stands for an arbitrary number of types $(\Kdual\ft,0)$. On the other hand, $\Kext R(\Kdual\ft)$ is given by all multisets which can be written as
\[
\emptyset,\,
[(\Kdual\ft,0),\ldots],\,
[(\ft,0),(\Kdual\ft,0),\ldots],
\,\text{ or }\,
[(\ft,0),(\ft,0),     (\Kdual\ft,0),\ldots].
\]
\end{example}

In order to continue, we recall from \cite[Rem.~5.17]{BrunedHairerZambotti2016} that given $\tilde \theta>0$ one can assume without loss of generality that the function $\reg$ satisfies the bound 
\begin{align}\label{eq:reg:optimal}
\reg(\ft) > \min\{ |\tau|_\fs : \tau \in \TT_\ft^F \backslash \bar\CT\} + |\ft|_\fs - \tilde \theta
\end{align}
for any $\ft \in \FL_+$, where $\TT_\ft^F$ denotes the set of trees $\tau$ such that $\CJ_\ft^0\tau \in \TT$ and $\tau$ is $\ft$-non-vanishing as in Section \ref{sec:singularSPDEs}.

\begin{remark}
There are some subtleties here, since in \cite[Rem.~5.17]{BrunedHairerZambotti2016} this identity was only shown with $\TT_\ft^F$ replaced by the larger set $\TT_\ft$. In general (\ref{eq:reg:optimal}) might simply not be true, since the rule might be chosen larger then necessary to deal with the singular SPDE at hand. However, this problem can easily be circumvented by assuming without loss of generality that $R$ is given as the completion of the "naive" rule $R_\naive$, which is defined 
in such a way that the set of trees $\tau \in \TT_\ft$ that strongly conform to $R_\naive$ coincide with the set of trees $\tau \in \TT_\ft$ that are $\ft$-non-vanishing. One can then apply \cite[Rem.~5.17]{BrunedHairerZambotti2016} to $R_\naive$ in order to obtain (\ref{eq:reg:optimal}).
\end{remark}

We assume from now on that (\ref{eq:reg:optimal}) holds for some $\tilde\theta>0$ small enough (to be determined later), and with this convention we have the following lemma.
\begin{lemma}\label{lem:Kext}
Assume that $\theta>0$ is small enough and that (\ref{eq:reg:optimal}) holds. Then the rule $\Kext R$ is a subcritical rule with respect to $\reg$. In particular, there exists a subcritical completion of $\Kext R$ defined via \cite[Prop.~5.21]{BrunedHairerZambotti2016}, which we again denote by $\Kext R$, and we can define the extended regularity structure $\Kext\CT^\ex$ as in \cite[Sec.~5.5]{BrunedHairerZambotti2016}. 
\end{lemma}
\begin{proof}
In order to see that $\Kext R$ is subcritical, note first that for $\ft \in \FL_+$ one has 
\[
\inf_{N \in \Kext R(\ft)} \reg(N) = \inf_{N \in R(\ft)} \reg(N).
\]
Let now $\ft \in \FL_+$ and $N \in R(\Kdual \ft)$, and let $\fl \in \FL_+$ and $k\in N^{d+1}$ such that $N \sqcup \{(\fl,k)\} \in R(\ft)$.
By (\ref{eq:reg:optimal}) we can choose for any kernel-type $\fj \in \FL_+$ a tree $\tau_{\fj} \in \TT_{\fj}^F\backslash{\bar\CT}$ such that $\reg( \fj ) > | \tau_{ \fj } | _\fs + | \fj |_\fs - \tilde\theta$ for some $\tilde\theta>0$ small enough. We now consider the tree
\[
\tau := \prod_{(\fj,\alpha) \in N \sqcup \{ (\fl,k) \}}
	\CJ_{\fj}^\alpha [ \tau_\fj ].
\]
It follows that $\tau$ strongly conforms to $R$ (c.f. \cite[Def.~5.8]{BrunedHairerZambotti2016}) and it follows from the definition of $\CT^\ex$ in \cite[Def.~5.26]{BrunedHairerZambotti2016} that $\tau \in \TT_\ft^F$. 
We also define the subtree $\tilde\tau$ of $\tau$ by setting
\[
\tilde \tau := \prod_{(\fj,\alpha) \in N }
	\CJ_{\fj}^\alpha [ \tau_\fj ].
\]
Then $\tilde\tau\in \TT_\ft^F$ is a proper subtree of $\tau$ with identically root, and trees like this satisfy $|\tilde\tau|_\fs> - |\ft|_\fs + 2\theta$ for $\theta>0$ small enough by  \cite[Ass.~2.10]{BrunedChandraChevyrecHairer2017}.

On the other hand, provided that $\tilde\theta$ is smaller then $\frac{\theta}{\# N}$, one has $\reg(N)> |\tilde\tau|_\fs - \theta> -|\ft|_\fs+\theta = -|\ft|_\fs + \reg(\Kdual\ft)$, and this concludes the proof.
\end{proof}

For $X\in \{G,K,R\}$ we write $X_{\Kdual\ft}(z):= X_{\ft}(-z)$, so that in particular $G_{\Kdual\ft}$ is the Greens function for $-\partial_\ft-\CL_\ft^*$ and the compactly supported kernels $K_{\Kdual\ft}$ satisfy \cite[Ass.~5.1 Ass.~5.4]{Hairer2014}. Given the kernel assignment $K_\fl$ for any ${\fl \in \Kext\FL_+}$ we write $\Kext\CM_\infty$ for the set of smooth, reduced, admissible models for $\Kext\CT^\ex$ and we write $\Kext\CM_0$ for the closure of this set. For $f \in \Omega_\infty(\FL_-)$ and $\eta \in \SM_\infty(\FL_-)$ we write $\Kext Z^f$ and $\Kext{\hat Z}^\eta_\BPHZ$ for the canonical lift of $f$ and the BPHZ-renormalized canonical lift of $\eta$, respectively, defined as usual via \cite[Rem.~6.12]{BrunedHairerZambotti2016} and \cite[Thm.~6.17]{BrunedHairerZambotti2016}. We also write $\Kext\CD^{\gamma}$ for the space of modelled distributions as in Section \ref{sec:modelled:dist} with $\CT^\ex$ replaced by $\Kext\CT^\ex$. We will later on need to work with modelled distributions that are only defined on a domain of the form $(\theta,T)\times \T^d$ for some $T>\theta>0$, and we write $\Kext \CD^{\gamma,(\theta,T)}$ for the space of functions $W:(\theta,T) \times \T^d \to \Kext \CT_{<\gamma}$ that satisfy (\ref{eq:bound:D:gamma:eta}).

\subsection{An Abstract Fixed Point Problem for the Dual Equation}
\label{sec:dual:FPP}

Given $U \in \CD^{\gamma,\eta}_{\sector}$, with $\gamma,\eta,\sector$ as in Proposition \ref{prop:app:singular:SPDEs}, and $\phi \in \bigoplus_{\ft \in \FL_+}\CC_c^\infty(\domain)$ we want to consider the abstract point problem in $\Kext\CD^{\Kdual\gamma}_\CW$ for some families $\Kdual\gamma_{ \Kdual \ft }>0$, for $\ft \in \FL_+$, and some sector $\CW=\bigoplus_{\Kdual\ft \in \Kdual\FL_+}\CW_{\Kdual \ft}$, given by
\begin{align}\label{eq:FPP:dual}
W_{\Kdual\ft} = \CQ_{< \gamma_{\Kdual\ft} } \CP_{ \Kdual\ft}
	 \I_+ 
	 [ DF_\ft(U) W 
		+\sum_{\Xi \in \FL_-} DF_\ft^\Xi(U) W \Xi + \phi_\ft ]
\end{align}
for any $\ft \in \FL_+$. The purpose of this section it to find the right sectors and exponents for this fixed point problem to be well posed.
In order to unify notation, we define
\begin{align}\label{eq:F:dual}
F_{\Kdual \ft}( u , w ) := DF_\ft(u) w
\quad \text{ and } \quad
F_{\Kdual \ft}^\Xi( u , w ) := DF_\ft^\Xi( u ) w,
\end{align}
for any $\ft \in \FL_+$, $\Xi \in \FL_-$ and $u,w \in \domain^{\FL_+}$. With this convention $F_\fl$ and $F_\fl^\Xi$ are well defined for any $\fl \in \Kext\FL_+$ and $\Xi \in \FL_-$. We will sometimes write $F^\bullet_\ft := F_\ft$ to avoid case distinctions.

\begin{remark}\label{rem:index:simplified}
In \cite[Sec.3-6]{BrunedChandraChevyrecHairer2017} the authors were working in a more general setting, in the sense that they allowed for derivatives hitting noises and noises being multiplied together. Additionally, they were considering non-linearities that can depend on the extended decoration $\fo$. 
Our setting can easily be embedded into this more general setting, by defining $\tilde F$ via
\[
\tilde F_\ft^{0,0} := F_\ft
\quad \text{ and } \quad
\tilde F_\ft^{\I_\Xi, 0 }:= F_\ft^\Xi
\]
for any $\ft \in \FL_+$ and $\Xi \in \FL_-$, and $F_\ft^\fl = 0 $ otherwise. Whenever we refer to results from \cite{BrunedChandraChevyrecHairer2017} in the sequel we will make these identification implicitly. 
\end{remark}

In the sequel we will need results of \cite{BrunedChandraChevyrecHairer2017} applied to $(F_\fl)_{\fl \in \Kext \FL_+}$. In order to do so, we need the following technical Lemma.
\begin{lemma}\label{lem:dual:obeys:rule}
Assume that $\theta>0$ is small enough. Then $(F_\fl)_{ \fl \in \Kext \FL_+ }$ obeys $\Kext R$ in the sense of \cite[Def.~3.5]{BrunedChandraChevyrecHairer2017}.
\end{lemma}
\begin{proof}
We use \cite[Prop.~3.8]{BrunedChandraChevyrecHairer2017}.
In our setting, which is a bit simplified compared to \cite{BrunedChandraChevyrecHairer2017}, the second conditions of \cite[Prop.~3.8]{BrunedChandraChevyrecHairer2017} reduces to the statement
\begin{align}
\alpha \sqcup \{ (\Xi,0) \} \notin \Kext R(\fl)
&\implies
D^\alpha F_\fl^\Xi =0 
\\
\alpha  \notin \Kext R(\fl)
&\implies
D^\alpha F_\fl = 0
\end{align}
for any $\fl \in \Kext \FL_+$, $\Xi \in \FL_-$ and $\alpha \in \N^{\Kext \Psi_+}$ where $\Kext \Psi_+:=\Kext \FL_+ \times \N^{d+1}$. For $\ft \in \FL_+$ this follows from the respective assumption on $R(\ft)$ and the definition of $\Kext R(\ft)$ in (\ref{eq:Kext:rule}). 

Let now that $\alpha \in \N^{\Kext \Psi}$ such that $D^\alpha F_{\Kdual\ft} \ne 0$. Assume first that $\alpha \in \N^{\FL_+ \times \N^{d+1}}$. By definition of $F_{\Kdual \ft}$ in (\ref{eq:F:dual}) it follows that there exists $\fl \in \FL_+$ such that
\[
D^{\alpha \sqcup \{(\fl,0)\} } F_\ft \ne 0,
\]
and by the first part of the proof this implies $\alpha \sqcup \{ (\fl, 0) \} \in \Kext R(\ft)$. From (\ref{eq:Kext:rule:dual}) we infer that $\alpha \in \Kext R(\Kdual \ft)$, as required. Assume now that $\alpha \notin \N^{\FL_+ \times \N^{d+1}}$. Then it follows that there exists a (unique) $\fl \in \FL_+$ such that $\alpha_{(\Kdual \fl,0)} \ne 0$, and we can write $\alpha =: \bar \alpha \sqcup \{(\Kdual \fl,0) \}$. It follows that one has
\[
0 \ne D^\alpha F_{\Kdual \ft} = D^{\bar \alpha \sqcup \{ (\fl,0) \}} F_\ft,
\]
and thus $\bar \alpha \sqcup \{ (\fl,0) \} \in \Kext R(\ft)$. Again by (\ref{eq:Kext:rule:dual}) we infer that $\bar \alpha \in \Kext R(\Kdual \ft)$, and since $\Kext R$ satisfies \cite[Ass.~3.7]{BrunedChandraChevyrecHairer2017} (or directly from the definition), we infer that $\alpha = \bar \alpha \sqcup \{ (\Kdual \fl,0)\}\in \Kext R(\Kdual \ft)$, as required.

The claim concerning $F_\fl^\Xi$ for any $\fl \in \Kext \FL_+$ and $\Xi \in \FL_-$ follows in the same way.
\end{proof}

In analogue to (\ref{eq:ft:non:vanishing}), given $\fl \in \Kext \FL_+$ we say that a tree $\tau = (T^{\fn,\fo}_\fe,\ft) \in \Kext \TT^\ex$ is $\fl$-non-vanishing for $F$ if 
$$\big( 
	\partial^{\fn(\rho_\tau)} \prod_{ e \in K(\tau), e^\downarrow = \rho_\tau} D_{(\ft(e), \fe(e))}
\big)
 F^\Xi_\ft \ne 0$$
 and $\tau_e$ is $\ft(e)$-non-vanishing for any $e \in K(\tau)$ with $e^\downarrow = \rho_\tau$.
Here we set $\Xi := \ft(e)$ if there exists a (necessarily unique) noise-type edge $e \in K(\tau)$ with $e^\downarrow = \rho_\tau$ and $\Xi := \bullet$ otherwise, and $\tau_e$ denotes the largest sub-tree of $\tau$ with root $e^\uparrow$.
We define $\Kext \TT_\fl^F$, $\Kext {\tilde \TT}_\fl^F$ and $\Kext \TT_{\fl,-}^F$ in analogue to Section \ref{sec:singularSPDEs}, so that one has
\[
{\Kext \TT}_\fl^F := 
\{
\tau \in \Kext \TT :
\tau \text{ is } \fl \text{ - non vanishing for } F
\text{ and } \CJ_{(\fl,0)}[\tau] \in \Kext \TT
\}
\]
(note that in particular $\Kext \TT_\ft^F = \TT^F_\ft$ for any $\ft \in \FL_+$), and the sets $\Kext {\tilde\TT}_\fl^F$ and $\Kext \TT_{\fl,-}^F$ consist of those trees $\tau \in \Kext \TT_{\fl}^F$ such that $|\tau|_\fs \le 0$ and $|\tau|_\fs <0$, respectively. We also set
$\Kext\CT^F_\fl := \linspace {\Kext \TT_\fl^F}$ for any $\fl \in \Kext \FL_+$. With this notation we define for any $\fl \in \Kdual\FL_+$ the sectors
\begin{align}\label{eq:sectors:W}
\CW_\fl := \bar\CT^\ex \oplus  \CJ_{(\fl,0)}[\Kext\CT^F_\fl]
\qquad \text{ and } \qquad
\bar \CW_\fl := \bar\CT^\ex \oplus \Kext\CT^F_\fl.
\end{align}

We write similar to above $\CW:=\bigoplus_{\fl \in \Kdual\FL_+} \CW_\fl$ and $\bar\CW:=\bigoplus_{\fl \in \Kdual\FL_+} \bar\CW_\fl$. 
We now have the following analogue to \cite[Lem.~6.9]{BrunedChandraChevyrecHairer2017}.
\begin{lemma}\label{lem:dual:nonlinearity:sectors}
For any $\fl \in \Kext\FL_+$ the spaces $\CW_\fl$ and $\bar \CW_\fl$ form sectors in $\Kext\CT^\ex$. Moreover, for any $U \in \CV$ and $W \in \CW$ and any $\Xi \in \FL_-$ one has that 
\[
 DF_\ft(U) W 
\qquad
\text{ and }
\qquad
DF_\ft^\Xi(U) W \Xi
\]
are elements of $\bar\CW_{\Kdual \ft}$ for any $\ft \in \FL_+$.
\end{lemma}
\begin{proof}
This is the content of \cite[Lem.~6.9]{BrunedChandraChevyrecHairer2017}. 
\end{proof}

In the sequel we need to understand structure of the sets $\Kext \TT^F _{\Kdual \ft}$ for $\ft \in \FL_+$. For this we introduce the following notation. Given a tree $\tau= (T^{\fn,\fo}_\fe,\ft) \in \CT^\ex$ and a node $u \in N(T)$ we write 
\[
\dualTree{\tau}{u} := (T^{\fn,\fo}_\fe , \tilde \ft_u),
\]
where $\tilde \ft_u : E(T) \to \Kext\FL$ is given by 
\begin{align}\label{eq:type:dual:tree}
\tilde \ft_u (e)
:=
\begin{cases}
\Kdual{\ft(e)}		\quad & \text{ if }e\text{ lies on the shortest path from }u\text{ to }\rho_T \\
\ft(e)				\quad & \text{ otherwise.}
\end{cases}
\end{align}
It follows from Lemma \ref{lem:Kext} that one has $\dualTree{\tau}{u} \in \Kext\CT^\ex$ for any $\tau \in \CT^\ex$ and any $u \in N(\tau)$. Given additionally an edge $e \in K(\tau)$ with $e^\downarrow = u$, then we write $\dualTreeEdge{\tau}{u}{e}$ for the tree obtained from $\dualTree{\tau}{u}$ by removing $e$ from the edge set, and removing furthermore all edges $\tilde e \in E(\dualTree{\tau}{u})$ and vertices $\tilde u \in V(\dualTree{\tau}{u})$ with the property that $e$ lies on the shortest path from $\tilde e$ respectively $\tilde u$ to the root $\rho_\tau$. It is clear that one obtains another decorated, typed tree in this way by simply restricting the corresponding maps to $N(\dualTreeEdge{\tau}{u}{e})$ and $E(\dualTreeEdge{\tau}{u}{e})$, respectively, and since $\Kext R$ is a normal rule, one has $\dualTreeEdge{\tau}{u}{e} \in \Kext\CT^\ex$.
We now have the following Lemma.

\begin{lemma}\label{lem:ext:set:t-non-vnaishing}
Assume that Assumption \ref{ass:dual:equation:simplicitiy} holds. Then for any $\ft \in \FL_+$ the set $\Kext \TT^F _{\Kdual \ft}$ agrees with the set of trees $\dualTreeEdge{\tau}{u}{e}$ for $\tau \in \TT^F_\ft$, $u \in N(\tau)$ and $e \in K(\tau)$ such that $e^\downarrow = u$.
\end{lemma}
\begin{proof}
Let first $\tau \in \TT^{F}_{\ft}$ and let $u \in N(\tau)$ and $f \in K(\tau)$ be such that $f^\downarrow = u$. It follows from the definition of $\Kext R$ in (\ref{eq:Kext:rule:dual}) that $\CJ_{\ft}^0[\tau] \in \TT^\ex$ implies $\CJ_{\Kdual \ft}^0[\dualTree{\tau}{u}] \in \Kext \TT^\ex$ and by completeness of the rule $\Kext R$ one also has $\CJ_{\Kdual \ft}^0[\dualTreeEdge{\tau}{u}{f}] \in \Kext \TT^\ex$. It thus remains to show that $\dualTreeEdge{\tau}{u}{f}$ is $\Kdual\ft$-non-vanishing. Proceeding inductively in the number of edges of $\tau$, it suffices to show that (\ref{eq:ft:non:vanishing}) does not vanish identically for the root $\tilde\rho:=\rho(\dualTreeEdge{\tau}{u}{f})$. For this let $\tilde \ft_u$ be as in (\ref{eq:type:dual:tree}) and let $\CE$ and $\tilde\CE$ denote the set of edges $e \in K(\tau)$ and $e \in K(\dualTreeEdge{\tau}{u}{f})$ such that $e^\downarrow = \rho_\tau$ and $e^\downarrow = \tilde\rho$, respectively. In case that $u = \rho_\tau$, one has $\tilde\ft_u = \ft$ and by definition of $\dualTreeEdge{\tau}{u}{f}$ one has $\CE = \tilde\CE \sqcup \{f\}$, so that it follows that
\begin{multline}\label{eq:ext:set:1}
\partial^{\fn(u)} 
(\prod _{ e \in \tilde\CE } D_{ ( \ft(e) , \fe(e) ) })
	F_{\Kdual \ft}^{ \Xi[u] }(u, \nabla u , \ldots;w, \nabla w, \ldots)
\\=
\sum_{(\fl,k) \in \FL_+ \times \N^{d+1}}
\partial^{\fn(u)} 
(\prod _{ e \in \tilde\CE } D_{ ( \ft(e) , \fe(e) ) })
	D_{(\fl,k)}F_{\ft}^{ \Xi[u] }(u, \nabla u , \ldots) \partial_k w_{\fl}
\end{multline}
Since this expression is linear in $\partial_k w_\fl$, in order to see that this expression does not vanish identically, it suffices to find one pair $(\fl,k) \in \FL_+ \times \N^{d+1}$ such that the coefficient of $\partial_l w_\fl$ is non vanishing. We choose $(\fl,k)=(\ft(f),\fe(f))$ and we note that the corresponding coefficient in (\ref{eq:ext:set:1}) is equal to
\[
\partial^{\fn(u)} 
(\prod _{ e \in \CE } D_{ ( \ft(e) , \fe(e) ) })
	F_{\ft}^{ \Xi[u] }(u, \nabla u , \ldots)
\]
which does not vanish identically by assumption.
In case $u \ne \rho_\tau$ one has $\CE = \tilde \CE$, and there exists a unique edge $\bar e \in \CE$ such that $\bar e$ lies on the unique shortest path from $\rho_\tau$ to $u$. It follows that $\tilde \ft_u(\bar e) = \Kdual \ft(\bar e)$ and $\ft(e) = \tilde\ft_u(e)$ for any $e \in \CE\backslash\{ \bar f \}$, and using the fact that $D_{(\Kdual\fl,k)} F_{\Kdual\ft}^\Xi(u, \nabla u, \ldots; w, \nabla w, \ldots) = D_{ (\fl,k) } F^\Xi_\ft(u, \nabla u, \ldots)$ for any $\ft,\fl \in \FL_+$, we obtain
\begin{multline}\label{eq:ext:set:2}
\partial^{\fn(u)} 
(\prod _{ e \in \CE } D_{ ( \ft(e) , \fe(e) ) })
	F_{\Kdual \ft}^{ \Xi[u] }(u, \nabla u , \ldots;w, \nabla w, \ldots)
\\=
\partial^{\fn(u)} 
(\prod _{ e \in \CE\backslash\{ \bar e \} } D_{ ( \ft(e) , \fe(e) ) })
	D_{ (\ft(\bar e), \fe(\bar e) ) }F_{ \ft}^{ \Xi[u] }(u, \nabla u)
\\=
\partial^{\fn(u)} 
(\prod _{ e \in \CE } D_{ ( \ft(e) , \fe(e) ) })
	F_{\ft}^{ \Xi[u] }(u, \nabla u , \ldots),
\end{multline}
which does not vanish identically by assumption.

Conversely, let $\sigma=(S^\fn_\fe,\ft) \in \Kext{ \TT }^F_{\Kdual \ft}$. It follows from the fact that $F_{\Kdual\fl}^\Xi$ is linear in $(w, \nabla w, \ldots)$ for any $\fl \in \FL_+$ and $\Xi \in \FL_-$ that there exists a (unique) vertex $\mu \in N(\sigma)$ such that $\ft(e) \in \Kdual\FL_+$ if and only if $e$ lies on the unique shortest path from $\rho_\sigma$ to $\mu$. Let $\CE$ be the set of edges $e \in K(\sigma)$ such that $e^\downarrow = \mu$, and define $\fj \in \Kext\FL_+$ by setting $\fj := \ft(u^\downarrow)$ if $u \ne \rho_\sigma$, and $\fj := \Kdual\ft$ otherwise. By definition of the rule $\Kext R$ in (\ref{eq:Kext:rule:dual}) it follows that there exists $(\fl,k) \in \FL_+ \times \N^{d+1}$ such that one has
\footnote{Recall the notation $[\cdot,\cdot]$ for multisets from (\ref{eq:notation:multiset}).}
\[
[\CE,(\ft,\fe)] \sqcup \{ (\fl,k) \} \in \Kext R(\fj).
\]
Choose an arbitrary tree $\tilde \tau \in \TT_{\fl}^F$ and define now the typed, decorated tree $(T^{\tilde \fn}_{\tilde \fe},\fl)$ by connecting $\rho(\tilde \tau)$ to $\mu$ via an edge $\bar e$ such that $\fl(\bar e) = \fl$ and $\tilde \fe(\bar e) = k$, and where $\tilde \fn$, $\tilde \fe$ and $\fl$ extend the decorations and type-maps of $\sigma$ and $\tilde\tau$ otherwise. It then follows that 
$\tau=(T^{\tilde\fn}_{\tilde\fe},\cq\fl) \in \TT_{\ft}$ and one has $\sigma = \dualTreeEdge{\tau}{\mu }{{\bar e}}$. The fact that $\tau$ if $\ft$-non-vanishing follows by reversing the arguments of the first part of the proof.
\end{proof}

A particular consequence of Lemma \ref{lem:ext:set:t-non-vnaishing} is that we can give a direct proof of the fact that the sectors $\CW_{\Kdual\ft}$ are function like. Note that such a statement would also follow directly from the analysis \cite{BrunedChandraChevyrecHairer2017} and the fact that $\reg(\Kdual\ft)>0$.

\begin{lemma}\label{lem:sector:dual:function:like}
For any $\ft \in \FL_+$ and any $\tau \in \Kext \TT^F_ {\Kdual \ft}$ one has $|\tau|_\fs>- (|\ft|_\fs \lor \fs_0)$. In particular, the regularity of the sector $\bar\CW_{\Kdual \ft}$ is larger then $-(|\ft|_\fs \lor \fs_0)$, and the sector $\CW_{\Kdual\ft}$ is function-like.
\end{lemma}
\begin{proof}
This is a direct consequence of Lemma \ref{lem:ext:set:t-non-vnaishing}, the fact that $|\dualTree{\tau}{u}|_\fs = |\tau|_\fs$ and \cite[Ass.~2.10]{BrunedChandraChevyrecHairer2017}.
\end{proof}

For $\ft \in \FL_+$ let $\alpha_{\Kdual \ft}\le 0$ denote the regularity of the sector $\bar \CW_{\Kdual\ft}$ and for $\Kdual\gamma>0$ let $\Kdual\gamma_{\Kdual\ft}:= \alpha_{\Kdual\ft}+\gamma+|\ft|_\fs$.

\begin{corollary}\label{cor:dual:FPP:existence}
Assume that $\gamma>0$ is large enough such that $\gamma_\ft>\Kdual \gamma_{\Kdual\ft}>0$ for any $\ft \in \FL_+$. Then for any $\theta>0$ and any $U \in \CD^{\gamma,\eta,T+\theta}_\CV$ and any tupel $\phi \in \bigoplus_{\ft \in \FL_+} \CC_c^\infty((0,T) \times \T^d)$, the fixed point problem (\ref{eq:FPP:dual}) has a unique solution $W \in \Kext \CD^{\Kdual\gamma, (\theta,T)}_\CW$.
\end{corollary}

\begin{proof}
We first note that as a corollary from the proof of  \cite[Lem.~6.9]{BrunedChandraChevyrecHairer2017}, in particular \cite[(6.15)]{BrunedChandraChevyrecHairer2017}, it follows that for any $U \in \CD^{\gamma,\eta,T+\theta}_\CV$, any $\fl \in \FL_+$ and any $\Xi \in \FL_-$ one has $\CQ_{\bar\gamma_\ft} (\partial_{\fl} F_\ft)(U)$ and $\CQ_{\bar\gamma_\ft} (\partial_{\fl} F^\Xi_\ft )(U)\Xi$ are elements of $\CD^{\bar\gamma_\ft,(\theta,T)}$.

Moreover, combing Lemma \ref{lem:sector:dual:function:like} and Lemma \ref{lem:dual:nonlinearity:sectors}, it follows that both $\CQ_{\bar\gamma_\ft} (\partial_{\fl} F_\ft)(U)$ and $\CQ_{\bar\gamma_\ft} (\partial_{\fl} F^\Xi_\ft )(U)\Xi$ take values in a sector of regularity bigger then $-|\ft|_\fs$. Consequently, using the results of \cite[Sec.~6]{Hairer2014} (see Proposition \ref{prop:Dgamma:loc:Lipschitz}), one has that 
\[
W \mapsto DF_\ft(U) W 
		+\sum_{\Xi \in \FL_-} DF_\ft^\Xi(U) W \Xi + \phi_\ft
\]
is a locally Lipschitz continuous map from $\Kext\CD^{\Kdual \gamma,(\theta,T)}_\CW$ to $\Kext \CD^{\Kdual \gamma_\ft - |\ft|_\fs + \kappa,(\theta,T)}_{\bar\CW_\ft}$ for some $\kappa>0$ small enough. At this point the unique existence of a solution to (\ref{eq:FPP:dual}) follows directly from \cite[Thm.~7.8]{Hairer2014}.
\end{proof}

\subsection{Identifying the Solution to the Dual Equation}
\label{sec:dual:id}

We fix from now on a regularization $\xi^\eps$ of $\xi$, and we write $\Kext{\hat Z}^\eps_\BPHZ:=\Kext {\hat Z} ^{ \xi^\eps } _\BPHZ$ for any $\eps>0$. We also write $W^{ \Kext { \hat Z } ^\eps _\BPHZ}$ for the solution of (\ref{eq:FPP:dual}) constructed in Corollary \ref{cor:dual:FPP:existence} for the model $\Kext { \hat Z } ^\eps _\BPHZ$ with $U = \bar U + \tilde U$ given as in Section \ref{sec:application:to:SPDEs} (recall that $U \in \CD^{\gamma,\eta}_V$ and $u=\CR U$ is the solution to (\ref{eq:singular:SPDE})).
As above we denote by $\Kext g_\BPHZ^\eps \in \Kext \CG_-^\ex$ the BPHZ-character of $\xi^\eps$ (for the extended regularity structure $\Kext \CT^\ex$) and we let $M^{\Kext g_\BPHZ^\eps} := (\Kext g_\BPHZ^\eps \otimes \Id) \Delta_-^\ex$. 
Our goal is to link the abstract dual equation (\ref{eq:FPP:dual}) to the dual tangent equation (\ref{eq:dual:equation}). In a first step we can use the machinery of \cite{BrunedChandraChevyrecHairer2017} to derive an equation for the reconstructed solution $\CR W_\ft$ to the abstract fixed-point problem (\ref{eq:FPP:dual}). This equation will be automatically of the form (\ref{eq:dual:equation}), but it is a-priori unclear whether the renormalization constants that one obtains in these two ways coincide (or at least differ by something of order 1 in a suitable sense). This however is necessary if we want to take the limit $\eps \to 0$ in the model. 
Thus, in order to continue, we introduce the following assumption that makes sure that the dual renormalization constants are given by what we would naively expect.
\begin{assumption}\label{ass:counter-terms}
For any $\ft \in \FL_+$ one has the identity
\begin{align}\label{eq:ass:counter-terms}
\sum_{ \tau \in \tilde \TT_\ft^F}
\frac{ g^\eps_\BPHZ(\tau) }{S(\tau)}
D\Upsilon_\ft^F[\tau](u, \nabla u)(w, \nabla w)
&=
\sum_{ \tau \in \Kext {\tilde \TT} ^F_{\Kdual \ft } }
	\frac{ \Kext g^\eps_\BPHZ (\tau) }{ S(\tau) } \Upsilon_{ \Kdual \ft } ^{F} [ \tau ]
		(u, \nabla u; w, \nabla w).
\end{align}
\end{assumption}

\begin{remark}\label{rem:c=BPHZ}
The simplicity of Assumption \ref{ass:counter-terms} is the main reason for assuming that $c^\eps_\tau$ is given by the BPHZ character. In general, in order to pass to the limit $\eps\to 0$ in $u^\eps$, one could choose $c_\tau^\eps = (h \circ g^\eps_\BPHZ)(\tau)$ where $h \in \CG_-$ is an arbitrary fixed group element and $\circ$ denotes the group product in the renormalization group $\CG_-$. In order to treat this more general situation, one would need to show that (\ref{eq:ass:counter-terms}) above implies a similar relation with $g^\eps_\BPHZ$ replaced by $h\circ g^\eps_\BPHZ$ and $\Kext g^\eps_\BPHZ$ replaced by $\Kext h \circ \Kext g^\eps_\BPHZ$ for some character $\Kext h \in \Kext \CG_-$ determined by $h$. We refrain from doing so for simplicity.
\end{remark}

With this assumption, the following Proposition is a straight-forward application of \cite[Thm.~6.7]{BrunedChandraChevyrecHairer2017}, which provides a convenient link between reconstructed solutions to abstract fixed point problems (for smooth models) and renormalized random PDEs.

\begin{proposition}\label{prop:w=RW}
Assume that Assumptions \ref{ass:dual:equation:simplicitiy} and \ref{ass:counter-terms} hold. Then for any $\eps>0$ the smooth function $w_{\ft,\phi}^\eps$ given by
\begin{align}\label{eq:w=RW}
w_{\ft,\phi}^\eps := \CR^{\Text {\hat Z}^\eps_\BPHZ} W_\ft^{\Text {\hat Z}^\eps_\BPHZ}
\end{align}
solves (\ref{eq:dual:equation}).
\end{proposition}
\begin{proof}
We are going to apply \cite[Thm.~6.7]{BrunedChandraChevyrecHairer2017}. First note that we are indeed in the setting of this theorem, since by Lemma \ref{lem:dual:obeys:rule} one has that the right hand side of (\ref{eq:FPP:dual}) obeys $\Kext R$, the assumption on $\eta$ follows trivially, since we stay away from the initial time, the condition on $\gamma$ can always be achieved by increasing $\gamma$ if necessary, and the fact that $\I_+ DF_\ft(U)W$ and $\I_+ DF_\ft^\Xi(U)W \Xi$ take values in $\CD^{\bar\gamma,\bar\eta}$ for some $\bar\gamma>0$ and $\bar\eta>-\fs_0$ follows from Corollary \ref{cor:dual:FPP:existence}. Note also that in the proof of \cite[Thm.~6.7]{BrunedChandraChevyrecHairer2017} the equation is derived via its mild formulation and the equation in its derivative form is only obtained in the last step, so that is all arguments go through in the time reversed setting verbatim. Denoting by $\tilde w_{\ft,\phi}^\eps$ the right hand side of (\ref{eq:w=RW}), it follows now from \cite[Thm.~6.7]{BrunedChandraChevyrecHairer2017} that one has
\[
-\partial_t \tilde w_{\ft,\phi}^\eps = \CL_\ft^* \tilde w_{\ft,\phi}^\eps
+
(MF)_{\Kdual \ft}
	(u^\eps, \tilde w_{\phi}^\eps, \nabla \tilde w_{\phi}^\eps, \ldots) 
+
\sum _{ \Xi \in \FL_- } (MF)_{\Kdual \ft} ^\Xi 
	(u^\eps, \tilde w_{\phi}^\eps, \nabla \tilde w_{\phi}^\eps, \ldots) 
		\xi_\Xi^\eps
+\phi_\ft
\]
for any $\ft \in \FL_+$. Here, the function $(MF)_\fl^\Xi$ for ${\fl \in \Kext \FL_+}$ and $\Xi \in \FL_-\sqcup\{\bullet \}$ was defined in \cite[(3.9)]{BrunedChandraChevyrecHairer2017}, and is given by
\begin{align}\label{eq:dual:MF:2}
(MF)_{ \Kdual \ft }
	(u^\eps, \tilde w_{\phi}^\eps, \nabla \tilde w_{\phi}^\eps)
= 
	\sum_{\tau \in \Kext {\tilde \TT}_{\Kdual \ft}^F}
		\frac{ 1 }{ S(\tau) }
		\Kext g_\BPHZ^\eps (\tau)
		\Upsilon_{ \Kdual \ft } ^F[\tau]
		(u^\eps, \tilde w_{\phi}^\eps, \nabla \tilde w_{\phi}^\eps),
\end{align}
while $(MF)_{\Kdual \ft}^\Xi(u^\eps, w^\eps) = F_{\Kdual \ft}^\Xi(u^\eps, w^\eps) = DF_\ft^\Xi(u^\eps) w^\eps$.

On the other hand, the right hand side of (\ref{eq:dual:equation}) can be written in the form
\begin{multline}\label{eq:dual:MF}
\sum_{\Xi \in \FL_-} DF^\Xi_\ft(u^\eps) w^\eps_\phi \xi^\eps_\xi + \phi_\ft
+
DF_\ft(u^\eps)w^\eps_\phi
	+
	\sum_{\tau \in \TT_{\ft,-}^{F,\nond}} c^\eps_\tau 
D\Upsilon_\ft^\tau(u^\eps) w_{\phi}^\eps 
-
\sum_{i=1}^d \sum_{\tau \in \TT_{\ft,-}^{F,i}} c^\eps_\tau \partial_i w_{\ft,\phi}^\eps,
\\
=
\sum_{ \tau \in \tilde \TT_\ft^F}
\frac{ g^\eps_\BPHZ(\tau) }{S(\tau)}
D\Upsilon_\ft^F[\tau]( u^\eps ) ( w^\eps , \nabla w^\eps ).
\end{multline}
We conclude by applying Assumption \ref{ass:counter-terms}.
\end{proof}

Using Proposition \ref{prop:w=RW} and performing the limit $\eps \to 0$ in (\ref{eq:dual:identity:reg}) now gives the following corollary.

\begin{corollary}\label{cor:dual:identity}
Assume that Assumptions \ref{ass:dual:equation:simplicitiy} and \ref{ass:counter-terms} holds, and let $w_{\ft,\phi}:= \CR^{\Text {\hat Z}^\xi _\BPHZ} W_\ft^{\Text {\hat Z}^\xi _\BPHZ}$. Then one has the identity
\begin{align}\label{eq:dual:identity}
\dual{ v_h , \phi }_{L^2} 
&= \dual{ \sum_{\Xi \in \FL_-} F^\Xi(u) h_\Xi  , w_\phi }_{L^2}
+
\dual{v_h(0,\cdot), w_\phi(0, \cdot )}_{L^2(\T^d)}
\end{align}
\end{corollary}

Assumption \ref{ass:counter-terms} is not straight-forward to show in general. However, 
an important special case in which we can show directly that Assumption \ref{ass:counter-terms} holds is the case that we consider only a single equation, that is, in case that $\#\FL_+=1$.
\begin{proposition}\label{prop:dual:single:equation}
Under Assumption \ref{ass:dual:equation:simplicitiy} assume that $\FL_+=\{\ft\}$. Then Assumption \ref{ass:counter-terms} holds.
\end{proposition}

\begin{proof}
With the aid of Lemmas \ref{lem:dual:identify:counterterms} - \ref{lem:dual:identify:constants} below, we obtain the following chain of equalities.
\begin{align}
\sum_{ \tau \in \tilde \TT_\ft^F}
\frac{ g^\eps_\BPHZ(\tau) }{S(\tau)}
D\Upsilon_\ft^F[\tau]
&=
\sum_{ \tau \in \tilde \TT^F_\ft }\frac{ g^\eps_\BPHZ (\tau) }{ S(\tau) }
	\sum_{u \in N(\tau)}  \Upsilon_{\tilde \ft}^F[ \dualTree{\tau}{u} ] 
\\
\label{eq:dual:single:equation:2}
&=
\sum_{ \tau \in \tilde \TT^F_\ft }\frac{1}{ S(\tau) }
	\sum_{u \in N(\tau)}  \Kext g^\eps_\BPHZ(\Phi \dualTree{\tau}{u} )\Upsilon_{\tilde \ft}^F[ \Phi \dualTree{\tau}{u} ] 
\\
&=
\sum_{ \tau \in \Kext {\tilde \TT} ^F_{\Kdual \ft } }
	\frac{ 1 }{ S(\tau) } 
	\Kext g^\eps_\BPHZ (\tau)
	\Upsilon_{ \Kdual \ft } ^{F} [ \tau ].
\end{align}
Note that the summand in (\ref{eq:dual:single:equation:2})  vanishes whenever $\FD_u(\tau) \notin \Kext {\tilde \TT} ^F_{\Kdual \ft }$, and otherwise one has $\Phi \FD_u(\tau) \in \Kext {\tilde \TT} ^F_{\Kdual \ft } \ssq \Kext \CT_-^\ex$ by Lemma \ref{lem:root:shift:well:defined}, so that $\Kext g^\eps_\BPHZ(\Phi \FD_u(\tau))$ is well defined.
\end{proof}

In order to demonstrate that we can also deal with some multi-component equations, we consider the following example.

\begin{example}
Consider a coupled system of $\Phi^4_3$-type equations, given by
\[
\partial_t u_i = \Delta u_i + \sum_{i,j,k} c_{i,j,k} u_i u_j u_k + \xi_i,
\qquad i\le n
\]
on $\R \times \T^3$ for some coefficients $c_{i,j,k} \in \R$ and independent Gaussian space-time white-noises $\xi_i$.
It is then not hard to see that Assumption \ref{ass:counter-terms} holds.
\end{example}

\subsection{Existence of Densities}
\label{sec:existence:densities}

\begin{assumption}\label{ass:density:space-time}
We assume that the smooth functions $F_\ft^\Xi \in \CC^\infty(\R^{\FL_-})$ and the solution~$u$ have the property that $(\sum_{\ft \in \FL_+} F_\ft^\Xi(u) w_\ft)_{\Xi \in \FL_-} \ne 0$ on $(0,\tau)\times \T^d$ for any $w \in \R^{\FL_+}\backslash\{0\}$ almost surely.
\end{assumption}

We now have the following theorem, the proof of which is at this stage a generalization of the proof of \cite[Prop.~5.3]{GassiatLabbe2017}.

\begin{theorem}
Under Assumptions \ref{ass:dual:equation:simplicitiy} to \ref{ass:density:space-time}, assume that additionally the Cameron-Martin space $H^\xi$ is dense in $L^2(\domain)^{\FL_-}$. Let $T>0$ and let $\phi^i$, $i\le n$ be a collection of linearly independent test function $\phi^i \in \CC_c^\infty((0,T)\times \domain)^{\FL_-}$. Then, the $\R^{\FL_+\times [n]}$-valued random variable
\[
( \dual{u_\ft, \phi^1_\ft}, \cdots,  \dual{ u_\ft , \phi^n_\ft}) _{\ft \in \FL_+}
\]
conditioned on the event $\{ T < \tau\}$ admits a density with respect to Lebesgue measure.
\end{theorem}
\begin{proof}
Let $X:=( \dual{u_\ft, \phi^1_\ft}, \cdots,  \dual{ u_\ft , \phi^n_\ft}) _{\ft \in \FL_+}$. By Theorem \ref{thm:main:malliavin} we know that $X$ is locally $H^\xi$-differentiable, so that by the Bouleau-Hirsch criterion, Corollary \ref{cor:Bouleau:Hirsch}, we are left to show that $DX$ is almost surely of full rank. Assume first that $n=1$. Then by Theorem \ref{thm:main:malliavin} and Corollary \ref{cor:dual:identity} one has for any $h \in H^\xi$ the identity
\[
D_h X 
= \dual{v_{h} , \phi} 
= \sum_{\Xi \in \FL_-}
	\dual{ h_\Xi  , \sum_{\ft \in \FL_+} F^\Xi_\ft(u) w_{\phi,\ft} } 
+
\dual{v_h(0,\cdot), w_\phi(0, \cdot )}_{L^2(\T^d)}
.
\]
Using the assumption that $H^\xi$ is dense in $L^2(\domain)^{\FL_-}$, it suffices to show that one has
\[
\sum_{\ft \in \FL_+} F^\Xi_\ft(u) w_{\phi,\ft} \ne 0,
\]
which together with Assumption \ref{ass:density:space-time} is equivalent to showing that 
$
w_{\phi,\ft} \ne 0
$ for at least one $\ft \in \FL_+$. On the other hand, by assumption there exists $\ft \in \FL_+$ such that $\phi_\ft  \ne 0$, and it follows directly from (\ref{eq:FPP:dual}) that $W_{\Kdual \fl} \ne 0$ for at least one $\fl \in \FL_+$. It thus suffices to argue that whenever $W$ is a solution to (\ref{eq:FPP:dual}) on some time interval $(\theta,T)$ such that the reconstruction $\CR W$ vanishes on $(\theta,T)\times \T^d$, then this implies that one also has $W=0$ on $(\theta,T)\times \T^d$. Since $W$ takes values in a function-like sector by Lemma \ref{lem:sector:dual:function:like}, one hat $0=\CR W = \dual{W, \bold 1}$, and thus by \cite[Prop.~3.29]{Hairer2014} it suffices to show that $\dual{W_{\Kdual \ft}, \tau}=0$ for any $\ft \in \FL_+$ and any non-polynomial tree $\tau \in \Kext \TT^\ex \backslash \bar \CT$. Assume this was not the case, and let $\fl \in \FL_+$ and $\hat \tau \in \Kext \TT^\ex \backslash \bar \CT$ be the tree of minimal homogeneity such that $\dual{W_{\Kdual \fl}, \hat \tau} \ne 0$.
It follows in particular from Lemma \ref{lem:sector:dual:function:like} that $DF_\fl(U)$ and $D F_\fl^\Xi(U)\Xi$ take values in a sector of regularity $\alpha_\fl>-|\fl|_\fs$. Plugging this in the fixed point equation (\ref{eq:FPP:dual}) implies that 
\[
\min
\Big\{ 
	|\tau|_\fs : \tau \in \Kext \TT^\ex \backslash \bar \CT 
	\text{ and } 
	\dual{W_{\Kdual \fl}, \tau} \ne 0
\Big\}
 = |\hat\tau|_\fs + \alpha_\fl + |\fl|_\fs > |\hat\tau|_\fs,
\]
which gives the desired contradiction.

The case $n>1$ can readily be reduced to the case $n=1$. To see this, assume that there exists $\phi_\ft^i \in \CC_c^\infty((0,T) \times \domain)^{\FL_-}$ such that $DX$ is not almost surely of full rank. This implies in particular that there exits $\lambda_\ft^i \in \R$ for $\ft \in \FL_+$ and $i\le n$ such that $\lambda$ is not identically zero and
\[
\sum_{i \le n} \sum_{\ft \in \FL_+} \lambda^i_\ft v_{h,\ft}(\phi_\ft^i)=0,
\]
which in turn implies that one has $\dual{ v_h, \psi }=0$ where $\psi_\ft := \sum_{i \le n} \lambda^i_\ft \phi_\ft^i$.
\end{proof}

\appendix

\section{Continuity of Maps between Spaces of Modelled Distributions}

\begin{proposition}\label{prop:Dgamma:loc:Lipschitz}
Let $V, \bar V$ be sectors in $\CT$ of regularity $\alpha,\bar\alpha$ respectively. The one has the following.
\begin{itemize}
\item \emph{Multiplication.} Let $\gamma, \bar \gamma>0$ and $\eta,\bar \eta \in \R$ and let $\hat \gamma:= (\gamma + \bar\alpha)\land(\bar\gamma+\alpha)$ and $\hat\eta:= (\eta + \bar\alpha) \land (\bar\eta+\alpha)\land (\eta + \bar\eta)$. Assume furthermore that $\hat \gamma \ge 0$ and that $(V,W)$ is $\hat \gamma$-regular (c.f. \cite[Def.~4.6]{Hairer2017}) and denote by $\star:\CT \to \CT$ the tree product. Then one has
\[
\star : \CD^{\gamma,\eta}_V \times \CD^{\bar\gamma,\bar\eta}_{\bar V} \to \CD^{\hat\gamma, \hat\eta}_{V\star\bar V}
\]
is locally Lipschitz continuous.
\item \emph{Differentiation.} Let $i\le d$ and let $\gamma>\fs_i$. Then
\[
\CD_i : \CD_{V}^{\gamma,\eta} \to \CD^{\gamma-\fs_i, \eta-\fs_i}_{W}
\]
with $W:=\{\CD \tau:\tau \in V \}$ is locally Lipschitz continuous.
\item \emph{Integration.} Let $\ft\in\FL_+$ and assume that $\eta<\gamma$ and $\eta \land \alpha>-|\ft|_\fs$. Then if $\gamma+ |\ft|_\fs\notin \N$ and $\eta + |\ft|_\fs \notin \N$ one has
\[
\CP_\ft : \CD^{\gamma,\eta}_V \to \CD^{\gamma + |\ft|_\fs, (\eta\land\alpha)+|\ft|_\fs}_{W}
\]
with $W:=\{\CP_\ft\tau:\tau \in V \}$ is locally Lipschitz continuous.
\item \emph{Composition with smooth functions.} Let $V$ be function-like and  $\gamma$-regular and let $F:\R^n \to \R$ be smooth. Define $\hat F:V_1 \times \ldots V_n \to V$ (c.f. \cite[(4.2)]{Hairer2017}) by	
\begin{align}\label{eq:composition:smooth:function}
\hat F(\tau):= 
\sum_{\alpha \in \N^n}
\frac{D^\alpha F (\langle\tau, \bold 1\rangle)}{\alpha!}
	(\tau - \langle \tau, \bold 1 \rangle)^{\star \alpha}.
\end{align}
Then if $0 \le \eta \le \gamma$ one has that
\[
\CQ_{<\gamma} \hat F: \bigoplus_{i \le n} \CD^{\gamma,\eta}_{V_i} \to \CD^{\gamma,\eta}_W
\]
for some sector $W$ is locally Lipschitz continuous.
\end{itemize}
\end{proposition}
\begin{proof}
See \cite[Sec.~6]{Hairer2017}.
\end{proof}

\section{Lift of the Abstract Fixed Point Problems coming from singular SPDEs}
\label{sec:application:to:SPDEs}

We show in this section that the right hand side of the fixed point problems considered in  \cite{BrunedChandraChevyrecHairer2017} admit $C^1$ lifts to the extended regularity structure. In order to state our results in a clean way, we introduce some notation from \cite{BrunedChandraChevyrecHairer2017}.
To begin with, we fix the subspace $\CM_{0,1}$ of $\CM_0$ considered in \cite[Def.~6.1]{BrunedChandraChevyrecHairer2017} with metric given by
\[
\fancynorm{Z, \bar Z} = \sup_{n \in \N}	
	\frac{1}{n^2+1}\fancynorm{Z, \bar Z}_{[-n-1,n+1]\times \T^d}.
\]
It is clear that this metric is stronger then $\fancynorm{\cdot;\cdot}_{\gamma,K}$ for any $\gamma>0$ and $K\ssq \domain$ compact. From this it follows easily that all statements derived above holds true with $\CM_0$ replaced by $\CM_{0,1}$.  

Let now $\tilde U_\ft \in \CD^{\infty,\infty}_{\sector_\ft}$ be the constant function defined in \cite[(6.10)]{BrunedChandraChevyrecHairer2017} and recall that from \cite[Prop.~6.18]{BrunedChandraChevyrecHairer2017} that $\CP \tilde U := (\CP_\ft \tilde U_\ft)_{\ft \in \FL_+} \in \CD^\infty_{\sector}$ is strongly Lipschitz continuous.

Recall from \cite[Prop.~6.22]{BrunedChandraChevyrecHairer2017} that under Assumption \ref{ass:main} there exists a solution $u$ to the singular SPDE (\ref{eq:singular:SPDE}) and it is given as the reconstruction of a modelled distribution $U \in \CD^{\gamma,\eta}_\sector$, which in turn can be written as $U = \tilde U + \bar U$ with $\tilde U$ as above and $\bar U$ satisfies the fixed point equation
\begin{align}\label{eq:FPP:singular:SPDE}
\bar U_\ft =
\CQ_{<\gamma_\ft} \CP_\ft
\I_+ [ \sum_{\Xi \in \FL_- \sqcup\{\bullet\}} F_\ft^\Xi(\bar U + \CP\tilde U, \nabla (\bar U + \CP\tilde U), \ldots)\Xi
-
\tilde U_\ft]. 
\end{align}
Since $\tilde U$ is a constant $\CD^{\infty,\infty}$ modelled distribution its reconstruction is trivially locally $H^\xi$-Fr\'echet differentiable, and it remains to show that the same is true for $\bar U$. This will follow from the general result of Theorem \ref{thm:local:H:diff}, once we show that the right hand side of (\ref{eq:FPP:singular:SPDE}) admits a $C^1$ lift. 
The main statement that we will show in this section is thus the following.

\begin{proposition} \label{prop:app:singular:SPDEs}
In the notation of \cite[Sec.~6.5]{BrunedChandraChevyrecHairer2017}, let $\gamma_\ft:= \gamma + \reg(\ft)$, $\eta_\ft:=\eta + \ireg(\ft)$, $\bar\gamma_\ft := \gamma_\ft - |\ft|_\fs + \kappa_\ft$ and $\bar \eta_\ft := \eta_\ft + \bar  n_\ft$. Define moreover the sector $V_\ft := \bar\CT^{\ex} \oplus \CI_{(\ft,0)}[\CT_{\ft,+}^F]$ and $\bar V_\ft := \bar\CT^{\ex} \oplus \CT_{\ft,+}^F$. Let finally $H_\ft$ be the non-linearity of \cite[Lem.~6.19]{BrunedChandraChevyrecHairer2017}, given by
\[
H_\ft(U) = \CQ_{<\bar \gamma_\ft}
\sum_{\fl \in \FL_- \sqcup \{0 \} }
	F_\ft^\fl (U + \CP \tilde U)\Xi_\fl - \tilde U_\ft.
\]
Then $H$ is strongly locally Lipschitz and admits a $C^1$ lift.
\end{proposition}

First note that in \cite[Lem.~6.19]{BrunedChandraChevyrecHairer2017} and the discussions below the proof of \cite[Lem.~6.19]{BrunedChandraChevyrecHairer2017} it was shown that the non-linearity $H$ is strongly locally Lipschitz continuous.

Recall from \cite[(3.7)]{BrunedChandraChevyrecHairer2017} that $F_\ft^\fl$ are given as composition with smooth functions $f_\ft^\fl$ so that one has
\[
F_\ft^\fl(U) = \sum_{\alpha \in \N^{\FL_+\times \N^{d+1}}}
	\frac{D^\alpha f^\fl_\ft(\langle \bold U, \bold 1 \rangle )}{\alpha!}
	(\bold U - \langle \bold U, \bold 1 \rangle)^\alpha,
\]
for some $f_\ft^\fl \in \CC^\infty(\R^{\FL_+})$, where we adopted the notation $\bold U = (\bold U_{\ft,k})_{(\ft,k)\in \FL_+ \times \N^{d+1}}$ where $\bold U_{\ft,k}:= D^k U_\ft$, and we write $\langle \bold U, \bold 1 \rangle := (\langle \bold U_{\ft,k}, \bold 1 \rangle)_{(\ft,k) \in \FL_+\times \N^{d+1}}$.
We then have a natural candidate for the lift $\ext H_\ft^{(r)}$ of $H_\ft$ which is given by
\begin{align}
\ext H_\ft^{(r)} (\ext U) := \CQ_{<\bar \gamma_\ft}
	\sum_{\fl \in \CD_\ft} {\ext F}_\ft^\fl(\ext U + \SS_r \CP \tilde U)\SS_r \Xi_\fl - \SS_r \tilde U_\ft.
\end{align}
Here, given an extended model $\ext Z \in \ext \CM_{0,1}$ we write $\ext F_\ft^\fl(\ext U + \SS_r \CP \tilde U)$ for the composition of the smooth function $f^\ft_\fl$ with $\ext U + \SS_r \CP \tilde U$ in the model $\ext Z$.

First note that (\ref{eq:lift}) is a direct consequence of the definitions.
We now sketch the proof that $\ext H_\ft^{(r)}$ is a strongly locally Lipschitz map from $\ext \CD^{\gamma,\eta}_{\ext V}$ into $\ext \CD^{\bar \gamma_\ft,\bar \eta_\ft}_{\ext {\bar V}_\ft}$. Since the proof is very similar to the one given in \cite[Lem.~6.19]{BrunedChandraChevyrecHairer2017}, we will not go into too much detail. The proof essentially boils down to an application of the results of \cite[Sec.~6]{Hairer2014}, which we have summarized in Proposition \ref{prop:Dgamma:loc:Lipschitz}, and the only thing to notice is that the arguments given in the proof of \cite[Lem.~6.19]{BrunedChandraChevyrecHairer2017} carry over to to $\ext{H}_\ft^{(r)}$ verbatim.

The main part of the proof consists in showing that the map $(r,U) \mapsto \ext H^{(r)}_\ft(U)$ is Fr\'echet differentiable. The strategy for this is to strengthen the results of \cite[Sec.~6]{Hairer2014} and show that the respective operations considered there are actually not just Lipschitz continuous but $C^1$.

\begin{proposition}\label{prop:D:gamma:C1}
Under the assumptions of Proposition \ref{prop:Dgamma:loc:Lipschitz} one has that the operations of multiplication, differentiation, integration and composition are $C^1$ between the respective spaces. Moreover, one has the following identity for the Fr\'echet derivative of the operation of composition with smooth functions
\begin{align}\label{eq:composition:derivative}
D (\hat F)(U)V
=
\sum_{k \le n} 
	\CQ_{<\gamma}
	\widehat{D^{e_k} F} (U)
	V_k
\end{align}
where $e_k$ denotes the $k$-th unit vector in $\R^n$ and $V_k$ denotes the $k$-th component of $V$,
\end{proposition}
\begin{proof}
The fact that multiplication, differentiation and integration are $C^\infty$ follows simply from the fact that these operations are continuous and (multi-)linear.

We now show Fr\'echet differentiability of composition with smooth functions, together with (\ref{eq:composition:derivative}). For this we write
\begin{multline}\label{eq:composition:1}
\hat F(V) - \hat F(U) - \sum_{k \le n}\widehat{D^{e_k} F}(U)(V_k-U_k)\\
=
\sum_{\alpha} \frac{1}{\alpha!}
\Big(
D^\alpha F(\bar V) \tilde V^\alpha
-
D^\alpha F(\bar U) \tilde U^\alpha
-
\sum_{k \le n}  D^{\alpha + e_k} F(\bar U) \tilde U^\alpha (V-U)
\Big)
\end{multline}
where we write $\bar U:= \langle U,\bold 1, \rangle$ and $\tilde U:= U - \bar U$. We now define for $x \in \domain$ and $\alpha \in \N^{d+1}$ the function $g_\alpha$ by setting
\[
g_\alpha (t) := D^\alpha F(\bar U + t \bar W)(\tilde U + t \tilde W)^\alpha
\]
with $W:=V-U$. A direct computation shows that (\ref{eq:composition:1}) can be re-written into
\begin{align*}
\sum_{\alpha} \frac{1}{\alpha!}
(g_\alpha(1) - g_\alpha(0) - g_\alpha'(0))
&=
\sum_{\alpha} \frac{1}{\alpha!}
\int_0^1 (1-t) g_\alpha''(t) dt.
\end{align*}
Now note that one has
\[
\sum_{\alpha} \frac{1}{\alpha!} g_\alpha''(t)
=
\sum_{\alpha} \sum_{k,l\le n} \frac{1}{\alpha!} 
	D^{\alpha+e_k + e_l}  F(\bar U + t\bar W) 
	(\tilde U + t \tilde W)^\alpha W^2,
\]
from which we infer that (\ref{eq:composition:1}) can be re-written into
\[
\sum_{k,l\le n}
	\int_0^1 (1-t)	\widehat {F_{k,l}}
		(U + tW) W^2,
\]
where $F_{k,l}:=D^{e_k+e_l} F$. Now using that $\widehat{F_{k,l}}:\CD^{\gamma,\eta} \to \CD^{\gamma,\eta}$ is a Lipschitz continuous map, we can estimate the $\CD^{\gamma,\eta}$ norm of this expression by $\fancynorm{W^2}_{\gamma,\eta} \lesssim \fancynorm{W}_{\gamma,\eta} ^2$, which proves the claim.
\end{proof}

With this proposition, the proof that $\ext H^{(r)}_\ft(\ext U)$ is Fr\'echet differentiable is now straight forward, and consists in redoing the steps of \cite[Lem.~6.19]{BrunedChandraChevyrecHairer2017}, replacing all statements made about Lipschitz continuity with statements about Fr\'echet differentiability.

\section{The Case of a Single Equation}
\label{sec:dual:single:equation}

In the entire section we assume that $\FL_+=\{\ft\}$ contains a unique element. Our goal is to derive the identities necessary to show Proposition \ref{prop:dual:single:equation}.

The first lemma we are going to prove shows that functional derivatives of the original counter-terms are of the same form as the counter-terms of the dual equation.

\begin{lemma}\label{lem:dual:identify:counterterms}
Under Assumption \ref{ass:dual:equation:simplicitiy}, one has for any $\ft \in \FL_+$ and any $\tau \in \tilde \TT_{\ft}^F$ the identity
\begin{align}
D\Upsilon_{ \ft } ^F [ \tau ](u, \nabla u)(w, \nabla w) 
=
\sum_{ \mu \in N(\tau) }
	\Upsilon_{ \Kdual \ft} ^F [\dualTree{\tau}{\mu}]
		( u, w, \nabla w ).
\end{align}
\end{lemma}
\begin{remark}
Note that under Assumption \ref{ass:dual:equation:simplicitiy} the function $D\Upsilon^F_\ft[\tau](u, \nabla u)(w, \nabla w) $ does really only depend on $u, w, \nabla w$.
\end{remark}
\begin{proof}

We only show the statement in the case that $\tau=T^{\fn,\fo}_\fe \in \TT^{F,\nond}_{\ft,-}$. Under Assumption \ref{ass:dual:equation:simplicitiy} the case that $\tau \in \TT^{F,i}_{\ft,-}$ for some $i\in \{ 0, \ldots , d \}$ follows easily.
We first claim that $\fn\equiv 0$. Indeed, otherwise $\Upsilon^F_\ft[\tau]$ contains a factor of the form
\[
\partial^n 
\Big(
\prod_{i\le m} D_{(\ft_i,k_i)} 
\Big)
F_\fl^\Xi(u)
=
\sum_{\fl \in \FL_+}
	D_{(\fl,0)}
	\Big(
	\Big(\partial^{\bar n} 
	\prod_{i\le m} D_{(\ft_i,k_i)} 
	\Big)
	F_\fl^\Xi(u)
	\Big)
\partial_l u_{\fl}
\]
with $n = \bar n + e_l \in \N^{d+1}$ non-zero, where $e_l$ denotes the $l$-th unit vector on $\N^{d+1}$ for some $0 \le l\le d$. But since this factor does not depend on $\partial_l u_\fl$ explicitly by assumption, it must vanish identically, in contradiction to the assumption that $\tau$ is $\ft$-non-vanishing. Moreover, since $F$ only depends on the solution $u$ and not on its derivatives, one also has that $\fe \equiv 0$. Now, denoting by $e_j$ for $j=1,\ldots,n$ the distinct edges $e \in K(\tau)$ with $e^\downarrow = \rho_T$ and by $\tau_j$ the unique maximal subtree of $\tau$ such that $\rho(\tau_j) = e_j^\uparrow$, it follows by definition that one has
\begin{multline}
D\Upsilon_\ft^F[\tau] (u) w
= 
\sum_{i=1}^n
	D \Upsilon_{\ft_i}^F[\tau_i](u)w
	\prod_{j\le n, j\ne i} \Upsilon_{\ft_j}^F[\tau_j] (u)
	\big(
		\prod_{j=1}^n D_{(\ft_j,0)} (u)
	\big)
	F_\ft^\Xi(u)
\\
+
	\prod_{j\le n} \Upsilon_{\ft_j}^F[\tau_j] (u)
	D\big(
		\prod_{j=1}^n D_{(\ft_j,0)} (u)
	\big)
	F_\ft^\Xi(u) w
\end{multline}
where we set $\ft_j:= \ft(\tau_j)$, and $\Xi=\ft(f)$ if there exists a (necessarily unique) edge $f \in L(\tau)$ with $f^\downarrow = \rho(\tau)$, and $\Xi=\bullet$ otherwise.
We proceed inductively in the number of kernel-type edges of $\tau$. For $\# K(\tau) = 0$ the statement holds trivially, so assume from now on that $\# K(\tau) \ge 1$. It then follows from the induction hypothesis that we can write
\[
D\Upsilon_{\ft_i}^F[\tau_i](u)w = 
\sum_{\mu \in N(\tau_i)}
	\Upsilon_{\Kdual{\ft_i}}^F[\dualTree{\tau_i}{\mu}](u,w)
\]
for any $i \le n$. Moreover, by definition of $F_{\Kdual \ft} ^\Xi$ in (\ref{eq:F:dual}), it follows that $D_{(\ft_i,0)} F_\ft^\Xi = D_{(\Kdual{\ft_i},0)}F_{\Kdual \ft}^\Xi$, and hence
\[
\Upsilon_{\Kdual{\ft_i}}^F[\dualTree{\tau_i}{\mu}](u,w)
	\prod_{j\le n, j\ne i} \Upsilon_{\ft_j}^F[\tau_j] (u)
	\big(
		\prod_{j=1}^n D_{(\ft_j,0)}
	\big)
	F_\ft^\Xi(u)
=
\Upsilon_{\Kdual \ft}^F [\dualTree{\tau}{\mu}](u,w).
\]
Since $\bigsqcup_{i\le n} N(\tau_i) = N(\tau) \backslash \{\rho_\tau \}$, it remains to note that $\dualTree{\tau}{{\rho_\tau}}=\tau$ and
\[
\Upsilon_{\Kdual \ft}^F[\tau](u,w)
=
	\prod_{j\le n} \Upsilon_{\ft_j}^F[\tau_j] (u)
	\big(
		\prod_{j=1}^n D_{(\ft_j,0)} (u)
	\big)
	D F_\ft^\Xi(u) w
\]
by definition (\ref{eq:F:dual}).
\end{proof}

Next, we derive a useful identity for the symmetry factors appearing in (\ref{eq:dual:MF}) and (\ref{eq:dual:MF:2}). In order to state it, we introduce the set  $\FD(\tilde \TT_\ft^F):=\{ \dualTree{\tau}{u} : \tau \in \tilde \TT_\ft^F, u \in N(\tau) \}$.
\begin{lemma}\label{lem:symmetry:factors}
Let $G$ be an additive group and for fixed $\ft \in \FL_+$ let $f :  \FD(\tilde \TT_\ft^F) \to G$ be any map such that $f(\tau) = 0$ for any $\tau \in \FD( \tilde \TT_\ft ^F )\backslash \Kext {\tilde \TT} ^F _{ \Kdual \ft }$. Then one has 
\begin{align}\label{eq:symmetry:factors}
\sum_{ \tau \in \tilde \TT^F_\ft } 
	\sum_{ u \in N(\tau) }
	\frac{1}{S(\tau)} 
		f(\dualTree{\tau}{u})
=
\sum_{ \tau \in \Kext{ \tilde \TT}^F_\ft } \frac{1}{S(\tau)}f(\tau).
\end{align}
\end{lemma}
\begin{proof}
Note first that by Lemma \ref{lem:ext:set:t-non-vnaishing} the set $\Kext {\tilde \TT} ^F _{\Kdual \ft}$ is included in the set $\FD(\tilde \TT_\ft^F)$, so that the right hand side of (\ref{eq:symmetry:factors}) makes sense. Since moreover $f$ vanishes outside of $\Kext {\tilde \TT} ^F _{\Kdual \ft}$ by definition, it follows that we can rewrite the left hand side of (\ref{eq:symmetry:factors}) into
\[
\sum_{ \tau \in \Kext{ \tilde \TT}^F_\ft } \frac{m(\tau)}{S(\cq\tau)}f(\tau)
\]
where $m(\tau)\in \N$ is a symmetry factor given by 
\begin{align}\label{eq:identity:m:S}
m(\tau) := \# \{ u \in N(\cq\tau) : \dualTree{\cq\tau} {u} = \tau \}.
\end{align}
It remains to show the identity $m(\tau)S(\tau) = S(\cq\tau)$, which we show inductively in the number of kernel type edges of $\tau$. If $\# K(\tau)=0$ or $\tau = \cq\tau$ the identity is trivial, so that we exclude these cases in the sequel. In case that $\tau = \CJ_\ft^k[\tilde\tau]$ is planted and (\ref{eq:identity:m:S}) holds for $\tilde\tau$, this identity also holds for $\tau$ since $m(\tau) = m(\tilde\tau)$, $S(\tau) = S(\tilde\tau)$ and $S(\cq\tau) = S(\cq\tilde\tau)$. It remains to treat the case that $\cq\tau$ is of the form $\cq\tau = X^k \Xi \prod_{i=1}^n \CJ_{\ft_i}^{k_i}[\cq\tau_i]^{p_i}$ with $n\ge 1$, $p_i \ge 1$ and $(\ft_i,k_i,\tau_i)\ne(\ft_j, k_j, \tau_j)$ for $i\ne j$, $k\in \N^{d+1}$ and $\Xi \in \FL_- \sqcup \{\bullet\}$, and such that (\ref{eq:identity:m:S}) holds for $\CJ_{\ft_i}^{k_i}{\tau_i}$ for any $i \le n$. By assumption there exists $u \in N(\cq\tau)\backslash\{ \rho_\tau \}$ such that $\tau = \dualTree{\cq\tau}{u}$, and we assume without loss of generality that $\tau$ is of the form
\[
\tau = X^k \Xi \CJ_{\Kdual\ft_i}^{k_i}[\tau_1]
	\prod_{i=1}^n \CJ_{\ft_i}^{k_i}[\cq\tau_i]^{p_i- \delta_{i,1}},
\]
which can always be achieved by simply rearranging the order of the triples $(\ft_i,k_i,\tau_i)$. In this case one has
\[
S(\cq\tau)
= k! \prod_{i=1}^m S(\cq\tau_i)^{p_i} p_i!
= k! m(\tau_1)S(\tau_1)\prod_{i=1}^m S(\cq\tau_i)^{p_i - \delta_{i,1}} p_i!
\]

On the other hand, we obtain
\[
S(\tau) = k! S(\tau_1) \prod_{i=1}^m S(\cq\tau_i)^{p_i- \delta_{i,1}} {(p_i - \delta_{i,1})}!,
\]
and the proof is finished, noting that one has the identity $m(\tau) = p_1 m(\tau_1)$.
\end{proof}

Identifying the renormalization constants takes a bit more work. As a preparation, we introduce some notation. Given a rooted tree $T$ with vertex set $V(T)$, edge set $E(T)$ and root $\rho_T$, we can define for any $\nu \in V(T)$ another tree $\Phi_\nu(T)$ with identical edge and vertex sets, but where we set $\rho(\Phi_\nu(T)) := \nu$. Given additionally a type map $\ft : E(T) \to \FL$ and decorations $\fe : E(T) \to \N^{d+1}$ and $\fn : V(T) \to \N^{d+1}$ we obtain another typed, decorated tree $(\Phi_\nu(T)^{\tilde \fn}_\fe,\ft)$ by simply letting the maps $\fe$, $\fn$, and $\ft$ unchanged. For a tree $\tau = T^{\fn}_\fe$ we also define $\hat \Phi_u \tau$ by setting
\begin{align}
\hat\Phi_u \tau
:=
\sum _{ \fm : N(\tau) \to \N^{d+1} } (-1) ^{|\fm|}
\binom{\fn}{\fm}
(\Phi_u T) ^{\fn - \fm + \sum\fm \I_{\rho(\tau)}}_\fe.
\end{align}

Let now $\sigma \in \Kext{ \tilde \TT}_{\Kdual \ft}^F$ be a tree and let $\tau = \cq\sigma \in \tilde\TT_{\ft}^F$. Then by Lemma \ref{lem:ext:set:t-non-vnaishing} there exists $\nu \in N(\tau)$ such that $\sigma = \dualTree{\tau}{\nu}$. Note now that we can naturally identify the node set $N(\tau)$ of $\tau$ with the node set $N(\sigma)$ of $\sigma$. If we do this identification, then the vertex $\nu$ is distinguished in $\sigma$ by the property that it is the unique maximal vertex (with respect to the tree order) that has the property that $\ft(u^\downarrow) \in \Kdual\FL_+$.
In particular, the node $\nu \in N(\sigma)$ is uniquely determined by $\sigma$, and as a consequence it makes sense to define 
\[
\hat\Phi \sigma := \hat\Phi_\nu(\sigma).
\]
The point of this definition is that for any $\tau \in \tilde \TT_{\ft}^F$ and $\nu \in N(\tau)$ we can show an identity between the renormalization constants associated to $\tau$ and $\hat\Phi\dualTree{\tau}{\nu}$, compare Lemma \ref{lem:dual:identify:constants} below.

Before we state any precise statement, we need to deal with the subtlety that it is in general \emph{not} the case that $\Phi_\nu \tau \in \Kext \CT^\ex$ for any $\tau \in \Kext \CT^\ex$ and any $\nu \in N(\tau)$. However, we have the following lemma.

\begin{lemma}\label{lem:root:shift:well:defined}
Assume that Assumption \ref{ass:dual:equation:simplicitiy} holds and that $\FL_+:=\{\ft\}$. Then, the map $\hat\Phi$ is an involutory bijection from $\Kext{ \TT}_{\Kdual \ft}^F$ onto itself, and for any tree $\tau \in \Kext{  \TT}_{\Kdual \ft}^F$ one has the identity $\Upsilon_{\Kdual \ft} ^F[ \Phi \tau ] = \Upsilon_{ \Kdual \ft}^F[ \tau ]$ and $S(\Phi\tau)=S(\tau)$.
\end{lemma}
\begin{proof}
It suffices to show $\Upsilon_{\Kdual \ft} ^F[ \Phi \tau] = \Upsilon_{ \Kdual \ft}^F[ \tau ]$. The fact that $\Phi$ maps $\Kext{  \TT}_{\Kdual \ft}^F$ into itself is then a consequence of the definition of the latter and Lemma \ref{lem:ext:set:t-non-vnaishing}. Moreover, since $\Phi$ is involutory by definition, preserves homogeneity, and the set of $\tau\in\Kext{  \TT}_{\Kdual \ft}^F$ with $|\tau|_\fs = \gamma$ for some fixed $\gamma \in \R$ is finite, it is also a bijection. The identity $S(\Phi\tau)=S(\tau)$ follows easily from the definition.

Concerning the identity  $\Upsilon_{\Kdual \ft} ^F[ \Phi \tau] = \Upsilon_{ \Kdual \ft}^F[ \tau ]$, we use the expression
\[
\Upsilon_{\Kdual \ft}^F [ \tau](u,w)
=
\prod_{\mu \in N(\tau)} 
	\Upsilon_{\Kdual \ft}^F [  \tau , \mu]
:=
\prod_{ \mu \in N(\tau) }
\partial^{\fn(\mu)}
\big(
	\prod_{j=1}^{n[\mu]} D_{(\ft_j[\mu],0)}
\big)
F_{\ft[\mu]}^{\Xi[\mu]}(u,w).
\]
Let $\nu \in N(\tau)$ be the unique vertex such that $\tau = \dualTree{\cq\tau}{\nu}$.
It then follows that for any $\mu \notin \{ \rho_\tau, \nu\}$ the vertex $\mu$ has the same incoming edges (i.e. edges $e\in E(\tau)$ with $e^\downarrow = \mu$) and the same outgoing edge (i.e. the edge $e\in E(\tau)$ with $e^\uparrow = \mu$) when viewed as an element of $N(\tau)$ and $N(\hat\Phi\tau)$, respectively, and moreover the same polynomial label $\fn(\mu)$, so that it follows that one has $\Upsilon_{\Kdual \ft}^F [  \tau , \mu] = \Upsilon_{\Kdual \ft}^F [ \hat\Phi \tau , \mu]$. Moreover, if $e_1, \ldots , e_n$ are the incoming edges of $\nu$ and $e$ is the outgoing edge of $\nu$ when viewed as an element of $N(\tau)$, then $e_1, \ldots , e_n, e$ are the incoming edges of $\nu$ when viewed as an element of $N(\hat\Phi\tau)$, and by construction one has $\ft(e)=\Kdual \ft$. Assume first that $\tau \in \Kext{ \tilde \TT }^{F,\nond}_{\Kdual\ft}$. It then follows that $\fn \equiv 0$, and using the fact that $F_{\Kdual \ft}^{\Xi[\nu]}(u,w)$ is linear in $w_{ \ft}$, so that $(D_{(\Kdual \ft,0)} F_{\Kdual \ft}^{\Xi[\nu]}(u_\ft,w_\ft))w_\ft = F_{\Kdual \ft}^{\Xi[\nu]}(u_\ft, w_\ft)$, we obtain
\begin{align}
(\Upsilon_{\Kdual \ft}^F [ \hat\Phi \tau , \nu](u_\ft)) w_{\ft}
&=
\Big(
\big(
	D_{(\Kdual \ft,0)}\prod_{j=1}^{n[\nu]} D_{(\ft_j[\nu],0)}
\big)
F_{\Kdual \ft}^{\Xi[\nu]}(u_\ft ,w _\ft)
\Big)
w_{\ft}
\\
&=
\big(
	\prod_{j=1}^{n[\nu]} D_{(\ft_j[\nu],0)}
\big)
F_{\Kdual \ft}^{\Xi[\nu]}(u _\ft , w _\ft )
\\
&=
\Upsilon_{\Kdual \ft}^F [ \tau , \nu](u_\ft, w_\ft).
\end{align}
An identical calculation shows that one has $\Upsilon_{\Kdual \ft}^F [ \Phi \tau , \rho_\tau](u_\ft, w_\ft) = (\Upsilon_{\Kdual \ft}^F [  \tau , \rho_\tau](u_\ft))w_\ft$, and this concludes the proof.

Finally, if $\tau \in \Kext{ \tilde \TT}^{F,i}_{\Kdual \ft}$, then one has
\[
\Upsilon_{\Kdual \ft}^F [  \tau , \rho_\tau ](u_\ft, w_\ft) \equiv c
\qquad \text{ and }\qquad
\Upsilon_{\Kdual \ft}^F [  \tau , \nu ](u_\ft, w_\ft) = \partial^{i} w_\ft,
\]
and a similar identity with the roles of $\nu$ and $\rho_\tau$ reversed holds for $\hat\Phi\tau$.
\end{proof}

Finally, we show the following lemma, which is the reason for introducing the map $\hat\Phi$.

\begin{lemma}\label{lem:dual:identify:constants}
For any $\ft \in \FL_+$, any $\tau \in \TT_{\ft,-}^F$ and any node $\nu \in N(\tau)$ with the property that $\FD_\nu(\tau) \in \Kext \TT_{\Kdual\ft,-}^F$ one has the identity
\[
g^\eps_\BPHZ (\tau)
=
\Kext g^\eps _\BPHZ ( \hat\Phi \dualTree{\tau}{\nu} ).
\]
\end{lemma}

We point out that Lemma \ref{lem:dual:identify:constants} does neither require Assumption \ref{ass:dual:equation:simplicitiy} nor the assumption that $\#\FL_+ = 1$.

\begin{proof}
In this proof we use the notation that given a kernel assignment $(L_{\fl})_{\fl \in \Kext\FL_+}$ satisfying \cite[Ass.~5.1, Ass.~5.4]{Hairer2014} and a smooth noise $\eta \in \SM_\infty$, we write $\bold\Pi^{\eta,L}$ and $g_\BPHZ^{\eta,L}$ for the canonical evaluation and the BPHZ character constructed as in \cite[Rmk.~6.12]{BrunedHairerZambotti2016} and \cite[(6.24)]{BrunedHairerZambotti2016} for the kernel assignment $L$ and the noise $\eta$. It follows that if we set $L_\ft := L_{\Kdual\ft}:= K_\ft$ for any $\ft \in \FL_+$, then the effect of $\FD_\nu$ is not seen on the analytic level, and we obtain for any $\tau \in \TT_\ft^F$ and any node $\nu \in N(\tau)$ the identity $\bold\Pi^{\eta,K} \tau = \bold \Pi^{\eta,L}\FD_\nu(\tau)$, and similarly $g_\BPHZ^{\eta,K} \tau = g_\BPHZ^{\eta,L} (\FD_\nu(\tau))$.

For the proof of Lemma \ref{lem:dual:identify:constants}, we are thus left to show that for any $\ft \in \FL_+$ and any $\sigma \in \Kext \TT_{\Kdual\ft,-}^F$ one has the identity
\begin{align}\label{eq:dual:constants:1}
g^{\eta,L}_\BPHZ (\sigma)
=
\Kext g^{\eta,K} _\BPHZ ( \hat\Phi \sigma).
\end{align}
We first deal with the issue that the image of $\Kext \CT^\ex$ under $\hat\Phi$ does in general not coincide with $\Kext \CT^\ex$, which is due to the fact that if $\tau$ is a tree that strongly conforms to the rule $\Kext R$, its image under $\hat\Phi$ might not. We will circumvent this issue by working in the Hopf algebra $\CH_1$ defined in \cite[(4.10)]{BrunedHairerZambotti2016} for the type set $\Kext \FL$. Actually, it suffices for us to work in the reduced Hopf algebra $ \CH$, where $\CH$ is obtained from $\CH_1$ by identifying any trees that only differ by the extended decoration and additionally factoring out any trees $\tau = (T^\fn_\fe,\ft)$ with the property that there exists $e \in E(T)$ such that $\ft(e) \in \FL_-$ and $e^\uparrow$ is either not a leaf or one has $\fn(e^\uparrow)\ne 0$ (or both). Following \cite[Rem.~4.16]{BrunedHairerZambotti2016}, this leads to the following space.
\begin{definition}
We denote by $\CH$ the unital algebra freely generated by typed, rooted, decorated trees $\tau = (T^\fn_\fe,\ft)$ such that $\tau \ne \bullet$ and such that $\ft : \tau \to \Kext \FL$, $\fn : N(\tau) \to \N^{d+1}$ and $\fe : E(\tau) \to \N^{d+1}$, and such that $e^\uparrow$ is a leave of $T$ for any noise type edge $e \in L(\tau)$.
\end{definition}

By \cite[Prop.~3.32]{BrunedHairerZambotti2016}, this space becomes a Hopf algebra when endowed with the co-product $\Delta_1$ defined in \cite[Def.~3.3]{BrunedHairerZambotti2016}. We denote this co-product on $\CH$ simply by~$\Delta$.

\begin{definition}
We define the ideal $\CI_+ \ssq \CH$ generated by all trees $\tau \in \CH$ such that $|\tau|_\fs>0$, and we define the factor algebra $\CH_- := \CH/ \CI_+$, with canonical embedding $\iex : \CH_- \to \CH$.
\end{definition}
It straight forward to see that $\CI_+$ is a Hopf ideal, so that $\CH_-$ is a factor Hopf algebra.
The following Proposition follows exactly as \cite[Prop.~6.5]{BrunedHairerZambotti2016}.
\begin{proposition}
There exists a unique algebra homomorphism $\CA : \CH_- \to \CH$ with the property that 	
\[
\CM(\CA \otimes \Id)\Delta \iex = \bold 1 \bold 1^\star
\]
on $\CH_-$.
\end{proposition}

We now define a subspace $\tilde\CH \ssq \CH$ with the property that $\hat\Phi$ is well defined on $\tilde\CH$ and an involutory bijection. 

\begin{definition}
We define $\tilde\KCH \ssq \KCH$ (respectively $\tilde\KCHhat \ssq \KCHhat$) as the unital sub algebra generated by all trees $\tau \in \KCH$ (respectively $\tau \in \KCHhat$) with the property that there exists a node $u \in N(\tau)$ such that for any edge $e \in E(\tau)$ one has $\ft(e) \in \Kdual{\FL_+}$ if and only if $e$ lies on the unique path from $u$ to the root $\rho_\tau$.
\end{definition}

It is readily checked from the definition of the co-product $\Delta$ and the operation $\Phi$ that $\tilde\KCH$ is closed under $\Delta$, in the sense that $\Delta :\tilde\CH_- \to \tilde\CH_- \otimes \tilde\CH_-$, so that $\tilde\CH_-$ is a Hopf algebra, and $\hat\Phi : \tilde\KCH \to \tilde\KCH$ is such that $\hat\Phi \circ\hat \Phi = \Id$.

On $\KCHhat$ (respectively $\KCH$) we define the character $g^{\eta,K}$ (respectively $g^{\eta,K}_\BPHZ$) by setting $g^{\eta,K}\tau := \E \bold \Pi^{\eta,K}\tau$ (respectively $g^{\eta,K}_\BPHZ(\tau) = g^{\eta,K} (\CA \tau)$) for any tree $\tau$  and extending this linearly and multiplicatively. Note that $\KCT \ssq \CH_-$ and $\KCThat \ssq \CH$, and on these subspaces this notation is consistent, in the sense that one has $\Kext g^{\eta,M}_- = g^{\eta,M}$ and $\Kext g^{\eta,M}_\BPHZ = g^{\eta,M}_\BPHZ$ on $\KCThat$ and $\KCT$, respectively. We are thus left to show that
$
g^{\eta,L}_\BPHZ = g^{\eta,K}_\BPHZ \circ \hat\Phi
$
on $\tilde\CH_-$. We now note that directly from the definition one has the identity $g^{\eta,L} = g^{\eta,K}\circ\cq$ on $\CH$, and since $\hat\Phi$ is an involutory bijection on $\tilde\CH_-$, it follows that it suffices to show 
\begin{align}\label{eq:dual:constants:2}
g^{\eta,K}_\BPHZ \circ \cq \circ \hat\Phi = g^{\eta,K}_\BPHZ 
\end{align}
on $\tilde\CH_-$.
We now identify ideals $\KCI \ssq \KCH$ and $\KCIhat \ssq \KCHhat$ with the property that $g^{\eta,M} _-$ vanishes on $\KCIhat$ for any smooth noise $\eta$ and any kernel assignment $M$, and such that the canonical embedding $\iex$ restricts to an embedding $\iex : \KCI \to \KCIhat$. 

We start with a definition which is completely analogous to the ideal defined in \cite[(2.16)-(2.18)]{Hairer2017} for Feynman diagrams.
\begin{definition}
We denote by $\KCI \ssq \KCH$ (respectively $\KCIhat \ssq \KCHhat$) the ideals generated by all elements which are can be written in form (\ref{eq:ideal:1}), (\ref{eq:ideal:2}), or (\ref{eq:ideal:3}) for some tree $\tau = T^\fn_\fe \in \KCH$ (respectively $\tau \in \KCHhat$), where 
\begin{itemize}
\item For any node $u \in N(\tau) \backslash \{\rho(\tau)\}$ and any $i \le d$ 
\begin{align}\label{eq:ideal:1}
\sum_{ \substack{ e \in E(\tau) \\ e^\downarrow = u } }
	T ^\fn _{\fe + e_i\I_{e}}	
-
T^\fn _{\fe + e_i\I_{u^\downarrow}}
+
\fn(u) T^{\fn - e_i\I_u}_\fe
\end{align}
\item For any $i \le d$ 
\begin{align}\label{eq:ideal:2}
\sum_{ \substack{ e \in E(\tau) \\ e^\downarrow = \rho(\tau) } }
	T ^\fn _{\fe + e_i \I_{e}}	
+
\sum_{ u \in N(\tau) }
	\fn(u) T^{\fn - e_i \I_u}_\fe
\end{align}
\item One has
\begin{align}\label{eq:ideal:3}
\tau - \cq\hat\Phi\tau.
\end{align}
\end{itemize}
We also define the factor algebras
\[
\KCK := \KCH / \KCI \qquad\text{ and }\qquad
\KCKhat := \KCHhat / \KCIhat.
\]
Here we write $e_i \in \N^{d+1}$ for the $i$-th unit vector.
\end{definition}

With a proof identical to \cite[Prop.~2.12]{Hairer2017}, we obtain the following.
\begin{lemma}\label{lem:KCH:Hopf:algebra}
The ideal $\KCI$ is a Hopf ideal in $\KCH$, so that in particular $\KCK$ is a factor Hopf algebra. Moreover, one has $\CA : \KCIhat \to \KCI$, so that in particular the space $\KCKhat$ forms a left co-module over $\KCK$.
\end{lemma}

Now note that by definition one has $\tau - \cq\hat\Phi \tau \in \KCI$ for any $\tau \in \tilde\KCH$, so that it remains to show that $g_\BPHZ^{\eta,K}$ vanishes on $\KCI$. It follows readily from Lemma \ref{lem:KCH:Hopf:algebra} and the recursive identity for the twisted antipode that one has $\CA : \KCI \to \KCIhat$. It thus remains to show that $g^{\eta,K}$ vanishes identically on $\KCIhat$. This however follows identically to \cite[Property~4]{Hairer2017}.
\end{proof}

\bibliography{bibtex}
\bibliographystyle{siam}

\end{document}